\documentclass{article}
\usepackage{amsmath}
\usepackage{amssymb}
\usepackage{amsthm}
\usepackage{cite}
\usepackage[arrow,curve,matrix]{xy}
\begin{document}
\def\e#1\e{\begin{equation}#1\end{equation}}
\def\ea#1\ea{\begin{align}#1\end{align}}
\def\eq#1{{\rm(\ref{#1})}}
\theoremstyle{plain}
\newtheorem{thm}{Theorem}[section]
\newtheorem{lem}[thm]{Lemma}
\newtheorem{prop}[thm]{Proposition}
\newtheorem{cor}[thm]{Corollary}
\theoremstyle{definition}
\newtheorem{dfn}[thm]{Definition}
\newtheorem{ex}[thm]{Example}
\newtheorem{ass}[thm]{Assumption}
\newtheorem{rem}[thm]{Remark}
\def\na{{\rm na}}
\def\stk{{\rm stk}}
\def\dim{\mathop{\rm dim}\nolimits}
\def\cha{\mathop{\rm char}}
\def\Ker{\mathop{\rm Ker}}
\def\Spec{\mathop{\rm Spec}\nolimits}
\def\Sch{\mathop{\rm Sch}\nolimits}
\def\coh{\mathop{\rm coh}}
\def\Hom{\mathop{\rm Hom}\nolimits}
\def\Iso{\mathop{\rm Iso}\nolimits}
\def\Aut{\mathop{\rm Aut}}
\def\End{\mathop{\rm End}}
\def\GL{\mathop{\rm GL}\nolimits}
\def\Mor{\mathop{\rm Mor}\nolimits}
\def\Stab{\mathop{\rm Stab}\nolimits}
\def\vi{{\rm vi}}
\def\rk{{\rm rk}}
\def\CF{\mathop{\rm CF}\nolimits}
\def\CFi{\mathop{\rm CF}\nolimits^{\rm ind}}
\def\LCF{\mathop{\rm LCF}\nolimits}
\def\SF{\mathop{\rm SF}\nolimits}
\def\SFa{\mathop{\rm SF}\nolimits_{\rm al}}
\def\SFai{\mathop{\rm SF}\nolimits_{\rm al}^{\rm ind}}
\def\uSF{\mathop{\smash{\underline{\rm SF\!}\,}}\nolimits}
\def\oSF{\mathop{\bar{\rm SF}}\nolimits}
\def\oSFa{\mathop{\bar{\rm SF}}\nolimits_{\rm al}}
\def\oSFai{{\ts\bar{\rm SF}{}_{\rm al}^{\rm ind}}}
\def\uoSF{\mathop{\bar{\underline{\rm SF\!}\,}}\nolimits}
\def\LSF{\mathop{\rm LSF}\nolimits}
\def\Ht{{\mathcal H}{}^{\rm to}}
\def\bHt{\bar{\mathcal H}{}^{\rm to}}
\def\Hp{{\mathcal H}{}^{\rm pa}}
\def\bHp{\bar{\mathcal H}{}^{\rm pa}}
\def\Lt{{\mathcal L}{}^{\rm to}}
\def\bLt{\bar{\mathcal L}{}^{\rm to}}
\def\Lp{{\mathcal L}{}^{\rm pa}}
\def\bLp{\bar{\mathcal L}{}^{\rm pa}}
\def\Ext{\mathop{\rm Ext}\nolimits}
\def\id{\mathop{\rm id}\nolimits}
\def\Obj{\mathop{\rm Obj\kern .1em}\nolimits}
\def\Oss{\mathop{\rm Obj\kern .1em}\nolimits_{\rm ss}}
\def\Ost{\mathop{\rm Obj\kern .1em}\nolimits_{\rm st}}
\def\Osi{\mathop{\rm Obj\kern .1em}\nolimits_{\rm si}}
\def\fObj{\mathop{\mathfrak{Obj}\kern .05em}\nolimits}
\def\bdim{{\mathbin{\bf dim\kern .1em}}}
\def\modA{\text{\rm mod-$A$}}
\def\modKQ{\text{\rm mod-$\K Q$}}
\def\modKQI{\text{\rm mod-$\K Q/I$}}
\def\nilKQ{\text{\rm nil-$\K Q$}}
\def\nilKQI{\text{\rm nil-$\K Q/I$}}
\def\nilCQ{\text{\rm nil-$\C Q$}}
\def\bs{\boldsymbol}
\def\ge{\geqslant}
\def\le{\leqslant\nobreak}
\def\pr{{\mathop{\preceq}\nolimits}}
\def\npr{{\mathop{\npreceq}\nolimits}}
\def\tl{\trianglelefteq\nobreak}
\def\ntl{\ntrianglelefteq\nobreak}
\def\ps{\precsim\nobreak}
\def\ls{{\mathop{\lesssim\kern .05em}\nolimits}}
\def\N{{\mathbin{\mathbb N}}}
\def\R{{\mathbin{\mathbb R}}}
\def\C{{\mathbin{\mathbb C}}}
\def\Z{{\mathbin{\mathbb Z}}}
\def\Q{{\mathbin{\mathbb Q}}}
\def\K{{\mathbin{\mathbb K\kern .05em}}}
\def\KP{{\mathbin{\mathbb{KP}}}}
\def\A{{\mathbin{\mathcal A}}}
\def\F{{\mathbin{\mathcal F}}}
\def\G{{\mathbin{\mathcal G}}}
\def\H{{\mathbin{\mathcal H}}}
\def\L{{\mathbin{\mathcal L}}}
\def\M{{\mathcal M}}
\def\cP{{\mathbin{\mathcal P}}}
\def\cQ{{\mathbin{\mathcal Q}}}
\def\cR{{\mathbin{\mathcal R}}}
\def\U{{\mathbin{\mathcal U}}}
\def\fD{{\mathbin{\mathfrak D}}}
\def\fE{{\mathbin{\mathfrak E}}}
\def\fF{{\mathbin{\mathfrak F}}}
\def\fG{{\mathbin{\mathfrak G}}}
\def\fH{{\mathbin{\mathfrak H}}}
\def\fM{{\mathbin{\mathfrak M}}}
\def\fR{{\mathbin{\mathfrak R}}}
\def\fS{{\mathbin{\mathfrak S}}}
\def\fT{{\mathbin{\mathfrak T}}}
\def\sIp{{\smash{\sst(I,\pr)}}}
\def\sJp{{\smash{\sst(J,\pr)}}}
\def\sIt{{\smash{\sst(I,\tl)}}}
\def\sKt{{\smash{\sst(K,\tl)}}}
\def\sJl{{\smash{\sst(J,\ls)}}}
\def\dss{\de_{\rm ss}}
\def\dssb{\de_{\smash{\rm ss}}^{\,\rm b}}
\def\dst{\de_{\rm st}}
\def\dstb{\de_{\smash{\rm st}}^{\,\rm b}}
\def\dsi{\de_{\rm si}}
\def\dsib{\de_{\smash{\rm si}}^{\,\rm b}}
\def\bdss{\bar\de_{\rm ss}}
\def\bdssb{\bar\de_{\smash{\rm ss}}^{\,\rm b}}
\def\bdst{\bar\de_{\rm st}}
\def\bdstb{\bar\de_{\smash{\rm st}}^{\,\rm b}}
\def\bdsi{\bar\de_{\rm si}}
\def\bdsib{\bar\de_{\smash{\rm si}}^{\,\rm b}}
\def\Mss{{\mathcal M}_{\rm ss}}
\def\Mst{{\mathcal M}_{\rm st}}
\def\Msi{{\mathcal M}_{\rm si}}
\def\Mssb{{\mathcal M}_{\rm ss}^{\,\rm b}}
\def\Mstb{{\mathcal M}_{\rm st}^{\,\rm b}}
\def\Msib{{\mathcal M}_{\rm si}^{\,\rm b}}
\def\al{\alpha}
\def\be{\beta}
\def\ga{\gamma}
\def\de{\delta}
\def\bde{\bar\delta}
\def\io{\iota}
\def\ep{\epsilon}
\def\bep{\bar\epsilon}
\def\la{\lambda}
\def\ka{\kappa}
\def\th{\theta}
\def\up{\upsilon}
\def\si{\sigma}
\def\om{\omega}
\def\De{\Delta}
\def\La{\Lambda}
\def\Om{\Omega}
\def\Ga{\Gamma}
\def\Th{\Theta}
\def\Up{\Upsilon}
\def\ts{\textstyle}
\def\sst{\scriptscriptstyle}
\def\sm{\setminus}
\def\bu{\bullet}
\def\op{\oplus}
\def\ot{\otimes}
\def\bigop{\bigoplus}
\def\bigot{\bigotimes}
\def\lt{\ltimes}
\def\iy{\infty}
\def\ra{\rightarrow}
\def\ab{\allowbreak}
\def\longra{\longrightarrow}
\def\t{\times}
\def\ci{\circ}
\def\ti{\tilde}
\def\md#1{\vert #1 \vert}
\def\bmd#1{\big\vert #1 \big\vert}
\title{Configurations in abelian categories. III. \\
Stability conditions and identities}
\author{Dominic Joyce}
\date{}
\maketitle

\baselineskip 11.8pt plus .2pt

\begin{abstract}
This is the third in a series on {\it configurations\/} in
an abelian category $\A$. Given a finite poset $(I,\pr)$, an
$(I,\pr)$-{\it configuration\/} $(\si,\io,\pi)$ is a finite
collection of objects $\si(J)$ and morphisms $\io(J,K)$ or
$\pi(J,K):\si(J)\ra\si(K)$ in $\A$ satisfying some axioms,
where $J,K$ are subsets of $I$. Configurations describe how
an object $X$ in $\A$ decomposes into subobjects.

The first paper defined configurations and studied moduli spaces
of configurations in $\A$, using the theory of Artin stacks. It
showed well-behaved moduli stacks $\fObj_\A,\fM(I,\pr)_\A$ of
objects and configurations in $\A$ exist when $\A$ is the
abelian category $\coh(P)$ of coherent sheaves on a projective
scheme $P$, or $\modKQ$ of representations of a quiver $Q$.
The second studied algebras of {\it constructible functions\/}
and {\it stack functions} on~$\fObj_\A$.

This paper introduces ({\it weak\/}) {\it stability conditions\/}
$(\tau,T,\le)$ on $\A$. We show the moduli spaces
$\Oss^\al,\Osi^\al,\Ost^\al(\tau)$ of $\tau$-semistable,
indecomposable $\tau$-semistable and $\tau$-stable objects in class
$\al$ are {\it constructible sets} in $\fObj_\A$, and some
associated configuration moduli spaces $\Mss,\Msi,\ab\Mst,\ab
\Mssb,\ab\Msib,\ab\Mstb(I,\pr,\ka, \tau)_\A$ constructible in
$\fM(I,\pr)_\A$, so their characteristic functions
$\dss^\al,\dsi^\al,\dst^\al(\tau)$ and $\dss,\ldots,\dstb
(I,\pr,\ka,\tau)$ are constructible.

We prove many identities relating these constructible functions, and
their stack function analogues, under pushforwards. We introduce
interesting algebras $\Hp_\tau,\Ht_\tau,\bHp_\tau,\bHt_\tau$ of
constructible and stack functions, and study their structure. In the
fourth paper we show $\Hp_\tau,\ldots,\bHt_\tau$ are independent of
$(\tau,T,\le)$, and construct {\it invariants} of~$\A,(\tau,T,\le)$.
\end{abstract}

\section{Introduction}
\label{as1}

This is the third in a series of papers \cite{Joyc3,Joyc4,Joyc5} on
{\it configurations}. Given an abelian category $\A$ and a finite
partially ordered set (poset) $(I,\pr)$, we define an $(I,\pr)$-{\it
configuration} $(\si,\io,\pi)$ in $\A$ to be a collection of objects
$\si(J)$ and morphisms $\io(J,K)$ or $\pi(J,K):\si(J)\ra\si(K)$ in
$\A$ satisfying certain axioms, for~$J,K\subseteq I$.

The first paper \cite{Joyc3} defined configurations, developed their
basic properties, and studied moduli spaces of configurations in
$\A$, using the theory of Artin stacks. It proved well-behaved
moduli stacks $\fObj_\A,\fM(I,\pr)_\A$ of objects and configurations
exist when $\A$ is the abelian category $\coh(P)$ of coherent
sheaves on a projective $\K$-scheme $P$, or $\modKQ$ of
representations of a quiver $Q$. The second \cite{Joyc4} defined and
studied infinite-dimensional algebras of {\it constructible
functions} and {\it stack functions} on $\fObj_\A$, motivated by
{\it Ringel--Hall algebras}.

Configurations are a tool for describing how an object $X$ in $\A$
decomposes into subobjects. They are especially useful for studying
{\it stability conditions} on $\A$, which are the subject of this
paper. Given a stability condition $(\tau,T,\le)$ on $\A$, objects
$X$ in $\A$ are called $\tau$-{\it semistable}, $\tau$-{\it stable}
or $\tau$-{\it unstable} according to whether subobjects $S\subset
X$ with $S\ne 0,X$ have $\tau([S])\le\tau([X])$, $\tau([S])<\tau
([X])$, or $\tau([S])>\tau([X])$. Examples of stability conditions
include slope functions, and Gieseker stability of coherent sheaves.

We also define {\it weak stability conditions}, which include
$\mu$-stability and purity for coherent sheaves. When
$(\tau,T,\le)$ is a weak stability condition each $X\in\A$
has a unique {\it Harder--Narasimhan filtration} by subobjects
$0\!=\!A_0\!\subset\!\cdots\!\subset\!A_n\!=\!X$ whose factors
$S_k=A_k/A_{k-1}$ are $\tau$-semistable with $\tau([S_1])>\cdots>
\tau([S_n])$. If $(\tau,T,\le)$ is also a stability condition each
$\tau$-semistable $X$ has a (nonunique) filtration with (unique)
$\tau$-stable factors $S_k$ with $\tau([S_k])=\tau([X])$. Thus,
$\tau$-stability is well-behaved for stability conditions but
badly behaved for weak stability conditions, though
$\tau$-semistability is well-behaved for both.

We form moduli spaces $\Oss^\al,\Osi^\al,\Ost^\al(\tau)$ of
$\tau$-semistable, $\tau$-semistable-indecomposable and
$\tau$-stable objects in class $\al$ in $K(\A)$, and moduli spaces
$\Mss,\Msi,\Mst,\Mssb,\Msib,\Mstb(I,\pr,\ka,\tau)_\A$ of
$(I,\pr)$-configurations $(\si,\io,\pi)$ in which the smallest
objects $\si(\{i\})$ for $i\in I$ lie in $\Oss^{\ka(i)}\!,
\Osi^{\ka(i)}\!,\Ost^{\ka(i)}(\tau)$, and $(\si,\io,\pi)$ is {\it
best\/} for ${\mathcal M}_*^{\rm b}(\cdots)_\A$. It is a central,
and unconventional, feature of our approach that we regard these not
as spaces in their own right, but as {\it constructible sets\/} in
the stacks $\fObj_\A,\fM(I,\pr)_\A$, so their characteristic
functions $\dss^\al,\dsi^\al,\dst^\al(\tau)$ and
$\dss,\ldots,\dstb(I,\pr,\ka,\tau)$ are {\it constructible
functions}.

This has a number of ramifications. Firstly, our approach is helpful
for comparing moduli spaces, and especially for understanding how
$\Oss^\al(\tau)$ changes when we vary $(\tau,T,\le)$, as we are not
comparing two different varieties, but two subsets of the same stack
$\fObj_\A$. Secondly, $\Oss^\al(\tau)$ is a set of isomorphism
classes, not of S-equivalence classes. This is better for studying
the family of ways a $\tau$-semistable $X$ may be broken into
$\tau$-stable factors. But it means $\Oss^\al(\tau)$ is {\it not a
well-behaved topological space}, as it may not be Hausdorff, for
instance. Because of this, in \cite{Joyc5} we focus on `motivic'
invariants of constructible sets such as Euler characteristics and
virtual Poincar\'e polynomials.

We begin in \S\ref{as2} with background on abelian categories,
constructible sets and functions, and {\it stack functions} on Artin
$\K$-stacks, following \cite{Joyc1,Joyc2}. Stack functions are a
universal generalization of constructible functions, containing more
information. Section \ref{as3} reviews the previous papers
\cite{Joyc3,Joyc4}, and \S\ref{as4} defines ({\it weak\/}) {\it
stability conditions} $(\tau,T,\le)$ on $\A$. If $(\tau, T,\le)$ is
{\it permissible} $\Oss^\al,\Osi^\al,\Ost^\al(\tau)$ and
$\Mss,\ldots,\Mstb(I,\pr,\ka,\tau)_\A$ are {\it constructible sets}.
We give examples of permissible (weak) stability conditions on
$\A=\modKQ$ and~$\A=\coh(P)$.

Sections \ref{as5} and \ref{as6} prove identities relating the six
families of constructible functions $\dss,\dsi,\dst,
\dssb,\dsib,\dstb(I,\pr,\ka,\tau)$. These depend on theorems on the
Euler characteristics of parts of moduli spaces, and encode facts
about the family of ways of decomposing a $\tau$-semistable object
into $\tau$-stable factors, and so on. One conclusion is that each
of the six families determines the other five.

Section \ref{as7} studies the {\it algebras of constructible
functions} $\Hp_\tau,\Ht_\tau$ on $\fObj_\A$ generated by
$\CF^\stk(\bs\si(I))\dss(I,\pr,\ka,\tau),\dss^\al(\tau)$
respectively, for all $(I,\pr,\ka),\al$. Defining Lie algebras
$\Lp_\tau,\Lt_\tau$ to be the intersections of $\Hp_\tau,\Ht_\tau$
with the Lie subalgebra $\CFi(\fObj_\A)\subset\CF(\fObj_\A)$
supported on indecomposables in $\A$, we construct generators of
$\Hp_\tau,\Ht_\tau$ lying in $\Lp_\tau,\Lt_\tau$, and so show
$\Hp_\tau,\Ht_\tau$ are the {\it universal enveloping algebras}
of~$\Lp_\tau,\Lt_\tau$.

Finally, \S\ref{as8} generalizes the results of
\S\ref{as5}--\S\ref{as7} from constructible functions to the {\it
stack functions} of \cite{Joyc2}, giving {\it stack\/} ({\it Lie})
{\it algebras} $\bHp_\tau,\bHt_\tau,\bLp_\tau,\bLt_\tau$. The sequel
\cite{Joyc5} will show the (Lie) algebras $\Hp_\tau,\ldots,
\bLt_\tau$ are {\it independent of\/} $(\tau,T,\le)$, so that many
of our identities here and in \cite{Joyc5} can be regarded as {\it
change of basis formulae\/} in $\Hp_\tau,\ldots,\bLt_\tau$. It also
discusses systems of invariants of $\A,(\tau,T,\le)$ `counting'
$\tau$-semistable objects and configurations, and their identities
and transformation laws. These can often be interpreted using
morphisms from $\bHp_\tau,\ldots,\bLt_\tau$ to an explicit (Lie)
algebra, as in~\cite[\S 6]{Joyc4}.

A subsequent paper \cite{Joyc6} explains how to encode some of the
invariants of \cite{Joyc5} into {\it holomorphic generating
functions\/} on the complex manifold of stability conditions. These
satisfy an interesting p.d.e., that can be interpreted as the
flatness of a connection. The material of \S\ref{as7} will be
important in~\cite{Joyc6}.
\medskip

\noindent{\it Acknowledgements.} I would like to thank
Tom Bridgeland for many inspiring conversations and for
being interested, Frances Kirwan, Andrew Kresch and Burt
Totaro for help with moduli spaces and stacks, and
Alastair King and Richard Thomas for useful conversations.
I was supported by an EPSRC Advanced Research Fellowship
whilst writing this paper.

\section{Background material}
\label{as2}

We begin with some background material on abelian categories
in \S\ref{as21}, and Artin stacks, constructible functions
and stack functions in~\S\ref{as22}--\S\ref{as24}.

\subsection{Abelian categories}
\label{as21}

Here is the definition of abelian category, taken
from~\cite[\S II.5]{GeMa}.

\begin{dfn} A category $\A$ is called {\it abelian} if
\begin{itemize}
\setlength{\itemsep}{0pt}
\setlength{\parsep}{0pt}
\item[(i)] $\Hom(X,Y)$ is an abelian group for all $X,Y\in\A$,
and composition of morphisms is biadditive.
\item[(ii)] There exists a zero object $0\in\A$ such
that~$\Hom(0,0)=0$.
\item[(iii)] For any $X,Y\in\A$ there exists $Z\in\A$ and
morphisms $\io_X:X\ra Z$, $\io_Y:Y\ra Z$, $\pi_X:Z\ra X$,
$\pi_Y:Z\ra Y$ with $\pi_X\ci\io_X=\id_X$, $\pi_Y\ci\io_Y=
\id_Y$, $\io_X\ci\pi_X+\io_Y\ci\pi_Y=\id_Z$ and
$\pi_X\ci\io_Y=\pi_Y\ci\io_X=0$. We write $Z=X\op Y$,
the {\it direct sum} of $X$ and~$Y$.
\item[(iv)] For any morphism $f:X\ra Y$ there is a sequence
$\smash{K{\buildrel k\over\ra}X{\buildrel i\over\ra}I{\buildrel
j\over\ra}Y{\buildrel c\over\ra}C}$ in $\A$ such that $j\ci i=f$,
and $K$ is the kernel of $f$, and $C$ the cokernel of $f$, and
$I$ is both the cokernel of $k$ and the kernel of~$c$.
\end{itemize}
\label{as2def1}
\end{dfn}

In an abelian category we can define {\it exact sequences} as in
\cite[\S II.6]{GeMa}. A short exact sequence $0\ra X\smash{{\buildrel
f\over\ra}Y{\buildrel g\over\ra}}Z\ra 0$ in $\A$ is called
{\it split\/} if there exists a compatible isomorphism $h:X\op Z\ra Y$.
The {\it Grothendieck group} $K_0(\A)$ of $\A$ is the abelian group
generated by $\Obj(\A)$, with a relation $[Y]=[X]+[Z]$ for each short
exact sequence $0\!\ra\!X\!\ra\!Y\!\ra\!Z\!\ra\!0$ in $\A$. Throughout
the paper $K(\A)$ will mean {\it the quotient of\/ $K_0(\A)$ by some
fixed subgroup}. {\it Subobjects} of objects in $\A$ are analogous to
subgroups of an abelian group.

\begin{dfn} Let $\A$ be an abelian category and $X\in\A$.
Two injective morphisms $i:S\ra X$, $i':S'\ra X$ are
called equivalent if there exists an isomorphism
$h:S\ra S'$ with $i=i'\ci h$. A {\it subobject\/} of $X$
is an equivalence class of injective $i:S\ra X$. Usually
we refer to $S$ as the subobject, taking $i$ and the
equivalence class to be implicitly given, and write
$S\subset X$ to mean $S$ is a subobject of $X$. If
$S,T\subset X$ are represented by $i:S\ra X$ and
$j:T\ra X$, we write $S\subset T\subset X$ if there
exists $a:S\ra T$ with~$i=j\ci a$.

We call $\A$ {\it artinian} if for all $X\in\A$, all
descending chains of subobjects $\cdots\!\subset\!A_2\!
\subset\!A_1\!\subset\!X$ stabilize, that is, $A_{n+1}=A_n$
for all $n\gg 0$. We call $\A$ {\it noetherian} if all
ascending chains of subobjects $A_1\!\subset\!A_2\!
\subset\!\cdots\!\subset\!X$ stabilize.
\label{as2def2}
\end{dfn}

\subsection{Introduction to Artin $\K$-stacks}
\label{as22}

Fix an algebraically closed field $\K$ throughout. There
are four main classes of `spaces' over $\K$ used in
algebraic geometry, in increasing order of generality:
\begin{equation*}
\text{$\K$-varieties}\subset
\text{$\K$-schemes}\subset
\text{algebraic $\K$-spaces}\subset
\text{algebraic $\K$-stacks}.
\end{equation*}

{\it Algebraic stacks} (also known as Artin stacks) were
introduced by Artin, generalizing {\it Deligne--Mumford stacks}.
For a good introduction to algebraic stacks see G\'omez
\cite{Gome}, and for a thorough treatment see Laumon and
Moret-Bailly \cite{LaMo}. We make the convention that all
algebraic $\K$-stacks in this paper are {\it locally of
finite type}, and $\K$-substacks are {\it locally closed}.

Algebraic $\K$-stacks form a 2-{\it category}. That is, we have {\it
objects} which are $\K$-stacks $\fF,\fG$, and also two kinds of
morphisms, 1-{\it morphisms} $\phi,\psi:\fF\ra\fG$ between
$\K$-stacks, and 2-{\it morphisms} $A:\phi\ra\psi$ between
1-morphisms. An analogy to keep in mind is a 2-category of
categories, where objects are categories, 1-morphisms are functors
between the categories, and 2-morphisms are isomorphisms (natural
transformations) between functors.

We define the set of $\K$-{\it points} of a stack.

\begin{dfn} Let $\fF$ be a $\K$-stack. Write $\fF(\K)$ for the set of
2-isomorphism classes $[x]$ of 1-morphisms $x:\Spec\K\ra\fF$.
Elements of $\fF(\K)$ are called $\K$-{\it points}, or {\it
geometric points}, of $\fF$. If $\phi:\fF\ra\fG$ is a 1-morphism
then composition with $\phi$ induces a map of
sets~$\phi_*:\fF(\K)\ra\fG(\K)$.

For a 1-morphism $x:\Spec\K\ra\fF$, the {\it stabilizer group}
$\Iso_\K(x)$ is the group of 2-morphisms $x\ra x$. When $\fF$ is an
algebraic $\K$-stack, $\Iso_\K(x)$ is an {\it algebraic $\K$-group}.
We say that $\fF$ {\it has affine geometric stabilizers} if
$\Iso_\K(x)$ is an affine algebraic $\K$-group for all 1-morphisms
$x:\Spec\K\ra\fF$.

As an algebraic $\K$-group up to isomorphism, $\Iso_\K(x)$
depends only on the isomorphism class $[x]\in\fF(\K)$ of $x$
in $\Hom(\Spec\K,\fF)$. If $\phi:\fF\ra\fG$ is a 1-morphism,
composition induces a morphism of algebraic $\K$-groups
$\phi_*:\Iso_\K([x])\ra\Iso_\K\bigr(\phi_*([x])\bigr)$,
for~$[x]\in\fF(\K)$.
\label{as2def3}
\end{dfn}

One important difference in working with 2-categories rather than
ordinary categories is that in diagram-chasing one only requires
1-morphisms to be 2-{\it isomorphic} rather than {\it equal}. The
simplest kind of {\it commutative diagram} is:
\begin{equation*}
\xymatrix@R=6pt{
& \fG \ar@{=>}[d]^{\,F} \ar[dr]^\psi \\
\fF \ar[ur]^\phi \ar[rr]_\chi && \fH, }
\end{equation*}
by which we mean that $\fF,\fG,\fH$ are $\K$-stacks,
$\phi,\psi,\chi$ are 1-morphisms, and $F:\psi\ci\phi\ra\chi$ is a
2-isomorphism. Usually we omit $F$, and mean
that~$\psi\ci\phi\cong\chi$.

\begin{dfn} Let $\phi:\fF\ra\fH$, $\psi:\fG\ra\fH$ be 1-morphisms
of $\K$-stacks. Then one can define the {\it fibre product stack\/}
$\fF\t_{\phi,\fH,\psi}\fG$, or $\fF\t_\fH\fG$ for short, with
1-morphisms $\pi_\fF,\pi_\fG$ fitting into a commutative diagram:
\e
\begin{gathered}
\xymatrix@R=-4pt{
& \fF \ar[dr]^\phi \ar@{=>}[dd] \\
\fF\t_\fH\fG
\ar[dr]_{\pi_\fG} \ar[ur]^{\pi_\fF} && \fH.\\
& \fG \ar[ur]_\psi \\
}
\end{gathered}
\label{as2eq1}
\e
A commutative diagram
\begin{equation*}
\xymatrix@R=-4pt{
& \fF \ar[dr]^\phi \ar@{=>}[dd] \\
\fE
\ar[dr]_\eta \ar[ur]^\th && \fH\\
& \fG \ar[ur]_\psi \\
}
\end{equation*}
is a {\it Cartesian square} if it is isomorphic to \eq{as2eq1}, so
there is a 1-isomorphism $\fE\cong\fF\t_\fH\fG$. Cartesian squares
may also be characterized by a universal property.
\label{as2def4}
\end{dfn}

\subsection{Constructible functions on stacks}
\label{as23}

Next we discuss {\it constructible functions} on $\K$-stacks,
following \cite{Joyc1}. For this section we need $\K$ to have {\it
characteristic zero}.

\begin{dfn} Let $\fF$ be an algebraic $\K$-stack. We call
$C\subseteq\fF(\K)$ {\it constructible} if $C=\bigcup_{i\in I}
\fF_i(\K)$, where $\{\fF_i:i\in I\}$ is a finite collection of
finite type algebraic $\K$-substacks $\fF_i$ of $\fF$. We call
$S\subseteq\fF(\K)$ {\it locally constructible} if $S\cap C$
is constructible for all constructible~$C\subseteq\fF(\K)$.

A function $f:\fF(\K)\ra\Q$ is called {\it constructible} if
$f(\fF(\K))$ is finite and $f^{-1}(c)$ is a constructible set
in $\fF(\K)$ for each $c\in f(\fF(\K))\sm\{0\}$. A function
$f:\fF(\K)\ra\Q$ is called {\it locally constructible} if
$f\cdot\de_C$ is constructible for all constructible
$C\subseteq\fF(\K)$, where $\de_C$ is the characteristic
function of $C$. Write $\CF(\fF)$ and $\LCF(\fF)$ for the
$\Q$-vector spaces of $\Q$-valued constructible and
locally constructible functions on~$\fF$.
\label{as2def5}
\end{dfn}

Here \cite[\S 4]{Joyc1} are some important properties of
constructible sets.

\begin{prop} Let\/ $\fF,\fG$ be algebraic $\K$-stacks with
affine geometric stabilizers, $\phi:\fF\!\ra\!\fG$ a $1$-morphism,
and\/ $A,B\subseteq\fF(\K)$ constructible. Then $A\cup B,A\cap B$
and\/ $A\sm B$ are constructible in $\fF(\K)$, and\/
$\phi_*(A)$ is constructible in~$\fG(\K)$.
\label{as2prop}
\end{prop}

Following \cite[Def.s~4.8, 5.1 \& 5.5]{Joyc1} we define
{\it pushforwards} and {\it pullbacks} of constructible
functions along 1-morphisms.

\begin{dfn} In \cite[\S 3.3]{Joyc1} we define the {\it
Euler characteristic} $\chi(\cdots)$ of constructible subsets in
$\K$-schemes. In \S\ref{as5} we use the fact
\cite[Th.~3.10(vi)]{Joyc1} that
\e
\chi(\K^m)=1\quad\text{and}\quad \chi(\KP^m)=m+1
\quad\text{for all $m\ge 0$.}
\label{as2eq2}
\e
Let $\fF$ be an algebraic $\K$-stack with affine geometric
stabilizers, and $C\subseteq\fF(\K)$ a constructible subset.
Then \cite[Def.~4.8]{Joyc1} defines the {\it na\"\i ve Euler
characteristic} $\chi^\na(C)$ of $C$. It is called {\it
na\"\i ve} as it takes no account of stabilizer groups. For
$f\in\CF(\fF)$, define $\chi^\na(\fF,f)$ in $\Q$ by
$\chi^\na(\fF,f)=\sum_{c\in f(\fF(\K))\sm\{0\}}c\,\chi^\na
\bigl(f^{-1}(c)\bigr)$.

Let $\fF,\fG$ be algebraic $\K$-stacks with affine geometric
stabilizers, and $\phi:\fF\ra\fG$ a representable 1-morphism.
Then for any $x\in\fF(\K)$ we have an injective morphism
$\phi_*:\Iso_\K(x)\ra\Iso_\K\bigl(\phi_*(x)\bigr)$ of affine
algebraic $\K$-groups. The image $\phi_*\bigl(\Iso_\K(x)\bigr)$
is an affine algebraic $\K$-group closed in $\Iso_\K\bigl(
\phi_*(x)\bigr)$, so the quotient $\Iso_\K\bigl(\phi_*(x)\bigr)
/\phi_*\bigl(\Iso_\K(x)\bigr)$ exists as a quasiprojective
$\K$-variety. Define a function $m_\phi:\fF(\K)\ra\Z$ by
$m_\phi(x)=\chi\bigl(\Iso_\K(\phi_*(x))/\phi_*(\Iso_\K(x))
\bigr)$ for~$x\in\fF(\K)$.

For $f\in\CF(\fF)$, define $\CF^\stk(\phi)f:\fG(\K)\ra\Q$ by
\begin{equation*}
\CF^\stk(\phi)f(y)=\chi^\na\bigl(\fF,m_\phi\cdot f\cdot
\de_{\phi_*^{-1}(y)}\bigr) \quad\text{for $y\in\fG(\K)$,}
\end{equation*}
where $\de_{\smash{\phi_*^{-1}(y)}}$ is the characteristic function
of $\phi_*^{-1}(\{y\})\subseteq\fG(\K)$ on $\fG(\K)$. Then
$\CF^\stk(\phi):\CF(\fF)\ra\CF(\fG)$ is a $\Q$-linear map
called the {\it stack pushforward}.

Let $\th:\fF\ra\fG$ be a finite type 1-morphism. If $C\subseteq
\fG(\K)$ is constructible then so is $\th_*^{-1}(C)\subseteq\fF(\K)$.
It follows that if $f\in\CF(\fG)$ then $f\ci\th_*$ lies in $\CF(\fF)$.
Define the {\it pullback\/} $\th^*:\CF(\fG)\ra\CF(\fF)$ by $\th^*(f)=
f\ci\th_*$. It is a linear map.
\label{as2def6}
\end{dfn}

Here \cite[Th.s~5.4 \& 5.6 \& Def.~5.5]{Joyc1} are some properties
of these.

\begin{thm} Let\/ $\fE,\fF,\fG,\fH$ be algebraic $\K$-stacks with
affine geometric stabilizers, and\/ $\be:\fF\ra\fG$, $\ga:\fG\ra\fH$
be $1$-morphisms. Then
\ea
\CF^\stk(\ga\ci\be)&=\CF^\stk(\ga)\ci\CF^\stk(\be):\CF(\fF)\ra\CF(\fH),
\label{as2eq3}\\
(\ga\ci\be)^*&=\be^*\ci\ga^*:\CF(\fH)\ra\CF(\fF),
\label{as2eq4}
\ea
supposing $\be,\ga$ representable in \eq{as2eq3}, and of finite type
in \eq{as2eq4}. If
\e
\begin{gathered}
\xymatrix{
\fE \ar[r]_\eta \ar[d]^\th & \fG \ar[d]_\psi \\
\fF \ar[r]^\phi & \fH
}
\end{gathered}
\quad
\begin{gathered}
\text{is a Cartesian square with}\\
\text{$\eta,\phi$ representable and}\\
\text{$\th,\psi$ of finite type, then}\\
\text{the following commutes:}
\end{gathered}
\quad
\begin{gathered}
\xymatrix@C=35pt{
\CF(\fE) \ar[r]_{\CF^\stk(\eta)} & \CF(\fG) \\
\CF(\fF) \ar[r]^{\CF^\stk(\phi)} \ar[u]_{\th^*}
& \CF(\fH). \ar[u]^{\psi^*}
}
\end{gathered}
\label{as2eq5}
\e
\label{as2thm1}
\end{thm}

As discussed in \cite[\S 3.3]{Joyc1} for the $\K$-scheme case,
equation \eq{as2eq3} is {\it false} for algebraically closed fields
$\K$ of characteristic $p>0$. The definitions and results above all
have analogues for locally constructible functions,~\cite[\S
5.3]{Joyc1}.

\subsection{Stack functions}
\label{as24}

{\it Stack functions} are a universal generalization of
constructible functions introduced in \cite{Joyc2}. Here
\cite[Def.~3.1]{Joyc2} is the basic definition. Throughout $\K$ is
algebraically closed of arbitrary characteristic, except when we
specify~$\cha\K=0$.

\begin{dfn} Let $\fF$ be an algebraic $\K$-stack with affine
geometric stabilizers. Consider pairs $(\fR,\rho)$, where $\fR$ is a
finite type algebraic $\K$-stack with affine geometric stabilizers
and $\rho:\fR\ra\fF$ is a representable 1-morphism. We call two
pairs $(\fR,\rho)$, $(\fR',\rho')$ {\it equivalent\/} if there
exists a 1-isomorphism $\io:\fR\ra \fR'$ such that $\rho' \ci\io$
and $\rho$ are 2-isomorphic 1-morphisms $\fR\ra\fF$. Write
$[(\fR,\rho)]$ for the equivalence class of $(\fR,\rho)$. If
$(\fR,\rho)$ is such a pair and $\fS$ is a closed $\K$-substack of
$\fR$ then $(\fS,\rho\vert_\fS)$,
$(\fR\sm\fS,\rho\vert_{\fR\sm\fS})$ are pairs of the same kind.
Define $\SF(\fF)$ to be the $\Q$-vector space generated by
equivalence classes $[(\fR,\rho)]$ as above, with for each closed
$\K$-substack $\fS$ of $\fR$ a relation
\begin{equation*}
[(\fR,\rho)]=[(\fS,\rho\vert_\fS)]+[(\fR\sm\fS,\rho\vert_{\fR\sm\fS})].
\end{equation*}
\label{as2def7}
\end{dfn}

In \cite[Def.~3.2]{Joyc2} we relate $\CF(\fF)$ and~$\SF(\fF)$.

\begin{dfn} Let $\fF$ be an algebraic $\K$-stack with affine
geometric stabilizers and $C\subseteq\fF(\K)$ be constructible. Then
$C=\coprod_{i=1}^n\fR_i(\K)$, for $\fR_1,\ldots,\fR_n$ finite type
$\K$-substacks of $\fF$. Let $\rho_i:\fR_i\ra\fF$ be the inclusion
1-morphism. Then $[(\fR_i,\rho_i)]\in\SF(\fF)$. Define
$\bde_C=\ts\sum_{i=1}^n[(\fR_i,\rho_i)]\in\SF(\fF)$. We think of
this stack function as the analogue of the characteristic function
$\de_C\in\CF(\fF)$ of $C$. Define a $\Q$-linear map
$\io_\fF:\CF(\fF)\ra\SF(\fF)$ by $\io_\fF(f)=\ts\sum_{0\ne c\in
f(\fF(\K))}c\cdot\bde_{f^{-1}(c)}$. For $\K$ of characteristic zero,
define a $\Q$-linear map $\pi_\fF^\stk:\SF(\fF)\ra\CF(\fF)$ by
\begin{equation*}
\pi_\fF^\stk\bigl(\ts\sum_{i=1}^nc_i[(\fR_i,\rho_i)]\bigr)=
\ts\sum_{i=1}^nc_i\CF^\stk(\rho_i)1_{\fR_i},
\end{equation*}
where $1_{\fR_i}$ is the function 1 in $\CF(\fR_i)$. Then
\cite[Prop.~3.3]{Joyc2} shows $\pi_\fF^\stk\ci\io_\fF$ is the
identity on $\CF(\fF)$. Thus, $\io_\fF$ is injective and
$\pi_\fF^\stk$ is surjective. In general $\io_\fF$ is far from being
surjective, and $\SF(\fF)$ is much larger than~$\CF(\fF)$.
\label{as2def8}
\end{dfn}

In \cite[Def.~3.4]{Joyc2} we define {\it pushforwards}, {\it
pullbacks} and {\it tensor products}.

\begin{dfn} Let $\phi:\fF\!\ra\!\fG$ be a 1-morphism of algebraic
$\K$-stacks with affine geometric stabilizers. For $\phi$ representable,
define the {\it pushforward\/} $\phi_*:\SF(\fF)\!\ra\!\SF(\fG)$ by
$\phi_*:\sum_{i=1}^nc_i[(\fR_i,\rho_i)]\longmapsto\sum_{i=1}^nc_i
[(\fR_i,\phi\ci\rho_i)]$. For $\phi$ of finite type, define the
{\it pullback\/} $\phi^*:\SF(\fG)\!\ra\!\SF(\fF)$~by
\begin{equation*}
\phi^*:\ts\sum_{i=1}^nc_i[(\fR_i,\rho_i)]\longmapsto
\ts\sum_{i=1}^nc_i[(\fR_i\t_{\rho_i,\fG,\phi}\fF,\pi_\fF)].
\end{equation*}
The {\it tensor product\/} $\ot:\SF(\fF)\t\SF(\fG)\ra\SF(\fF\t\fG)$
is
\begin{equation*}
\bigl(\ts\sum_{i=1}^mc_i[(\fR_i,\rho_i)]\bigr)\!\ot\!
\bigl(\ts\sum_{j=1}^nd_j[(\fS_j,\si_j)]\bigr)\!=\!\ts
\sum_{i,j}c_id_j[(\fR_i\!\t\!\fS_j,\rho_i\!\t\!\si_j)].
\end{equation*}
\label{as2def9}
\end{dfn}

Here \cite[Th.~3.5]{Joyc2} is the analogue of Theorem~\ref{as2thm1}.

\begin{thm} Let\/ $\fE,\fF,\fG,\fH$ be algebraic $\K$-stacks with
affine geometric stabilizers, and\/ $\be:\fF\ra\fG$, $\ga:\fG\ra\fH$
be $1$-morphisms. Then
\begin{equation*}
(\ga\ci\be)_*=\ga_*\ci\be_*:\SF(\fF)\ra\SF(\fH),\quad
(\ga\ci\be)^*=\be^*\ci\ga^*:\SF(\fH)\ra\SF(\fF),
\end{equation*}
for $\be,\ga$ representable in the first equation, and of
finite type in the second. If
\begin{equation*}
\begin{gathered}
\xymatrix{
\fE \ar[r]_\eta \ar[d]^{\,\th} & \fG \ar[d]_{\psi\,} \\
\fF \ar[r]^\phi & \fH
}
\end{gathered}
\quad
\begin{gathered}
\text{is a Cartesian square with}\\
\text{$\eta,\phi$ representable and}\\
\text{$\th,\psi$ of finite type, then}\\
\text{the following commutes:}
\end{gathered}
\quad
\begin{gathered}
\xymatrix@C=35pt{
\SF(\fE) \ar[r]_{\eta_*} & \SF(\fG) \\
\SF(\fF) \ar[r]^{\phi_*} \ar[u]_{\,\th^*}
& \SF(\fH). \ar[u]^{\psi^*\,}
}
\end{gathered}
\end{equation*}
\label{as2thm2}
\end{thm}

In \cite[Prop.~3.7 \& Th.~3.8]{Joyc2} we relate pushforwards and
pullbacks of stack and constructible functions
using~$\io_\fF,\pi_\fF^\stk$.

\begin{thm} Let\/ $\K$ have characteristic zero, $\fF,\fG$ be
algebraic $\K$-stacks with affine geometric stabilizers, and\/
$\phi:\fF\ra\fG$ be a $1$-morphism. Then
\begin{itemize}
\setlength{\itemsep}{0pt}
\setlength{\parsep}{0pt}
\item[{\rm(a)}] $\phi^*\!\ci\!\io_\fG\!=\!\io_\fF\!\ci\!\phi^*:
\CF(\fG)\!\ra\!\SF(\fF)$ if\/ $\phi$ is of finite type;
\item[{\rm(b)}] $\pi^\stk_\fG\ci\phi_*=\CF^\stk(\phi)\ci\pi_\fF^\stk:
\SF(\fF)\ra\CF(\fG)$ if\/ $\phi$ is representable; and
\item[{\rm(c)}] $\pi^\stk_\fF\ci\phi^*=\phi^*\ci\pi_\fG^\stk:
\SF(\fG)\ra\CF(\fF)$ if\/ $\phi$ is of finite type.
\end{itemize}
\label{as2thm3}
\end{thm}

In \cite[\S 5.2]{Joyc2} we define {\it projections}
$\Pi^\vi_n:\SF(\fF)\ra\SF(\fF)$ which project to stack functions
whose stabilizer groups have `virtual rank'~$n$.

In \cite[\S 3]{Joyc2} we define {\it local stack functions}
$\LSF(\fF)$, the analogue of locally constructible functions.
Analogues of Definitions \ref{as2def8}--\ref{as2def9} and Theorems
\ref{as2thm2}--\ref{as2thm3} hold for $\LSF(\fF)$, with differences
in which 1-morphisms are required to be of finite type. We also
study enlarged versions $\underline{\rm SF\!}\,(\fF),\underline{\rm
LSF\!}\,(\fF)$ of $\SF(\fF),\LSF(\fF)$ in which the 1-morphisms
$\rho$ of Definition \ref{as2def7} are not supposed representable.

In \cite[\S 4--\S 6]{Joyc2} we define other classes of stack
functions $\uSF,\uoSF,\oSF(\fF,\Up,\La),\ab \uoSF,\ab
\oSF(\fF,\Up,\La^\ci),\uoSF,\oSF(\fF,\Th,\Om)$ `twisted' by a {\it
motivic invariant\/} $\Up$ or $\Th$ of $\K$-varieties, taking values
in a $\Q$-algebra $\La,\La^\ci$ or $\Om$; the basic facts are
explained in \cite[\S 2.4--\S 2.5]{Joyc4}. All the above material on
$\SF(\cdots)$ applies to these spaces, except that $\pi_\fF^\stk,
\Pi^\vi_n$ are not always defined. For the purposes of this paper
the differences between these spaces are unimportant, so we shall
not explain them.

\section{Background on configurations from \cite{Joyc3,Joyc4}}
\label{as3}

We now recall in \S\ref{as31} and \S\ref{as32} the main definitions
and results from \cite{Joyc3} on $(I,\pr)$-configurations and their
moduli stacks that we will need later, and in \S\ref{as33} some
facts about algebras of constructible and stack functions
from~\cite{Joyc4}.

\subsection{Basic definitions}
\label{as31}

Here is some notation for {\it finite posets}, taken from
\cite[Def.s~3.2, 4.1 \& 6.1]{Joyc3}.

\begin{dfn} A {\it finite partially ordered set\/} or {\it
finite poset\/} $(I,\pr)$ is a finite set $I$ with a partial
order $I$. Define $J\subseteq I$ to be an {\it f-set\/} if
$i\in I$ and $h,j\in J$ and $h\pr i\pr j$ implies $i\in J$.
Define $\F_\sIp$ to be the set of f-sets of $I$. Define
$\G_\sIp$ to be the subset of $(J,K)\in\F_\sIp\t\F_\sIp$
such that $J\subseteq K$, and if $j\in J$ and $k\in K$
with $k\pr j$, then $k\in J$. Define $\H_\sIp$ to be the
subset of $(J,K)\in\F_\sIp\t\F_\sIp$ such that
$K\subseteq J$, and if $j\in J$ and $k\in K$ with
$k\pr j$, then~$j\in K$.

Let $I$ be a finite set and $\pr,\tl$ partial orders on $I$
such that if $i\pr j$ then $i\tl j$ for $i,j\in I$. Then we
say that $\tl$ {\it dominates} $\pr$. Let $s$ be the number
of pairs $(i,j)\in I\t I$ with $i\tl j$ but $i\npr j$. Then
we say that $\tl$ {\it dominates $\pr$ by $s$ steps}.

A partial order $\tl$ on $I$ is called a {\it total order}
if $i\tl j$ or $j\tl i$ for all $i,j\in I$. Then $(I,\tl)$
is canonically isomorphic to $(\{1,\ldots,n\},\le)$ for
$n=\md{I}$. Every partial order $\pr$ on $I$ is dominated
by a total order~$\tl$.
\label{as3def1}
\end{dfn}

We define $(I,\pr)$-{\it configurations},~\cite[Def.~4.1]{Joyc3}.

\begin{dfn} Let $(I,\pr)$ be a finite poset, and use the
notation of Definition \ref{as3def1}. Define an $(I,\pr)$-{\it
configuration} $(\si,\io,\pi)$ in an abelian category $\A$ to be
maps $\si:\F_\sIp\ra\Obj(\A)$, $\io:\G_\sIp\ra\Mor(\A)$, and
$\pi:\H_\sIp\ra\Mor(\A)$, where
\begin{itemize}
\setlength{\itemsep}{0pt}
\setlength{\parsep}{0pt}
\item[(i)] $\si(J)$ is an object in $\A$ for $J\in\F_\sIp$,
with~$\si(\emptyset)=0$.
\item[(ii)] $\io(J,K):\si(J)\!\ra\!\si(K)$ is injective
for $(J,K)\!\in\!\G_\sIp$, and~$\io(J,J)\!=\!\id_{\si(J)}$.
\item[(iii)] $\pi(J,K)\!:\!\si(J)\!\ra\!\si(K)$ is surjective
for $(J,K)\!\in\!\H_\sIp$, and~$\pi(J,J)\!=\!\id_{\si(J)}$.
\end{itemize}
These should satisfy the conditions:
\begin{itemize}
\setlength{\itemsep}{0pt}
\setlength{\parsep}{0pt}
\item[(A)] Let $(J,K)\in\G_\sIp$ and set $L=K\sm J$. Then the
following is exact in~$\A$:
\begin{equation*}
\xymatrix@C=40pt{ 0 \ar[r] &\si(J) \ar[r]^{\io(J,K)} &\si(K)
\ar[r]^{\pi(K,L)} &\si(L) \ar[r] & 0. }
\end{equation*}
\item[(B)] If $(J,K)\in\G_\sIp$ and $(K,L)\in\G_\sIp$
then~$\io(J,L)=\io(K,L)\ci\io(J,K)$.
\item[(C)] If $(J,K)\in\H_\sIp$ and $(K,L)\in\H_\sIp$
then~$\pi(J,L)=\pi(K,L)\ci\pi(J,K)$.
\item[(D)] If $(J,K)\in\G_\sIp$ and $(K,L)\in\H_\sIp$ then
\begin{equation*}
\pi(K,L)\ci\io(J,K)=\io(J\cap L,L)\ci\pi(J,J\cap L).
\end{equation*}
\end{itemize}

A {\it morphism} $\al:(\si,\io,\pi)\ra(\si',\io',\pi')$ of
$(I,\pr)$-configurations in $\A$ is a collection of morphisms
$\al(J):\si(J)\ra\si'(J)$ for each $J\in\F_\sIp$ satisfying
\begin{align*}
\al(K)\ci\io(J,K)&=\io'(J,K)\ci\al(J)&&
\text{for all $(J,K)\in\G_\sIp$, and}\\
\al(K)\ci\pi(J,K)&=\pi'(J,K)\ci\al(J)&& \text{for all
$(J,K)\in\H_\sIp$.}
\end{align*}
It is an {\it isomorphism} if $\al(J)$ is an isomorphism for
all~$J\in\F_\sIp$.
\label{as3def2}
\end{dfn}

In \cite[Prop.~4.7]{Joyc3} we relate the classes $[\si(J)]$
in~$K_0(\A)$.

\begin{prop} Let\/ $(\si,\io,\pi)$ be an $(I,\pr)$-configuration
in an abelian category $\A$. Then there exists a unique map
$\ka:I\ra K_0(\A)$ such that\/ $[\si(J)]=\sum_{j\in J}\ka(j)$
in $K_0(\A)$ for all f-sets~$J\subseteq I$.
\label{as3prop1}
\end{prop}

Here \cite[Def.s~5.1, 5.2]{Joyc3} are two ways to construct new
configurations.

\begin{dfn} Let $(I,\pr)$ be a finite poset and $J\in\F_\sIp$. Then
$(J,\pr)$ is also a finite poset, and $\F_\sJp,\G_\sJp,\H_\sJp\!
\subseteq\!\F_\sIp,\G_\sIp,\H_\sIp$. Let $(\si,\io,\pi)$ be an
$(I,\pr)$-configuration in an abelian category $\A$. Define the
$(J,\pr)$-{\it subconfiguration} $(\si',\io',\pi')$ of
$(\si,\io,\pi)$ by $\si'\!=\!\si\vert_{\F_\sJp}$,
$\io'\!=\!\io\vert_{\G_\sJp}$ and~$\pi'\!=\!\pi\vert_{\H_\sJp}$.

Let $(I,\pr),(K,\tl)$ be finite posets, and $\phi:I\!\ra\!K$ be
surjective with $i\pr j$ implies $\phi(i)\!\tl \!\phi(j)$. Then
$\phi^{-1}$ maps $\F_\sKt,\G_\sKt,\ab\H_\sKt\!\ra\!\F_\sIp,
\G_\sIp,\H_\sIp$. Let $(\si,\io,\pi)$ be an $(I,\pr)$-configuration
in an abelian category $\A$. Define the {\it quotient\/
$(K,\tl)$-configuration} $(\ti\si,\ti\io,\ti\pi)$ by
$\ti\si(A)\!=\!\si(\phi^{-1}(A))$ for $A\!\in\!\F_\sKt$,
$\ti\io(A,B)\!=\!\io(\phi^{-1}(A),\phi^{-1}(B))$ for
$(A,B)\!\in\!\G_\sKt$, and
$\ti\pi(A,B)\!=\!\pi(\phi^{-1}(A),\ab\phi^{-1}(B))$ for
$(A,B)\!\in\!\H_\sKt$. We call $(\si,\io,\pi)$ a {\it refinement\/}
of~$(\ti\si,\ti\io,\ti\pi)$.
\label{as3def3}
\end{dfn}

Following \cite[Def.~6.1]{Joyc3} we define {\it improvements} and
{\it best configurations}.

\begin{dfn} Let $(I,\pr)$ be a finite poset and $\tl$ a partial
order on $I$ dominating $\pr$, as in Definition \ref{as3def1}. Let
$\A$ be an abelian category. For each $(I,\pr)$-configuration
$(\si,\io,\pi)$ in $\A$ we have a quotient $(I,\tl)$-configuration
$(\ti\si,\ti\io,\ti\pi)$, as in Definition \ref{as3def3} with
$\phi=\id:I\ra I$. We call $(\si,\io,\pi)$ an {\it improvement\/} or
an $(I,\pr)$-{\it improvement} of $(\ti\si,\ti\io,\ti\pi)$, and a
{\it strict improvement\/} if $\pr,\tl$ are distinct. If $\tl$
dominates $\pr$ by $s$ steps we also call $(\si,\io,\pi)$ an $s$
{\it step improvement\/} of $(\ti\si,\ti\io,\ti\pi)$. We call an
$(I,\tl)$-configuration $(\ti\si,\ti\io,\ti\pi)$ {\it best\/} if
there exists no strict improvement $(\si,\io,\pi)$ of
$(\ti\si,\ti\io,\ti\pi)$. Note that improvements are a special kind
of {\it refinement}.
\label{as3def4}
\end{dfn}

In \cite[Prop.~6.9 \& Th.~6.10]{Joyc3} we classify one step
improvements and prove a criterion for best
$(I,\tl)$-configurations. Recall that a short exact sequence $0\ra
X\ra Y\ra Z\ra 0$ in $\A$ is {\it split\/} if there is a compatible
isomorphism~$Y\cong X\op Z$.

\begin{thm} Let\/ $(I,\tl)$ be a finite poset. Call\/ $i,j\in I$
{\bf consecutive} if\/ $i\!\tl\!j$ with\/ $i\!\ne\!j$, but there
exists no\/ $k\in I$ with\/ $i\!\ne\!k\!\ne\!j$ and\/
$i\!\tl\!k\!\tl\!j$. That is, $i,j$ are distinct with\/ $i\tl j$,
and no other\/ $k\in I$ lies between\/ $i,j$ in the order\/~$\tl$.

An $(I,\tl)$-configuration $(\si,\io,\pi)$ in an abelian category
$\A$ is best if and only if for all consecutive\/ $i,j$ in $I$, the
following short exact sequence is not split:
\e
\xymatrix@C=23pt{ 0 \ar[r] & \si\bigl(\{i\}\bigr)
\ar[rr]^(0.45){\io(\{i\},\{i,j\})} && \si\bigl(\{i,j\}\bigr)
\ar[rr]^{\pi(\{i,j\},\{j\})} && \si\bigl(\{j\}\bigr) \ar[r] & 0. }
\label{as3eq1}
\e

Suppose $i,j$ are consecutive, and\/ \eq{as3eq1} is split. Define
$\pr$ on $I$ by $a\pr b$ if\/ $a\tl b$ and\/ $a\!\ne\!i$,
$b\!\ne\!j$, so that\/ $\tl$ dominates $\pr$ by one step. Then the
$(I,\pr)$-improvements of\/ $(\si,\io,\pi)$ are in $1$-$1$
correspondence with\/~$\Hom(\si(\{j\}),\si(\{i\}))$.
\label{as3thm1}
\end{thm}

\subsection{Moduli stacks of configurations}
\label{as32}

Here \cite[Assumptions 7.1 \& 8.1]{Joyc3} is the data we require.

\begin{ass} Let $\K$ be an algebraically closed field and $\A$
a $\K$-linear noetherian abelian category with $\Ext^i(X,Y)$
finite-dimensional $\K$-vector spaces for all $X,Y\in\A$ and $i\ge
0$. Let $K(\A)$ be the quotient of the Grothendieck group $K_0(\A)$
by some fixed subgroup. Suppose that if $X\in\A$ with $[X]=0$ in
$K(\A)$ then~$X\cong 0$.

To define moduli stacks of objects or configurations in $\A$, we
need some {\it extra data}, to tell us about algebraic families of
objects and morphisms in $\A$, parametrized by a base scheme $U$. We
encode this extra data as a {\it stack in exact categories\/}
$\fF_\A$ on the {\it category of\/ $\K$-schemes\/} $\Sch_\K$, made
into a {\it site\/} with the {\it \'etale topology}. The
$\K,\A,K(\A),\fF_\A$ must satisfy some complex additional conditions
\cite[Assumptions 7.1 \& 8.1]{Joyc3}, which we do not give.
\label{as3ass}
\end{ass}

Note that \cite{Joyc3,Joyc4} did not assume $\A$ {\it noetherian},
but we need this to make $\tau$-semistability well-behaved, so we
suppose it from the outset. All the examples of \cite[\S 9--\S
10]{Joyc3} have $\A$ noetherian. Here is some new notation.

\begin{dfn} We work in the situation of Assumption \ref{as3ass}.
Define
\e
C(\A)=\bigl\{[X]\in K(\A):X\in\A,\;\> X\not\cong 0\bigr\}
\subset K(\A).
\label{as3eq2}
\e
That is, $C(\A)$ is the collection of classes in $K(\A)$ of {\it
nonzero objects} $X\in\A$. Note that $C(\A)$ is {\it closed under
addition}, as $[X\op Y]=[X]+[Y]$. Note also that $0\notin C(\A)$, as by
Assumption \ref{as3ass} if $X\not\cong 0$ then $[X]\ne 0$ in~$K(\A)$.

In \cite{Joyc3,Joyc4} we worked mostly with $\bar
C(\A)=C(\A)\cup\{0\}$, the collection of classes in $K(\A)$ of all
objects $X\in\A$. But here and in \cite{Joyc5} we find $C(\A)$ more
useful, as stability conditions will be defined only on nonzero
objects. We think of $C(\A)$ as the `positive cone' and $\bar C(\A)$
as the `closed positive cone' in~$K(\A)$.

Define a set of $\A$-{\it data\/} to be a triple $(I,\pr,\ka)$ such
that $(I,\pr)$ is a finite poset and $\ka:I\ra C(\A)$ a map. We {\it
extend\/ $\ka$ to the set of subsets of\/} $I$ by defining
$\ka(J)=\sum_{j\in J}\ka(j)$. Then $\ka(J)\in C(\A)$ for all
$\emptyset\ne J\subseteq I$, as $C(\A)$ is closed under addition.
Define an $(I,\pr,\ka)$-{\it configuration\/} to be an
$(I,\pr)$-configuration $(\si,\io,\pi)$ in $\A$ with
$[\si(\{i\})]=\ka(i)$ in $K(\A)$ for all $i\in I$. Then $[\si(J)]
=\ka(J)$ for all $J\in\F_\sIp$, by Proposition~\ref{as3prop1}.
\label{as3def5}
\end{dfn}

In the situation above, we define the following $\K$-stacks
\cite[Def.s 7.2 \& 7.4]{Joyc3}:
\begin{itemize}
\setlength{\itemsep}{0pt}
\setlength{\parsep}{0pt}
\item The {\it moduli stacks} $\fObj_\A$ of {\it objects in} $\A$,
and $\fObj^\al_\A$ of {\it objects in $\A$ with class $\al$ in}
$K(\A)$, for each $\al\in\bar C(\A)$. They are algebraic $\K$-stacks,
locally of finite type, with $\fObj_\A^\al$ an open and closed
$\K$-substack of $\fObj_\A$. The underlying geometric spaces
$\fObj_\A(\K),\fObj_\A^\al(\K)$ are the sets of isomorphism
classes of objects $X$ in $\A$, with $[X]=\al$ for~$\fObj_\A^\al(\K)$.
\item The {\it moduli stacks\/} $\fM(I,\pr)_\A$ of $(I,\pr)$-{\it
configurations} and $\fM(I,\pr,\ka)_\A$ of $(I,\pr,\ka)$-{\it
configurations in} $\A$, for all finite posets $(I,\pr)$ and
$\ka:I\ra\bar C(\A)$. They are algebraic $\K$-stacks, locally of
finite type, with $\fM(I,\pr,\ka)_\A$ an open and closed
$\K$-substack of $\fM(I,\pr)_\A$. Write
$\M(I,\pr)_\A,\M(I,\pr,\ka)_\A$ for the underlying geometric spaces
$\fM(I,\pr)_\A(\K),\fM(I,\pr,\ka)_\A(\K)$. Then $\M(I,\pr)_\A$ and
$\M(I,\pr,\ka)_\A$ are the {\it sets of isomorphism classes of\/
$(I,\pr)$- and\/ $(I,\pr,\ka)$-configurations in} $\A$,
by~\cite[Prop.~7.6]{Joyc3}.
\end{itemize}
Each stabilizer group $\Iso_\K([X])$ or $\Iso_\K\bigl([(\si,
\io,\pi)]\bigr)$ in $\fObj_\A$ or $\fM(I,\pr)_\A$ is the group
of invertible elements in the finite-dimensional $\K$-algebra
$\End(X)$ or $\End\bigl((\si,\io,\pi)\bigr)$. Thus $\fObj_\A,
\fObj_\A^\al,\fM(I,\pr)_\A,\fM(I,\pr,\ka)_\A$ have {\it affine
geometric stabilizers}, which is required to use the results
of~\S\ref{as23}--\S\ref{as24}.

In \cite[Def.~7.7 \& Prop.~7.8]{Joyc3} we define 1-morphisms of
$\K$-stacks, as follows:
\begin{itemize}
\setlength{\itemsep}{0pt}
\setlength{\parsep}{0pt}
\item For $(I,\pr)$ a finite poset, $\ka:I\ra\bar C(\A)$ and
$J\in\F_\sIp$, we define $\bs\si(J):\fM(I,\pr)_\A\ra\fObj_\A$
or $\bs\si(J):\fM(I,\pr,\ka)_\A\ra\fObj_\A^{\ka(J)}$. The
induced maps $\bs\si(J)_*:\M(I,\pr)_\A\ra\fObj_\A(\K)$ or
$\M(I,\pr,\ka)_\A\ra\fObj_\A^{\ka(J)}(\K)$ act
by~$\bs\si(J)_*:[(\si,\io,\pi)]\mapsto[\si(J)]$.
\item For $(I,\pr)$ a finite poset, $\ka:I\ra\bar C(\A)$ and
$J\in\F_\sIp$, we define the $(J,\pr)$-{\it subconfiguration
$1$-morphism} $S(I,\pr,J):\fM(I,\pr)_\A\ra\fM(J,\pr)_\A$ or
$S(I,\pr,J):\fM(I,\pr,\ka)_\A\ra\fM(J,\pr,\ka\vert_J)_\A$.
The induced maps $S(I,\pr,J)_*$ act by $S(I,\pr,J)_*:
[(\si,\io,\pi)]\mapsto[(\si',\io',\pi')]$, where
$(\si,\io,\pi)$ is an $(I,\pr)$-configuration in $\A$,
and $(\si',\io',\pi')$ its $(J,\pr)$-subconfiguration.
\item Let $(I,\pr)$, $(K,\tl)$ be finite posets, $\ka:I\!\ra\!\bar C(\A)$,
and $\phi:I\!\ra\!K$ be surjective with $i\pr j$ implies $\phi(i)
\!\tl\!\phi(j)$ for $i,j\in I$. Define $\mu:K\!\ra\!\bar C(\A)$ by
$\mu(k)\!=\!\ka(\phi^{-1}(k))$. The {\it quotient\/
$(K,\tl)$-configuration\/ $1$-morphisms\/} are
\ea
&Q(I,\pr,K,\tl,\phi):\fM(I,\pr)_\A\ra\fM(K,\tl)_\A,
\label{as3eq3}\\
&Q(I,\pr,K,\tl,\phi):\fM(I,\pr,\ka)_\A\ra\fM(K,\tl,\mu)_\A.
\label{as3eq4}
\ea
The induced maps $Q(I,\pr,K,\tl,\phi)_*$ act by $Q(I,\pr,K,\tl,\phi)_*:
[(\si,\io,\pi)]\!\mapsto\![(\ti\si,\ti\io,\ti\pi)]$, where $(\si,\io,\pi)$
is an $(I,\pr)$-configuration in $\A$, and $(\ti\si,\ti\io,\ti\pi)$
its quotient $(K,\tl)$-configuration from $\phi$. When $I\!=\!K$ and
$\phi:I\!\ra\!I$ is the identity $\id_I$, write $Q(I,\pr,\tl)\!=\!
Q(I,\pr,I,\tl,\id_I)$. Then $\mu\!=\!\ka$, so that
\ea
Q(I,\pr,\tl)&:\fM(I,\pr)_\A\ra\fM(I,\tl)_\A,
\label{as3eq5}\\
Q(I,\pr,\tl)&:\fM(I,\pr,\ka)_\A\ra\fM(I,\tl,\ka)_\A.
\label{as3eq6}
\ea
\end{itemize}
{}\hskip\parindent Here \cite[Th.~8.4]{Joyc3} are some properties of
these 1-morphisms:

\begin{thm} {\rm(a)} $Q(I,\pr,K,\tl,\phi),Q(I,\pr,\tl)$ in
\eq{as3eq3}--\eq{as3eq6} are representable, and \eq{as3eq4},
\eq{as3eq6} are of finite type.

\noindent{\rm(b)} $\bs\si(I):\fM(I,\pr)_\A\!\ra\!\fObj_\A$ is representable,
and\/ $\bs\si(I):\fM(I,\pr,\ka)_\A\!\ra\!\fObj^{\ka(I)}_\A$ is
representable and of finite type.

\noindent{\rm(c)} $\prod_{i\in I}\bs\si(\{i\}):\fM(I,\pr)_\A\!\ra\!
\prod_{i\in I}\fObj_\A$ and\/  $\prod_{i\in I}\bs\si(\{i\}):
\fM(I,\pr,\ka)_\A\!\ra\!\prod_{i\in I}\fObj_\A^{\ka(i)}$ are of finite type.
\label{as3thm2}
\end{thm}

We also define some more moduli spaces $\fM(X,I,\pr,\ka)_\A$ in
\cite[Def.~8.5]{Joyc3}, for which $\si(I)$ is a fixed
object~$X\in\A$.

\begin{dfn} In the situation above, let $X\in\A$. Then $X$
corresponds to a 1-morphism $X:\Spec\K\ra\fObj_\A^{[X]}$.
For $\A$-data $(I,\pr,\ka)$ with $\ka(I)=[X]$ in $K(\A)$, define
an algebraic $\K$-stack
\begin{equation*}
\fM(X,I,\pr,\ka)_\A=\fM(I,\pr,\ka)_\A\t_{\bs\si(I),\fObj^{\ka(I)}_\A,
X}\Spec\K.
\end{equation*}
Theorem \ref{as3thm2}(b) implies $\fM(X,I,\pr,\ka)_\A$ is
represented by a {\it finite type algebraic $\K$-space}. Write
$\Pi_X\!:\!\fM(X,I,\pr,\ka)_\A\!\ra\!\fM(I,\pr,\ka)_\A$ for the
1-morphism from the fibre product. It is of finite type. Write
$\M(X,I,\pr,\ka)_\A\!=\!\fM(X,I,\pr,\ka)_\A(\K)$ for the underlying
geometric space. Then \cite[Prop.~8.6]{Joyc3} identifies
$\M(X,I,\pr,\ka)_\A$ with the set of isomorphism classes of
$(I,\pr,\ka)$-configurations $(\si,\io,\pi)$ in $\A$ with
$\si(I)=X$, modulo isomorphisms $\al:(\si,\io,\pi)
\ra(\si',\io',\pi')$ of $(I,\pr)$-configurations
with~$\al(I)=\id_X$.

The 1-morphisms $Q(I,\pr,K,\tl,\phi)$, $Q(I,\pr,\tl)$
above on $\fM(I,\pr,\ka)_\A$ have analogues for
$\fM(X,I,\pr,\ka)_\A$, denoted the same way.
\label{as3def6}
\end{dfn}

In \cite[\S 9--\S 10]{Joyc3} we define the data $\A,K(\A),\fF_\A$ in
some large classes of examples, and prove Assumption \ref{as3ass}
holds in each case.

\subsection{Algebras of constructible and stack functions}
\label{as33}

Next we summarize parts of \cite{Joyc4}, which define and study
associative multiplications $*$ on $\CF(\fObj_\A)$ and
$\SF(\fObj_\A)$, based on {\it Ringel--Hall algebras}.

\begin{dfn} Let Assumption \ref{as3ass} hold with $\K$ of
characteristic zero. Write $\de_{[0]}\in\CF(\fObj_\A)$ for the
characteristic function of $[0]\in\fObj_\A(\K)$. Following
\cite[Def.~4.1]{Joyc4}, using the diagrams of 1-morphisms and
pullbacks, pushforwards
\begin{equation*}
\text{
\begin{footnotesize}
$\displaystyle
\begin{gathered}
\xymatrix@C=54pt@R=10pt{ \fObj_\A\t\fObj_\A & \fM(\{1,2\},\le)_\A
\ar[l]_{\bs\si(\{1\})\t\bs\si(\{2\})} \ar[r]^{\bs\si(\{1,2\})}
& \fObj_\A,\\
\CF\bigl(\fObj_\A\bigr)\t\CF\bigl(\fObj_\A\bigr)
\!\!\!\!\!\!\!\!\!\!\!\!\!\!\!
\ar@<.7ex>[dr]^(0.6){\qquad\qquad(\bs\si(\{1\}))^*\cdot(\bs\si(\{2\}))^*}
\ar[d]_{\ot}
\\
\CF\bigl(\fObj_\A\!\t\!\fObj_\A\bigr)
\ar[r]^(.47){(\bs\si(\{1\})\t\bs\si(\{2\}))^*}
& \CF\bigl(\fM(\{1,2\},\le)_\A\bigr)
\ar[r]^(.56){\CF^\stk(\bs\si(\{1,2\}))}
& \CF\bigl(\fObj_\A\bigr),
}
\end{gathered}
$
\end{footnotesize}}
\end{equation*}
define a bilinear operation $*:\CF(\fObj_\A)\t\CF(\fObj_\A)
\ra\CF(\fObj_\A)$ by
\e
f*g=\CF^\stk(\bs\si(\{1,2\}))\bigl[\bs\si(\{1\})^*(f)\cdot
\bs\si(\{2\})^*(g)\bigr].
\label{as3eq7}
\e
Then \cite[Th.~4.3]{Joyc4} shows $*$ is {\it associative}, and
$\CF(\fObj_\A)$ is a $\Q$-{\it algebra}, with identity $\de_{[0]}$
and multiplication~$*$.

Following \cite[Def.~4.8]{Joyc4}, write $\CFi(\fObj_\A)$ for the
vector subspace of $f$ in $\CF(\fObj_\A)$ {\it supported on
indecomposables}, that is, $f\bigl([X]\bigr)\ne 0$ implies
$0\not\cong X$ is indecomposable. Define a bilinear bracket
$[\,,\,]:\CF(\fObj_\A)\t\CF(\fObj_\A)\ra\CF(\fObj_\A)$ by
$[f,g]=f*g-g*f$. Since $*$ is associative, $[\,,\,]$ satisfies the
{\it Jacobi identity}, and makes $\CF(\fObj_\A)$ into a $\Q$-{\it
Lie algebra}. Then \cite[Th.~4.9]{Joyc4} shows $\CFi(\fObj_\A)$ is
closed under $[\,,\,],$ and so is also a $\Q$-Lie algebra.
\label{as3def7}
\end{dfn}

The next result follows from \cite[Def.~4.13 \& Prop.~4.14]{Joyc4}.
The important point is that $\Phi$ is an {\it isomorphism}, not just
a homomorphism.

\begin{prop} Suppose Assumption \ref{as3ass} holds with\/ $\K$ of
characteristic zero, and use the notation of Definition
\ref{as3def7}. Let\/ $\L$ be a $\Q$-Lie subalgebra of\/
$\CFi(\fObj_\A)$, and\/ $\H_\L$ the $\Q$-subalgebra of\/
$\CF(\fObj_\A)$ with identity generated by $\L$. Write $U(\L)$ for
the universal enveloping algebra of\/ $\L$. Then the inclusion
$\L\subseteq\H_\L$ induces a unique $\Q$-algebra isomorphism
$\Phi:U(\L)\ra\H_\L$ with\/ $\Phi(1)=\de_{[0]}$ and\/
$\Phi(f_1\cdots f_n)= f_1*\cdots*f_n$ for~$f_1,\ldots,f_n\in\L$.
\label{as3prop2}
\end{prop}

In \cite[\S 5]{Joyc4} we extend much of the above to {\it stack
functions}, as in \S\ref{as24}. Here are a few of the basic
definitions and results.

\begin{dfn} Suppose Assumption \ref{as3ass} holds. If
$[(\fR,\rho)]\in\SF(\fObj_\A)$ and $r\in\fR(\K)$ with
$\rho_*(r)=[X]\in\fObj_\A(\K)$ for $X\in\A$, then $\rho$ induces an
injective morphism of stabilizer $\K$-groups
$\rho_*:\Iso_\K(r)\ra\Iso_\K([X])\cong\Aut(X)$, which induces an
isomorphism of $\Iso_\K(r)$ with a $\K$-subgroup of $\Aut(X)$. Now
$\Aut(X)$ is the $\K$-group of invertible elements in the
$\K$-algebra~$\End(X)=\Hom(X,X)$.

As in \cite[Def.~5.5]{Joyc4} define $\SFa(\fObj_\A)$ to be the
subspace of $\SF(\fObj_\A)$ spanned by $[(\fR,\rho)]$ such that for
all $r\in\fR(\K)$ with $\rho_*(r)=[X]$, the $\K$-subgroup
$\rho_*\bigl(\Iso_\K(r)\bigr)$ in $\Aut(X)$ is the $\K$-group of
invertible elements in a $\K$-subalgebra of $\End(X)$. Then
$\io_{\fObj_\A}$ in Definition \ref{as2def8}
maps~$\CF(\fObj_\A)\!\ra\!\SFa(\fObj_\A)$.

By analogy with \eq{as3eq7}, using $\fM(\{1,2\},\le)_\A$ define
\cite[Def.~5.1]{Joyc4} a bilinear operation $*:\SF(\fObj_\A)\t
\SF(\fObj_\A)\ra\SF(\fObj_\A)$~by
\e
f*g=\bs\si(\{1,2\})_*\bigl[(\bs\si(\{1\})\t\bs\si(\{2\}))^*(f\ot
g)\bigr].
\label{as3eq8}
\e
Write $\bde_{[0]}\!\in\!\SFa(\fObj_\A)$ for $\bde_C$ in Definition
\ref{as2def8} with~$C\!=\!\{[0]\}$.

Then \cite[Th.~5.2 \& Prop.~5.6]{Joyc4} show that $\SF(\fObj_\A)$ is
a $\Q$-algebra with associative multiplication $*$ and identity
$\bde_{[0]}$, and $\SFa(\fObj_\A)$ is closed under $*$ and so is a
$\Q$-subalgebra. When $\K$ has characteristic zero,
\e
\pi_{\fObj_\A}^\stk:\SF(\fObj_\A)\longra\CF(\fObj_\A)
\label{as3eq9}
\e
is a $\Q$-algebra morphism, where $\CF(\fObj_\A)$ is an algebra as
in Definition~\ref{as3def7}.
\label{as3def8}
\end{dfn}

\begin{dfn} Let Assumption \ref{as3ass} hold. Following
\cite[Def.~5.13]{Joyc4}, define $\SFai(\fObj_\A)$ to be the subspace
of $f\in\SFa(\fObj_\A)$ with $\Pi^\vi_1(f)=f$, where $\Pi^\vi_1$ is
the operator of \cite[\S 5.2]{Joyc2}, interpreted as projecting to
stack functions `supported on virtual indecomposables'. Write
$[f,g]=f*g-g*f$ for $f,g\in\SFa(\fObj_\A)$. As $*$ is associative
$[\,,\,]$ satisfies the {\it Jacobi identity}, and makes
$\SFa(\fObj_\A)$ into a $\Q$-{\it Lie algebra}. Then
\cite[Th.~5.17]{Joyc4} shows $\SFai(\fObj_\A)$ is closed under
$[\,,\,]$, and is a {\it Lie subalgebra}. When $\cha\K=0$,
\eq{as3eq9} restricts to a Lie algebra morphism
\e
\pi_{\fObj_\A}^\stk:\SFai(\fObj_\A)\longra\CFi(\fObj_\A).
\label{as3eq10}
\e
\label{as3def9}
\end{dfn}

The above material also works for the other stack function spaces on
$\fObj_\A$, in particular for $\oSF(\fObj_\A,\Up,\La)$,
$\oSF(\fObj_\A,\Up,\La^\ci)$ and $\oSF(\fObj_\A,\Th,\Om)$, giving
algebras $\oSF,\oSFa(\fObj_\A,*,*)$ and Lie algebras
$\oSFai(\fObj_\A,*,*)$.

\section{Stability conditions}
\label{as4}

We now introduce our concepts of ({\it weak\/}) {\it stability
condition} $(\tau,T,\le)$ on $\A$, which are based on the
stability conditions of Rudakov \cite{Ruda}. Perhaps their most
important properties are Theorems \ref{as4thm1} and \ref{as4thm2}
below. These show that for a weak stability condition
$(\tau,T,\le)$ with $\A$ noetherian and $\tau$-artinian, each
$X\in\A$ may be decomposed into $\tau$-semistable factors $S_k$
in a unique way, and if $(\tau,T,\le)$ is a stability condition
the $S_k$ can be further split into $\tau$-stable pieces. One
moral of this is that $\tau$-{\it stability is well-behaved
for stability conditions, but badly behaved for weak stability
conditions}.

\subsection{Definitions and basic properties}
\label{as41}

Here is our notion of ({\it weak\/}) {\it stability condition},
generalizing Rudakov~\cite{Ruda}.

\begin{dfn} Let $\A$ be an abelian category, $K(\A)$ be
the quotient of $K_0(\A)$ by some fixed subgroup, and
$C(\A)$ as in \eq{as3eq2}. Suppose $(T,\le)$ is a totally
ordered set, and $\tau:C(\A)\ra T$ a map. We call $(\tau,T,\le)$
a {\it stability condition} on $\A$ if whenever $\al,\be,\ga\in
C(\A)$ with $\be=\al+\ga$ then either $\tau(\al)\!<\!\tau(\be)
\!<\!\tau(\ga)$, or $\tau(\al)\!>\!\tau(\be)\!>\!\tau(\ga)$, or
$\tau(\al)\!=\!\tau(\be)\!=\!\tau(\ga)$. We call $(\tau,T,\le)$
a {\it weak stability condition} on $\A$ if whenever $\al,\be,
\ga\in C(\A)$ with $\be=\al+\ga$ then either $\tau(\al)\!\le\!
\tau(\be)\!\le\!\tau(\ga)$, or $\tau(\al)\!\ge\!\tau(\be)\!\ge
\!\tau(\ga)$. Clearly, a stability condition is a weak stability
condition, but not necessarily vice versa.
\label{as4def1}
\end{dfn}

Our stability conditions are motivated by, and more-or-less
equivalent to, Rudakov's \cite[Def.~1.1]{Ruda}. The difference
is that Rudakov's stability conditions are {\it preorders} on
the nonzero objects of $\A$. In effect our definition requires
Rudakov's preorder to factor through the map $\Obj(\A)\ra K(\A)$,
$X\mapsto[X]$, and so amounts to a preorder on $C(\A)$. Rudakov
calls the trichotomy $\tau(\al)\!<\!\tau(\be)\!<\!\tau(\ga)$ or
$\tau(\al)\!>\!\tau(\be)\!>\!\tau(\ga)$ or $\tau(\al)\!=\!\tau
(\be)\!=\!\tau(\ga)$ the {\it seesaw inequality}.

In the same way, we call the alternative $\tau(\al)\!\le\!\tau
(\be)\!\le\!\tau(\ga)$ or $\tau(\al)\!\ge\!\tau(\be)\!\ge\!
\tau(\ga)$ the {\it weak seesaw inequality}. As far as I know this
abstract idea of weak stability condition is new. I believe it is a
useful innovation, since as we shall see in \S\ref{as44} important
concepts such as the torsion filtration and $\mu$-(semi)stability of
sheaves are examples of weak stability conditions which are not
stability conditions. Also, to transform between two stability
conditions in \cite{Joyc5} we will need to go via a weak stability
condition.

We use many ordered sets in the paper: finite posets $(I,\pr),
(J,\ls),(K,\tl)$ for $(I,\pr)$-configurations, and now total orders
$(T,\le)$ for stability conditions. As the number of order symbols
is limited, here and in \cite{Joyc5} we will always use `$\le$' for
the total order, so that $(\tau,T,\le)$, $(\ti\tau,\ti T,\le)$ may
denote two different stability conditions, with two {\it
different\/} total orders on $T,\ti T$ both denoted by~`$\le$'.

We define $\tau$-{\it semistable}, $\tau$-{\it stable} and
$\tau$-{\it unstable} objects.

\begin{dfn} Let $(\tau,T,\le)$ be a weak stability condition on
$\A,K(\A)$ as above. Then we say that a nonzero object $X$ in $\A$ is
\begin{itemize}
\setlength{\itemsep}{0pt}
\setlength{\parsep}{0pt}
\item[(i)] $\tau$-{\it semistable} if for all $S\subset X$ with
$S\not\cong 0,X$ we have $\tau([S])\le\tau([X/S])$;
\item[(ii)] $\tau$-{\it stable} if for all $S\subset X$ with
$S\not\cong 0,X$ we have $\tau([S])<\tau([X/S])$; and
\item[(iii)] $\tau$-{\it unstable} if it is not $\tau$-semistable.
\end{itemize}
If $S\subset X$ is a subobject with $S\ne 0,X$ then $[S],[X],[X/S]
\in C(\A)$ with $[X]=[S]+[X/S]$. Thus, if $(\tau,T,\le)$ is a {\it
stability condition} then $\tau([S])\!\le\!\tau([X/S])$ in (i) is
equivalent to $\tau([S])\!\le\!\tau([X])$ and to $\tau([X])\!\le\!
\tau([X/S])$, and $\tau([S])\!<\!\tau([X/S])$ in (ii) is equivalent
to $\tau([S])\!<\!\tau([X])$ and to~$\tau([X])\!<\!\tau([X/S])$.
\label{as4def2}
\end{dfn}

We will need the following weakening of $\A$ {\it artinian} in
Definition~\ref{as2def2}.

\begin{dfn} Let $(\tau,T,\le)$ be a weak stability condition on
$\A,K(\A)$. We say $\A$ is $\tau$-{\it artinian} if there exist
no infinite chains of subobjects $\cdots\!\subset\!A_2\!\subset\!
A_1\!\subset\!X$ in $\A$ with $A_{n+1}\!\ne\!A_n$ and
$\tau([A_{n+1}])\!\ge\!\tau([A_n/A_{n+1}])$ for all $n$. If
$(\tau,T,\le)$ is a {\it stability condition} $\tau([A_{n+1}])\!
\ge\!\tau([A_n/A_{n+1}])$ is equivalent to $\tau([A_{n+1}])\!\ge
\!\tau([A_n])$, and the definition reduces to~\cite[Def.~1.7]{Ruda}.
\label{as4def3}
\end{dfn}

In the next theorem we call $0=A_0\subset A_1\subset\cdots\subset
A_n=X$ a {\it Harder--Narasimhan filtration}, as it generalizes
the filtrations constructed by Harder and Narasimhan \cite{HaNa}
for vector bundles over algebraic curves. The proof is adapted
from Rudakov \cite[Th.~2]{Ruda}, which implies the result for
stability conditions.

\begin{thm} Let\/ $(\tau,T,\le)$ be a weak stability condition
on an abelian category $\A$. Suppose $\A$ is noetherian and\/
$\tau$-artinian. Then each\/ $X\in\A$ admits a unique filtration
$0\!=\!A_0\!\subset\!\cdots\!\subset\!A_n\!=\!X$ for $n\ge 0$, such
that\/ $S_k\!=\!A_k/A_{k-1}$ is $\tau$-semistable for $k=1,\ldots,n$,
and\/~$\tau([S_1])>\tau([S_2])>\cdots>\tau([S_n])$.
\label{as4thm1}
\end{thm}

\begin{proof} For $X\cong 0$ the result is trivial with $n=0$,
so fix $X\in\A$ with $X\not\cong 0$. We divide the proof into the
following seven steps:
\begin{list}{}{
\setlength{\itemsep}{1pt}
\setlength{\parsep}{1pt}
\setlength{\labelwidth}{50pt}
\setlength{\leftmargin}{50pt}
}
\item[{\bf Step 1.}] Given $0\ne B\subset X$, there exists
$0\ne A\subset B\subset X$ with $A$ $\tau$-semistable
and~$\tau([A])\ge\tau([B])$.
\item[{\bf Step 2.}] Suppose $0\ne A,B\subseteq X$ with $A$
$\tau$-semistable and $\tau([A])\ge\tau([B])$. Then~$\tau([A+B])
\ge\tau([B])$.
\item[{\bf Step 3.}] Call $0\ne C\subset X$ {\it greedy in} $X$
if $0\ne A\subset X$ with $A$ $\tau$-semistable and $\tau([A])\ge
\tau([C])$ implies $A\subset C$. Then for any $0\ne B\subset X$
there exists $C\subset X$ greedy in $X$ with~$\tau([C])\ge\tau([B])$.
\item[{\bf Step 4.}] There exist unique $\tau^{\rm max}\in T$ and
(not necessarily unique) $0\ne B\subset X$ with $\tau([B])=\tau^{\rm
max}$, such that if $0\ne A\subset X$ with $A$ $\tau$-semistable then
$\tau([A])\le\tau^{\rm max}$. We can choose $B$ $\tau$-semistable.
\item[{\bf Step 5.}] If $0\ne A,B\subset X$ are $\tau$-semistable
with $\tau([A])=\tau([B])=\tau^{\rm max}$, then $A+B\subset X$
is $\tau$-semistable with~$\tau([A+B])=\tau^{\rm max}$.
\item[{\bf Step 6.}] There exists a unique $\tau$-semistable
$0\ne S_1\subset X$ with $\tau([S_1])=\tau^{\rm max}$, such
that if $A\subset X$ is $\tau$-semistable with $\tau([A])=
\tau^{\rm max}$ then~$A\subset S_1$.
\item[{\bf Step 7.}] Complete the proof.
\end{list}

\noindent{\bf Step 1.} Suppose for a contradiction there exists
no such $A$. Set $B_1=B$, and construct by induction a sequence
$\cdots B_2\subset B_1\subset X$ with $B_{j+1}\ne 0,B_j$ and
$\tau([B_{j+1}])\ge\tau([B_j/B_{j+1}])$ as follows. Having
chosen $B_j$, if $j>1$ then $\tau([B_j])\ge\tau([B_{j-1}/B_j])$
implies $\tau([B_j])\ge\tau([B_{j-1}])$ by the weak seesaw inequality.
So $\tau([B_j])\ge\cdots\ge\tau([B_1])=\tau([B])$. As $A=B_j$ will
not do, $B_j$ cannot be $\tau$-semistable. Thus $B_{j+1}$ exists
as we want by Definition \ref{as4def2}(i). But the sequence
$\cdots B_2\subset B_1\subset X$ contradicts $\A$ $\tau$-{\it
artinian} in Definition~\ref{as4def3}.
\smallskip

\noindent{\bf Step 2.} Let $A,B$ be as above. If $A\subset B$ then
$A+B=B$ and $\tau([A+B])=\tau([B])$, so suppose $A\not\subset B$.
Then $A\cap B$ is a proper subobject of $A$, so $A/(A\cap B)\not
\cong 0$, and $A$ $\tau$-semistable implies $\tau\bigl(\bigl[A/
(A\cap B)\bigr]\bigr)\ge\tau([A])$. But $(A+B)/B\cong A/(A\cap B)$
by properties of subobjects in an abelian category. Thus $\tau\bigl(
\bigl[(A+B)/B\bigr]\bigr)\ge\tau([A])\ge\tau([B])$, so $\tau([A+B])
\ge\tau([B])$ by the weak seesaw inequality.
\smallskip

\noindent{\bf Step 3.} Suppose for a contradiction there exists
no such $C$. Construct by induction a sequence $B=B_1\!\subset\!B_2\!
\subset\!\cdots\!\subset\!X$ with $B_j\ne B_{j+1}$ and $\tau([B_{j+1}])
\ge\tau([B_j])$, as follows. Having chosen $B_j$, as
$\tau([B_j])\ge\cdots\ge\tau([B_1])=\tau([B])$, and $C=B_j$ will
not do, $B_j$ cannot be greedy. Thus there exists $\tau$-semistable
$A\subset X$ with $\tau([A])\ge\tau([B_j])$ but $A\not\subset B_j$.
Define $B_{j+1}=A+B_j$. Then $B_{j+1}\ne B_j$ as $A\not\subset B_j$,
and $\tau([B_{j+1}])\ge\tau([B_j])$ by Step 2, completing the
induction. But $B_1\!\subset\!B_2\!\subset\!\cdots\!\subset\!X$
contradicts $\A$ {\it noetherian} in Definition~\ref{as2def2}.
\smallskip

\noindent{\bf Step 4.} Suppose for a contradiction that no such
(not yet unique) $\tau^{\rm max}$ and (not necessarily
$\tau$-semistable) $B$ exist. Construct by induction a sequence
$\cdots\subset C_2\subset C_1\subset X$ with $C_j$ greedy and
$\tau([C_{j+1}])>\tau([C_j])$ for all $j$, as follows. Set
$C_1=X$, which is greedy. Having chosen $C_j$, as $\tau^{\rm
max}=\tau([C_j])$ and $B=C_j$ will not do, there exists a
$\tau$-semistable $A\subset X$ with~$\tau([A])>\tau([C_j])$.

Then $A\subset C_j$, as $C_j$ is greedy. By Step 3 with $C_j$
in place of $X$, there exists $C_{j+1}\subset C_j$ greedy in
$C_j$ with $\tau([C_{j+1}])\ge\tau([A])>\tau([C_j])$. Suppose
$0\ne A'\subset X$ is $\tau$-semistable with $\tau([A'])\ge
\tau([C_{j+1}])$. Then $\tau([A'])\ge\tau([C_j])$, so $A'\subset
C_j$ as $C_j$ is greedy in $X$, and thus $A'\subset C_{j+1}$ as
$C_{j+1}$ is greedy in $C_j$. Hence $C_{j+1}$ is greedy in $X$,
completing the inductive step.

But $\tau([C_{j+1}])>\tau([C_j])$ implies $\tau([C_{j+1}])
>\tau([C_j/C_{j+1}])$ by the weak seesaw inequality, so
$\cdots\subset C_2\subset C_1\subset X$ contradicts $\A$
$\tau$-artinian. Thus $\tau^{\rm max},B$ exist. Step 1
shows there exists $0\ne A\subset B\subset X$ with $A$
$\tau$-semistable and $\tau([A])\ge\tau([B])=\tau^{\rm max}$.
But by definition $\tau([A])\le\tau^{\rm max}$, so
$\tau([A])=\tau^{\rm max}$. Therefore $\tau^{\rm max}$ is
the maximum value in $T$ of $\tau([A])$ for $\tau$-semistable
$0\ne A\subset X$, so $\tau^{\rm max}$ is unique, and replacing
$B$ by $A$, we can choose $B$ $\tau$-semistable.
\smallskip

\noindent{\bf Step 5.} Suppose $S\subset A+B$ with $S\ne A+B$.
Properties of subobjects in abelian categories give isomorphisms
$(S+A)/S\cong A/(S\cap A)$ and $(A+B)/(S+A)\cong B/((S+A)\cap B)$.
Thus from the exact sequence $0\ra(S+A)/S\ra(A+B)/S\ra(A+B)/(S+A)
\ra 0$ we obtain an exact sequence
\e
0\longra A/(S\cap A) \longra (A+B)/S \longra B/((S+A)\cap B)
\longra 0.
\label{as4eq1}
\e
The weak seesaw and $A,B$ $\tau$-semistable give
$\tau\bigl(\bigl[A/(S\cap A)\bigr]\bigr)\ge\tau([A])=
\tau^{\rm max}$ if $A/(S\cap A)\ne0$, and $\tau\bigl(
\bigl[B/((S+A)\cap B)\bigr]\bigr)\ge\tau([B])=
\tau^{\rm max}$ if $B/((S+A)\cap B)\ne0$. From
\eq{as4eq1} and the weak seesaw we deduce~$\tau\bigl([(A+B)/S]
\bigr)\ge\tau^{\rm max}$.

In particular, for $S=0$ we have $\tau([A+B])\ge\tau^{\rm max}$.
But $S\subset X$ and $A+B\subset X$, so by Steps 1 and 4 we see
that $\tau([S])\le\tau^{\rm max}$ if $S\ne 0$, and $\tau([A+B])\le
\tau^{\rm max}$. Hence $\tau([A+B])=\tau^{\rm max}$, as we want,
and if $S\ne 0$ then $\tau([S])\le\tau^{\rm max}\le\tau\bigl([
(A+B)/S]\bigr)$, which implies $A+B$ is $\tau$-semistable.
\smallskip

\noindent{\bf Step 6.} Suppose for a contradiction no such $S_1$
exists. Construct by induction a sequence $B_1\subset B_2\subset
\cdots\subset X$ with $B_j\ne B_{j+1}$ and $B_j$ $\tau$-semistable
with $\tau([B_j])=\tau^{\rm max}$, as follows. Set $B_1=B$ from
Step 4, chosen $\tau$-semistable. Having chosen $B_j$, as $S_1=B_j$
will not do there exists $\tau$-semistable $A\subset X$ with
$\tau([A])=\tau^{\rm max}$ and~$A\not\subset B_j$.

Define $B_{j+1}=A+B_j$. Then $B_{j+1}$ is $\tau$-semistable with
$\tau([B_{j+1}])=\tau^{\rm max}$ by Step 5, and $B_{j+1}\ne B_j$
as $A\not\subset B_j$, completing the inductive step. But $B_1
\subset B_2\subset\cdots\subset X$ contradicts $\A$ noetherian,
so $S_1$ exists. If $S_1,S_1'$ satisfy the conditions then
$S_1\subset S_1'$ and $S_1'\subset S_1$, so $S_1=S_1'$ and
$S_1$ is unique.
\smallskip

\noindent{\bf Step 7.} By induction construct a sequence $0=A_0
\subset A_1\subset\cdots\subset A_n=X$ with $0\ne S_j=A_j/A_{j-1}$
$\tau$-semistable, as follows. Set $A_0=0$ and $A_1\subset X$ to
be $S_1$ from Step 6. Then $S_1=A_1/A_0$ is $\tau$-semistable.
Having constructed $A_j$, if $A_j=X$ then set $n=j$ and finish.
Otherwise define $A_{j+1}$ such that $A_j\subset A_{j+1}\subset X$
and $A_{j+1}/A_j\subset X/A_j$ is the subobject $S_1$ given by
Step 6 with $X/A_j$ in place of $X$. Then $S_{j+1}=A_{j+1}/A_j$
is nonzero and $\tau$-semistable.

As $A_{j+1}\ne A_j$ and $\A$ is noetherian the sequence must
terminate at some $n$, so $0=A_0\subset\cdots\subset A_n=X$ is
well-defined. Suppose for a contradiction that $\tau([S_j])\le
\tau([S_{j+1}])$. Then we have subobjects $S_j=A_j/A_{j-1}\subset
X/A_{j-1}$ and $A_{j+1}/A_{j-1}\subset X/A_{j-1}$, with
$(A_{j+1}/A_{j-1})/S_j\cong S_{j+1}$. Write $\tau^{\rm max}_j=
\tau([S_j])$. Then $\tau([S_{j+1}])\ge\tau^{\rm max}_j$, so
the weak seesaw implies $\tau([A_{j+1}/A_{j-1}])\ge\tau^{\rm
max}_j$, and an argument similar to Step 5 shows $A_{j+1}/A_{j-1}$
is $\tau$-semistable. Hence $A_{j+1}/A_{j-1}\subset S_j$ by
definition of $S_j$, giving $S_{j+1}=0$, a contradiction.

Therefore $\tau([S_1])>\tau([S_2])>\cdots>\tau([S_n])$, as we
want. It remains only to prove $0=A_0\subset\cdots\subset A_n=X$
is unique. But it is easy to show that for a filtration satisfying
the conditions of the theorem, the subobject $S_j\subset X/A_{j-1}$
satisfies the conditions of $S_1$ in Step 6 with $X/A_{j-1}$ in
place of $X$. Thus, having chosen $A_{j-1}$, uniqueness in Step 6
implies $S_j=A_j/A_{j-1}$ and $A_j$ are uniquely determined, so
uniqueness follows by induction on~$j$.
\end{proof}

Theorem \ref{as4thm1} justifies the weak case in Definition
\ref{as4def1}, as it shows that $\tau$-{\it semistability is
well-behaved for weak stability conditions}. However the next
result, which follows from Rudakov \cite[Th.~3]{Ruda}, is {\it
false} for weak stability conditions $(\tau,T,\le)$, as one can show
by example. One moral is that $\tau$-{\it stability is well-behaved
for stability conditions, but badly behaved for weak stability
conditions}. Therefore in \S\ref{as5} and \S\ref{as6} below, which
deal with $\tau$-stability, we will consider only stability
conditions, not weak stability conditions.

\begin{thm} Let\/ $\A$ be an abelian category, and\/ $(\tau,T,\le)$
a stability condition on $\A,K(\A)$. Suppose $\A$ is noetherian
and\/ $\tau$-artinian. Then each\/ $\tau$-semistable $X\in\A$ admits
a filtration\/ $0\!=\!A_0\!\subset\!A_1\!\subset\!\cdots\!\subset\!
A_n\!=\!X$ for $n\ge 1$, such that\/ $S_k\!=\!A_k/A_{k-1}$ is
$\tau$-stable for $1\le k\le n$, with\/ $\tau([S_1])\!=\!\cdots\!=\!
\tau([S_n])\!=\!\tau([X])$. Suppose $0\!=\!A_0\!\subset\!\cdots\!\subset
\!A_n\!=\!X$ and\/ $0\!=\!B_0\!\subset\!\cdots\!\subset\!B_m\!=\!X$
are two such filtrations with\/ $\tau$-stable factors $S_k\!=\!A_k/A_{k-1}$
and\/ $T_k\!=\!B_k/B_{k-1}$. Then $n=m$, and for some permutation $\si$
of\/ $1,\ldots,n$ we have $S_k\cong T_{\si(k)}$ for~$1\le k\le n$.
\label{as4thm2}
\end{thm}

The restriction to {\it noetherian} $\A$ in these two theorems
is unnecessarily strong. Rudakov only assumes $\A$ is `weakly
noetherian' \cite[Def.~1.12]{Ruda}. But Rudakov's condition
seems unsatisfactory to the author, so we shall not use it.

\subsection{Permissible stability conditions}
\label{as42}

The following notation will be used throughout the rest of the paper.

\begin{dfn} Let Assumption \ref{as3ass} hold and $(\tau,T,\le)$
be a weak stability condition on $\A$. Then $\fObj_\A^\al$ is an
algebraic $\K$-stack for $\al\in C(\A)$, with $\fObj_\A^\al(\K)$
the set of isomorphism classes of $X\in\A$ with class $\al$ in
$K(\A)$. Define
\e
\begin{split}
\Oss^\al(\tau)&=\bigl\{[X]\in\fObj_\A^\al(\K):
\text{$X$ is $\tau$-semistable}\bigr\},\\
\Osi^\al(\tau)&=\bigl\{[X]\in\fObj_\A^\al(\K): \text{$X$ is
$\tau$-semistable and indecomposable}\bigr\},\\
\Ost^\al(\tau)&=\bigl\{[X]\in\fObj_\A^\al(\K): \text{$X$ is
$\tau$-stable}\bigr\}.
\end{split}
\label{as4eq2}
\e

Let $(I,\pr,\ka)$ be $\A$-data, as in Definition \ref{as3def5}, and
$X\in\A$ with $[X]=\ka(I)$. From \S\ref{as32} we have algebraic
$\K$-stacks $\fM(I,\pr, \ka)_\A,\fM(X,I,\pr,\ka)_\A$ such that
$\M(I,\pr,\ka)_\A\!=\! \fM(I,\pr,\ka)_\A(\K)$,
$\M(X,I,\pr,\ka)_\A\!=\!\fM(X,I,\pr, \ka)_\A(\K)$ are sets of
isomorphism classes $[(\si,\io,\pi)]$ of
$(I,\pr,\ka)$-configurations $(\si,\io,\pi)$ in $\A$, with
$\si(I)=X$ in the second case. Define an $(I,\pr,\ka)$-configuration
$(\si,\io,\pi)$ to be $\tau$-{\it semistable} if $\si(\{i\})$ is
$\tau$-semistable, $\tau$-{\it semistable-indecomposable} if
$\si(\{i\})$ is $\tau$-semistable and indecomposable, and
$\tau$-{\it stable} if $\si(\{i\})$ is $\tau$-stable, for all $i\in
I$. Define
\e
\begin{split}
&\Mss,\Msi,\Mst,\Mssb,\Msib,\Mstb(I,\pr,\ka,\tau)_\A
\subseteq\M(I,\pr,\ka)_\A
\quad\text{and}\\
&\Mss,\Msi,\Mst,\Mssb,\Msib,\Mstb(X,I,\pr,\ka,\tau)_\A
\subseteq\M(X,I,\pr,\ka)_\A
\end{split}
\label{as4eq3}
\e
to be the subsets of $[(\si,\io,\pi)]$ with $(\si,\io,\pi)$
$\tau$-{\it semistable} in the $\Mss,\Mssb(\cdots)_\A$ cases, and
$\tau$-{\it semistable-indecomposable} in the
$\Msi,\Msib(\cdots)_\A$ cases, and $\tau$-{\it stable} in the
$\Mst,\Mstb(\cdots)_\A$ cases, and {\it best\/} in the
$\Mssb,\Msib,\Mstb(\cdots)_\A$ cases, as in Definition
\ref{as3def4}. Write
$\dss^\al,\dsi^\al,\dst^\al(\tau):\fObj_\A^\al(\K)$ or
$\fObj_\A(\K)\ra\{0,1\}$ for the {\it characteristic functions} of
$\Oss^\al,\Osi^\al,\Ost^\al(\tau)$. Write
$\dss,\ab\dsi,\ab\dst,\ab\dssb,\ab\dsib,\ab\dstb(I,\pr,\ka,\tau):
\M(I,\pr,\ka)_\A\ab\ra\ab\{0,1\}$ for the characteristic functions
of $\Mss$, $\Msi$, $\Mst$, $\Mssb$, $\Msib$,
$\Mstb(I,\pr,\ka,\tau)_\A$, and
$\dss,\ldots,\dstb(X,I,\pr,\ka,\tau):\ab\M(X,I,\pr,\ka)_\A\ra\{0,1\}$
for those of~$\Mss,\ldots,\Mstb(X,I,\pr,\ka,\tau)_\A$.

Using \cite[Ass.~7.1(iii)]{Joyc3} we see
$\Oss^\al,\Osi^\al,\Ost^\al(\tau)$ are {\it open} in the natural
topology on $\fObj_\A(\K)$, and so are {\it locally constructible}.
Being best is also an open condition on configurations. Therefore
$\Mss,\ldots,\ab\Mstb\ab(I,\pr,\ka,\tau)_\A$ and
$\Mss,\ldots,\Mstb(X,I,\pr,\ka,\tau)_\A$ are locally constructible,
so that
\e
\begin{gathered}
\dss^\al,\dsi^\al,\dst^\al(\tau)\!\in\!\LCF(\fObj_\A),\;\>
\dss,\ldots,\dstb(I,\pr,\ka,\tau)\!\in\!
\LCF\bigl(\fM(I,\pr,\ka)_\A\bigr),\\
\text{and}\qquad
\dss,\dsi,\dst,\dssb,\dsib,\dstb(X,I,\pr,\ka,\tau)\in
\LCF\bigl(\fM(X,I,\pr,\ka)_\A\bigr).
\end{gathered}
\label{as4eq4}
\e
\label{as4def4}
\end{dfn}

We want \eq{as4eq2} and \eq{as4eq3} to be constructible sets, so
that \eq{as4eq4} are constructible functions, as in \S\ref{as23}.
To do this we must impose some assumptions on~$(\tau,T,\le)$.

\begin{dfn} Let Assumption \ref{as3ass} hold and $(\tau,T,\le)$ be
a weak stability condition on $\A$. We call $(\tau,T,\le)$ {\it
permissible} if:
\begin{itemize}
\setlength{\itemsep}{0pt}
\setlength{\parsep}{0pt}
\item[(i)] $\A$ is $\tau$-artinian, in the sense of Definition
\ref{as4def3}, and
\item[(ii)] $\Oss^\al(\tau)$ is a constructible subset in
$\fObj_\A^\al$ for all~$\al\in C(\A)$.
\end{itemize}
\label{as4def5}
\end{dfn}

\begin{thm} Let Assumption \ref{as3ass} hold and\/ $(\tau,T,\le)$
be a permissible weak stability condition on $\A$. Then
$\Osi^\al,\Ost^\al(\tau)$ are constructible sets in $\fObj_\A$ for
all\/ $\al\!\in\!C(\A)$. Suppose\/ $(I,\pr,\ka)$ is\/ $\A$-data
and\/ $X\!\in\!\A$ with\/ $[X]=\ka(I)$ in $K(\A)$. Then\/
$\Mss,\Msi,\Mst,\Mssb,\Msib,\Mstb(I,\pr,\ka,\tau)_\A$ are
constructible in\/ $\fM(I, \pr,\ka)_\A$, and\/
$\Mss,\ldots,\Mstb(X,I,\pr,\ka,\tau)_\A$ in $\fM(X,I,\pr,\ka)_\A$.
Hence
\e
\begin{gathered}
\dss^\al,\dsi^\al,\dst^\al(\tau)\!\in\!\CF(\fObj_\A),\;\>
\dss,\ldots,\dstb(I,\pr,\ka,\tau)\!\in\!
\CF\bigl(\fM(I,\pr,\ka)_\A\bigr),\\
\text{and}\qquad
\dss,\dsi,\dst,\dssb,\dsib,\dstb(X,I,\pr,\ka,\tau)\in
\CF\bigl(\fM(X,I,\pr,\ka)_\A\bigr).
\end{gathered}
\label{as4eq5}
\e
\label{as4thm3}
\end{thm}

\begin{proof} $\prod_{i\in I}\Oss^{\ka(i)}(\tau)$ is constructible in
$\prod_{i\in I}\fObj_\A^{\ka(i)}$ by Definition \ref{as4def5}(ii).
But $\prod_{i\in I}\bs\si(\{i\})\!:\!\fM(I,\pr,\ka)_\A\!\ra\!
\prod_{i\in I}\fObj_\A^{\ka(i)}$ is finite type by Theorem
\ref{as3thm2}(c), and pulls back constructible sets to constructible
sets. Thus $\Mss(I,\pr,\ka,\tau)_\A\!=\!\bigl(\prod_{i\in
I}\bs\si(\{i\})\bigr)_*^{-1}\bigl(\prod_{i\in
I}\Oss^{\ka(i)}(\tau)\bigr)$ is constructible in
$\fM(I,\pr,\ka)_\A$. By Definition \ref{as4def4},
$\Osi^\al,\Ost^\al(\tau)$ are locally constructible subsets of
$\Oss^\al(\tau)$, which is constructible by Definition
\ref{as4def5}(ii), and $\Msi,\Mst,\Mssb,\Msib,\Mstb
(I,\pr,\ka,\tau)_\A$ are locally constructible subsets of $\Mss
(I,\pr,\ka,\tau)_\A$, which is constructible from above, so all
these sets are constructible. As $\Pi_X$ in Definition \ref{as3def6}
is of finite type and $\Mss,\ldots,\Mstb(X,I,\pr,\ka,\tau)_\A$ are
$\Pi_X^*$ of $\Mss,\ldots,\Mstb(I,\pr,\ka,\tau)_\A$, they too are
constructible. Equation \eq{as4eq5} is immediate.
\end{proof}

Here is a useful finiteness property of permissible stability conditions.

\begin{prop} In the situation above, let\/ $(\tau,T,\le)$ be
permissible. Then for each\/ $\al\in C(\A)$, there are only
finitely many pairs $\be,\ga\in C(\A)$ with\/ $\al=\be+\ga$,
$\tau(\al)=\tau(\be)=\tau(\ga)$ and\/~$\Oss^\be(\tau)\ne
\emptyset\ne\Oss^\ga(\tau)$.
\label{as4prop1}
\end{prop}

\begin{proof} Let $\al\in C(\A)$ and $X\in\A$ with $[X]=\al$. Then
as $\Hom(X,X)$ is a finite-dimensional $\K$-algebra by Assumption
\ref{as3ass}, general properties of abelian categories imply
$X\!\cong\!X_1\!\op\!\cdots\!\op\!X_n$, where the $0\not\cong
X_i\in\A$ are {\it indecomposable}, and are unique up to order and
isomorphism. Consider $\bigl\{[X_1],\ldots,[X_n]\bigr\}$ as a subset
of $C(\A)$ {\it with multiplicity}, that is, we remember how many
times each element of $C(\A)$ is repeated in $[X_1],\ldots,[X_n]$.
Then $\bigl\{[X_1],\ldots,[X_n]\bigr\}$ depends only on the
isomorphism class of $X$, that is, on~$[X]\in\fObj^\al_\A(\K)$.

Form the map $[X]\mapsto\bigl\{[X_1],\ldots,[X_n]\bigr\}$ from
$\fObj^\al_\A(\K)$ to the set of finite subsets of $C(\A)$ with
multiplicity. Using \cite[Ass.~7.1(iii)]{Joyc3}, it is not difficult
to see this map is {\it locally constructible}. As $\Oss^\al(\tau)$
is constructible by Definition \ref{as4def5}(ii), it follows that
this map takes {\it only finitely many values} on~$\Oss^\al(\tau)$.

Suppose $\be,\ga\!\in\!C(\A)$ with $\al\!=\!\be\!+\!\ga$, $\tau(\al)
\!=\!\tau(\be)\!=\!\tau(\ga)$ and $\Oss^\be(\tau)\!\ne\!\emptyset\!
\ne\!\Oss^\ga(\tau)$. Pick $[Y]\in\Oss^\be(\tau)$ and $[Z]\in\Oss^\ga
(\tau)$, and set $X=Y\op Z$. Then $X$ is $\tau$-semistable with
$[X]=\al$, so $[X]\in\Oss^\al(\tau)$. Let $Y\!\cong\!X_1\!\op\!
\cdots\!\op\!X_k$ and $Z\cong X_{k+1}\!\op\!\cdots\!\op\!X_n$ with
all $0\not\cong X_i\in\A$ indecomposable. Then $X\!\cong\!X_1\!\op
\!\cdots\!\op\!X_n$ splits $X$ into indecomposables. Hence there
are only finitely many possibilities for $\bigl\{[X_1],\ldots,
[X_n]\bigr\}$, as a subset of $C(\A)$ with multiplicity. But
$\be=[X_1]+\cdots+[X_k]$ and $\ga=[X_{k+1}]+\cdots+[X_n]$, so
we see there are only finitely many possibilities for~$\be,\ga$.
\end{proof}

In \cite{Joyc5} we will need the following notion.

\begin{dfn} Let $(\tau,T,\le)$ and $(\ti\tau,\ti T,\le)$ be
weak stability conditions on an abelian category $\A$, with
the same $K(\A)$. We say $(\ti\tau,\ti T,\le)$ {\it dominates}
$(\tau,T,\le)$ if $\tau(\al)\le\tau(\be)$ implies $\ti\tau(\al)
\le\ti\tau(\be)$ for all~$\al,\be\in C(\A)$.
\label{as4def6}
\end{dfn}

Many examples of this arise through the following construction:
if $(\tau,T,\le)$ is a weak stability condition, $(\ti T,\le)$
a total order, and $\pi:T\ra\ti T$ a map with $t\le t'$ implies
$\pi(t)\le\pi(t')$, then setting $\ti\tau=\pi\ci\tau$ we find
$(\ti\tau,\ti T,\le)$ is a weak stability condition dominating
$(\tau,T,\le)$. The next lemma is elementary.

\begin{lem} Let\/ $(\ti\tau,\ti T,\le)$ dominate $(\tau,T,\le)$ on
$\A$. Then $X$ $\ti\tau$-stable implies $X$ $\tau$-stable implies
$X$ $\tau$-semistable implies $X$ $\ti\tau$-semistable for $X\in\A$.
Also $\A$ $\ti\tau$-artinian implies $\A$ $\tau$-artinian, and if
Assumption \ref{as3ass} holds then $(\ti\tau,\ti T,\le)$ permissible
implies $(\tau,T,\le)$ permissible.
\label{as4lem1}
\end{lem}

\subsection{Stability conditions on quiver representations}
\label{as43}

We give examples of permissible stability conditions for the data
$\A,K(\A),\fF_\A$ of \cite[\S 10]{Joyc3}. Here is a criterion for
weak stability conditions to be permissible.

\begin{prop} If Assumption \ref{as3ass} holds and\/
$\fObj_\A^\al$ is of finite type for all\/ $\al\in C(\A)$ then all
weak stability conditions $(\tau,T,\le)$ on $\A$ are permissible.
\label{as4prop2}
\end{prop}

\begin{proof} Suppose $\cdots\subset A_2\subset A_1 \subset X$ is an
infinite chain of subobjects in $\A$ with $A_{n+1}\ne A_n$ for all
$n$. Set $\al=[X]$ in $C(\A)$. Consider the function
$\fObj_\A^\al(\K)\ra\N$ taking $[Y]\mapsto n$, where $Y\cong
Y_1\op\cdots\op Y_n$ has $n$ indecomposable factors $0\not\cong
Y_1,\ldots,Y_n$. This function is locally constructible, and so
takes only finitely many values on $\fObj_\A^\al(\K)$ as
$\fObj_\A^\al$ is of finite type. Thus it has a maximum value
$n^\al$. However, $Y=(X/A_2)\op(A_2/A_3)\op\cdots\op
(A_{n^\al}/A_{n^\al+1})\op A_{n^\al+1}$ has at least $n^\al+1$
indecomposable factors and $[Y]=[X]=\al$, a contradiction.

Thus there exist no such infinite chains $\cdots\subset A_2\subset
A_1\subset X$, so $\A$ is {\it artinian}, and therefore $\tau$-{\it
artinian} for any $(\tau,T,\le)$, proving Definition
\ref{as4def5}(i). For (ii), as $\Oss^\al(\tau)$ is locally
constructible by Definition \ref{as4def4} and a subset of
$\fObj_\A^\al(\K)$ which is constructible as $\fObj_\A^\al$ is of
finite type, $\Oss^\al(\tau)$ is constructible.
\end{proof}

In \cite[Ex.s 10.5--10.9]{Joyc3} we define data $\A,K(\A),\fF_\A$
satisfying Assumption \ref{as3ass} with $\A=\modKQ$ or $\nilKQ$ for
$Q=(Q_0,Q_1,b,e)$ a {\it quiver}, and $\A=\modKQI$ or $\nilKQI$ for
$(Q,I)$ a {\it quiver with relations}, and $\A=\modA$ for $A$ a {\it
finite-dimensional\/ $\K$-algebra}. For all of these $\fObj_\A^\al$
is of finite type by \cite[Th.~10.11]{Joyc3}, so Proposition
\ref{as4prop2} gives:

\begin{cor} For the data $\A,K(\A),\fF_\A$ defined using
quivers in {\rm\cite[Ex.s 10.5--10.9]{Joyc3}}, all weak stability
conditions $(\tau,T,\le)$ on $\A$ are permissible.
\label{as4cor}
\end{cor}

Stability conditions on categories of quiver representations
were first considered by King \cite{King}, who proved the
existence of coarse moduli schemes of semistable representations.
His definition of stability \cite[Def.~1.1]{King} is not of our
type, though it gives the same notions of (semi)stable object.
Instead, we define stability using {\it slope functions}
following \cite[\S 3]{Ruda}, based on much older ideas on
slope stability for vector bundles and coherent sheaves.

\begin{ex} Let $\K$ be an algebraically closed field, and
$\A,K(\A),\fF_\A$ be as in one of \cite[Ex.s 10.5--10.9]{Joyc3}. In
each case there is an isomorphism $\bdim:K(\A)\ra\Z^{Q_0}$, where
$Q_0$ is the finite set of vertices of a quiver $Q$. If $X\in\A$
then $\bdim[X]\in\N^{Q_0}\subset\Z^{Q_0}$ is the {\it dimension
vector} of $X$, so~$\bdim\,C(\A)=\N^{Q_0}\sm\{0\}$.

Let $c,r:K(\A)\ra\R$ be group homomorphisms with $r(\al)>0$
for all $\al\in C(\A)$. Using $\bdim:K(\A)\ra\Z^{Q_0}$
we see $c,r$ may be uniquely written
\begin{equation*}
c(\al)=\ts\sum_{v\in Q_0}c_v(\bdim\al)(v)
\quad\text{and}\quad
r(\al)=\ts\sum_{v\in Q_0}r_v(\bdim\al)(v),
\end{equation*}
where $c_v,r_v\in\R$ for $v\in Q_0$, and $r_v>0$ for all
$v\in Q_0$. It is common to take $r_v=1$ for all $v$, so that
$r(\al)$ is the {\it total dimension} of $\al$. Define
$\mu:C(\A)\ra\R$ by $\mu(\al)=c(\al)/r(\al)$ for $\al\in C(\A)$.
Then $\mu$ is called a {\it slope function} on $K(\A)$, as
$\mu(\al)$ is the {\it slope} of the vector $\bigl(r(\al),
c(\al)\bigr)$ in $\R^2$. It is easy to verify $(\mu,\R,\le)$
is a {\it stability condition} on $\A$, which is
{\it permissible} by Corollary~\ref{as4cor}.
\label{as4ex1}
\end{ex}

\subsection{(Weak) stability conditions on coherent sheaves}
\label{as44}

Next we define (weak) stability conditions $(\tau,T,\le)$ for the
examples of \cite[\S 9]{Joyc3}, in which $\A=\coh(P)$ is the abelian
category of {\it coherent sheaves} on a projective $\K$-scheme $P$.
Our first example is {\it Gieseker stability}, introduced by
Gieseker \cite{Gies} for vector bundles on algebraic surfaces, and
studied in \cite{HuLe}. We define some total orders $(G_m,\le)$ on
sets of {\it monic polynomials}.

\begin{dfn} Let $m\ge 0$ be an integer, and define
\e
G_m=\bigl\{p(t)=t^d+a_{d-1}t^{d-1}+\cdots+a_0:
\text{$0\le d\le m$, $a_0,\ldots,a_{d-1}\in\R$}\bigr\}.
\label{as4eq6}
\e
That is, $G_m$ is the set of {\it monic real polynomials
$p$ of degree at most\/} $m$. Here `monic' means
{\it with leading coefficient\/}~1.

Define a {\it total order} `$\le$' on $G_m$ by $p\le q$
for $p,q\in G_m$ if either
\begin{itemize}
\setlength{\itemsep}{0pt}
\setlength{\parsep}{0pt}
\item[(a)] $\deg p>\deg q$, or
\item[(b)] $\deg p=\deg q$ and $p(t)\le q(t)$ for all $t\gg 0$.
\end{itemize}
Explicitly, if $p(t)=t^d+a_{d-1}t^{d-1}+\cdots+a_0$ and
$q(t)=t^e+b_{e-1}t^{e-1}+\cdots+b_0$, we have $p\le q$
if either (a) $d>e$, or (b) $d=e$, and either $p=q$ or
for some $k=0,\ldots,d-1$ we have $a_k<b_k$ and $a_l=b_l$
for~$k<l<d$.

Note that (a) and (b) are {\it not\/} related in the way one
might expect. For if $\deg p>\deg q$ as in (a) then $p(t)>q(t)$
for all $t\gg 0$, which is the opposite of $p(t)\le q(t)$ for
all $t\gg 0$ in~(b).
\label{as4def7}
\end{dfn}

We define Gieseker stability conditions on $\coh(P)$,
following~\cite[\S 2]{Ruda}.

\begin{ex} Let $\K$ be an algebraically closed field, $P$
a projective $\K$-scheme of dimension $m$, $\A=\coh(P)$ the abelian
category of {\it coherent sheaves} on $P$, and $K(\A),\fF_\A$ as in
\cite[Ex.~9.1 or Ex.~9.2]{Joyc3}, supposing $P$ {\it smooth\/}
in~\cite[Ex.~9.1]{Joyc3}.

Let $E$ be an ample line bundle (invertible sheaf) on $P$. For
$X\in\coh(P)$, following \cite[\S 1.2]{HuLe} define the {\it Hilbert
polynomial\/} $p_X$ computed using $E$ by
\e
p_X(n)=\ts\sum_{i=0}^m(-1)^i\dim_\K H^k\bigl(P,X\ot E^n\bigr)
\quad\text{for $n\in\Z$,}
\label{as4eq7}
\e
where $H^*(P,\cdot)$ is {\it sheaf cohomology} on $P$. Then
\e
p_X(n)=\ts\sum_{i=0}^mb_in^i/i! \quad\text{for
$b_0,\ldots,b_m\in\Z$,}
\label{as4eq8}
\e
by \cite[p.~10]{HuLe}. So $p_X(t)$ is a polynomial with rational
coefficients, written $p_X(t)\in\Q[t]$, with degree no more than
$m$. It depends only on the class $[X]$ in $K(\A)$, so that
$p_X=\Pi([X])$ for a unique group homomorphism~$\Pi:K(\A)\ra\Q[t]$.

If $X\not\cong 0$ then the degree of $p_X$ is the dimension of the
support of $X$, and the leading coefficient of $p_X$ is positive.
Hence by~\eq{as4eq8},
\begin{equation*}
\Pi\bigl(C(\A)\bigr)\subseteq\bigl\{p(t)=
\ts\sum_{i=0}^kb_it^i/i!:
\text{$0\le k\le m$, $b_0,\ldots,b_k\in\Z$, $b_k>0$}\bigr\}.
\end{equation*}
Let $(G_m,\le)$ be as in Definition \ref{as4def7}, and define
$\ga:C(\A)\ra G_m$ by
\begin{equation*}
\ga(\al)=\sum_{i=0}^k\frac{k!\,b_i}{i!\,b_k}\,t^i
\quad\text{when $\Pi(\al)=\sum_{i=0}^k\frac{b_i}{i!}\,t^i$, $b_k>0$.}
\end{equation*}
That is, $\ga(\al)$ is $\Pi(\al)$ divided by the leading coefficient
$b^k/k!$ to make it {\it monic}, as in \eq{as4eq6}. So $\ga$ does
map~$C(\A)\ra G_m$.

By Rudakov \cite[Lem.~2.5]{Ruda}, $(\ga,G_m,\le)$ is a {\it
stability condition}, in the sense of Definition \ref{as4def1}. It
is {\it permissible} by Theorem \ref{as4thm4} below. By
construction, $\ga$-(semi)stability coincides with the definition of
Gieseker (semi)stability in \cite[Def.~1.2.4]{HuLe}, which refers to
it just as (semi)stability. Note that the restriction in
\cite[Def.~1.2.4]{HuLe} that (semi)stable sheaves must be {\it pure}
follows automatically from Definitions \ref{as4def2}
and~\ref{as4def7}(a).
\label{as4ex2}
\end{ex}

Huybrechts and Lehn also define $\mu$-({\it semi}){\it stability}
of coherent sheaves \cite[Def.~1.2.12]{HuLe}. We can express this
as a {\it weak\/} stability condition $(\mu,M_m,\le)$ on $\coh(P)$,
by {\it truncating} $p_X(t)$ at the second term.

\begin{ex} In the situation of Example \ref{as4ex2}, define
\e
M_m=\bigl\{p(t)=t^d+a_{d-1}t^{d-1}:\text{$0\le d\le m$,
$a_{d-1}\in\R$}\bigr\}\subseteq G_m,
\label{as4eq9}
\e
and restrict the total order $\le$ on $G_m$ to $M_m$. Define
$\pi_M:G_m\ra M_m$ by $\pi_M:t^d+a_{d-1}t^{d-1}+\cdots+a_0\mapsto
t^d+a_{d-1}t^{d-1}$. Define $\mu:C(\coh(P))\ra M_m$ by
$\mu=\pi_M\ci\ga$. Then $p\le q$ implies $\pi(p)\le\pi(q)$ for all
$p,q\in G_m$, so as $(\ga,G_m,\le)$ is a stability condition on
$\coh(P)$, the remark after Definition \ref{as4def6} shows
$(\mu,M_m,\le)$ is a {\it weak stability condition} on $\coh(P)$,
which dominates $(\ga,G_m,\le)$. It is {\it permissible} by Theorem
\ref{as4thm4} below.

It is easy to show that $X\in\coh(P)$ is $\mu$-(semi)stable in
our sense if and only if $X$ is {\it pure} and $\mu$-(semi)stable
in the sense of \cite[Def.~1.2.12]{HuLe}. Note that Huybrechts
and Lehn do not require $\mu$-semistable sheaves $X$ to be pure,
only that torsion subsheaves of $X$ have codimension at least two.
\label{as4ex3}
\end{ex}

When $m=\dim P\ge 2$ we can find $\al,\be,\ga\in C(\A)$ with
$\be=\al+\ga$ and $\Pi(\al)=t^2$, $\Pi(\be)=t^2+1$, $\Pi(\ga)=1$
for $\Pi$ as in Example \ref{as4ex2}. Then $\mu(\al)=\mu(\be)=t^2$
but $\mu(\ga)=1$, so that $\mu(\al)=\mu(\be)<\mu(\ga)$, which
violates the seesaw inequality. Therefore $(\mu,M_m,\le)$ is
{\it not\/} a stability condition.

We defined $(\mu,M_m,\le)$ by truncating Hilbert polynomials
$p_X(t)$ at the second term. Truncating after any number of terms
also gives a weak stability condition. In particular, we may
truncate after one term, which is related to {\it pure sheaves}
\cite[Def.~1.1.2]{HuLe} and the {\it torsion
filtration}~\cite[Def.~1.1.4]{HuLe}.

\begin{ex} In the situation of Examples \ref{as4ex2} and
\ref{as4ex3}, define
\begin{equation*}
D_m=\bigl\{p(t)=t^d:0\le d\le m\bigr\}\subseteq M_m\subseteq G_m,
\end{equation*}
and restrict $\le$ on $G_m$ to $D_m$, so that $t^d\le t^e$ if and
only if $d\ge e$. Define $\pi_D:G_m\ra D_m$ by $\pi_D:t^d+a_{d-1}
t^{d-1}+\cdots+a_0\mapsto t^d$. Define $\de:C(\coh(P))\ra M_m$ by
$\de=\pi_D\ci\ga$. Then $(\de,D_m,\le)$ is a {\it weak stability
condition} on $\coh(P)$ as in Example \ref{as4ex3}, which dominates
$(\ga,G_m,\le)$ and $(\mu,M_m,\le)$. It is easy to show $X\in\coh(P)$
is $\de$-semistable if and only if $X$ is {\it pure}. Note that
$\de([X])=t^{\dim X}$ for $X\in\coh(P)$, so $(\de,D_m,\le)$ is
independent of choice of ample line bundle~$E$.
\label{as4ex4}
\end{ex}

We show below that $\coh(P)$ is $\de$-{\it artinian}. Thus Theorem
\ref{as4thm1} shows every $X\in\coh(P)$ has a unique filtration
$0=A_0\subset\cdots\subset A_n=X$ with $S_k=A_k/A_{k-1}$ pure of
strictly increasing dimension. This is the {\it torsion filtration}
of $X$, with repeated terms omitted. Again, $(\de,D_m,\le)$ is not
a stability condition for $m\ge 1$. These examples suggest weak
stability conditions are a useful idea.

\begin{lem} $\coh(P)$ is $\de$-artinian in Example~\ref{as4ex4}.
\label{as4lem2}
\end{lem}

\begin{proof} Suppose for a contradiction that there exists
$\cdots\!\subset\!A_2\!\subset A_1\!\subset X$ in $\coh(P)$ with
$A_{n+1}\ne A_n$ and $\de([A_{n+1}])\ge\de([A_n/A_{n+1}])$ for
all $n$. Then $\de([A_{n+1}])\ge\de([A_n])$, so as $\de([A_n])
=t^{\deg\Pi([A_n])}$ we see $(\deg\Pi([A_n]))_{n\ge 1}$ is a
decreasing sequence of nonnegative integers. Thus $\deg\Pi([A_n])=d$
for some $N$ and all $n\ge N$. For $n\ge N$ we have $\Pi([A_n])=
a_{n,d}t^d/d!\!+\!\cdots\!+\!a_{n,0}$, and $\de([A_{n+1}])\ge
\de([A_n/A_{n+1}])$ implies $\de([A_n/A_{n+1}])$ also has degree
$d$, which forces $a_{n+1,d}<a_{n,d}$. Hence $(a_{n,d})_{n\ge N}$ is
a strictly decreasing sequence of positive integers, a contradiction.
\end{proof}

\begin{thm} $(\ga,G_m,\le)$ and\/ $(\mu,M_m,\le)$ above are permissible.
\label{as4thm4}
\end{thm}

\begin{proof} As $(\de,D_m,\le)$ dominates $(\ga,G_m,\le)$ and
$(\mu,M_m,\le)$, Lemmas \ref{as4lem1} and \ref{as4lem2} imply
$\coh(P)$ is $\ga$- and $\mu$-artinian, proving Definition
\ref{as4def5}(i). For (ii), as $\fObj_{\rm ss}^\al(\ga), \fObj_{\rm
ss}^\al(\mu)$ are locally constructible by Definition \ref{as4def4},
they are constructible if they are contained in a constructible set.
This is equivalent to the families of $\ga$- and $\mu$-semistable
sheaves in class $\al$ being {\it bounded\/} in the sense of
\cite[Def.~1.7.5]{HuLe}. This is proved for $\K$ of characteristic
zero by Huybrechts and Lehn \cite[Th.~3.3.7]{HuLe}, and for
arbitrary characteristic by Langer~\cite[Th.~4.2]{Lang}.
\end{proof}

Note that $(\de,D_m,\le)$ in Example \ref{as4ex4} is {\it not\/}
permissible when $m\!=\!\dim P\!\ge\!1$, as the pure sheaves in a
class $\al$ of nonzero degree are not bounded.

\section{Identities relating the $\dss,\dsi,\dst,\dssb,\dsib,
\dstb(*,\tau)$}
\label{as5}

Here and in \S\ref{as6} we will derive {\it universal identities}
relating the six families of constructible functions $\dss,\dsi,
\dst,\dssb,\dsib,\dstb(*,\tau)$. This section works using
constructible function techniques, mostly involving computing Euler
characteristics of pieces of moduli spaces. Section \ref{as6} then
uses combinatorial methods to invert the identities of this section.
As we are working with constructible functions, we assume $\K$ has
{\it characteristic zero} here and in~\S\ref{as6}.

In \S\ref{as51}--\S\ref{as52}, which relate configurations to best
configurations and semistables to semistable-indecomposables, we
work with a permissible {\it weak\/} stability condition
$(\tau,T,\le)$. But in \S\ref{as53}--\S\ref{as54}, which relate
$\tau$-stability and $\tau$-semistability, we take $(\tau,T,\le)$ to
be a stability condition, so that Theorem \ref{as4thm2} applies. Our
results show that to express $\dss^\al(\tau)$ in terms of
$\dst^\be(\tau)$ and vice versa, we have to use configuration moduli
stacks $\fM(I,\pr)_\A$ for all finite posets $(I,\pr)$. This is some
justification for the work of developing the configurations
formalism.

\subsection{Counting best improvements}
\label{as51}

Our first theorem says, in effect, that the family of all {\it best
improvements\/} of an $(I,\tl)$-configuration $(\si,\io,\pi)$ has
Euler characteristic~1.

\begin{thm} Let Assumption \ref{as3ass} hold and\/ $(\tau,T,\le)$ be
a permissible weak stability condition on $\A$. Suppose
$(I,\tl,\ka)$ is $\A$-data, as in Definition \ref{as3def5}, and\/
$X\in\A$ with\/ $[X]=\ka(I)$. Then
\ea
\sum_{\substack{\text{p.o.s $\pr$ on $I$:}\\
\text{$\tl$ dominates $\pr$}}}
\CF^\stk\bigl(Q(I,\pr,\tl)\bigr)\dssb(X,I,\pr,\ka,\tau)&=
\dss(X,I,\tl,\ka,\tau),
\label{as5eq1}\\[3pt]
\sum_{\substack{\text{p.o.s $\pr$ on $I$:}\\
\text{$\tl$ dominates $\pr$}}}
\CF^\stk\bigl(Q(I,\pr,\tl)\bigr)\dsib(X,I,\pr,\ka,\tau)&=
\dsi(X,I,\tl,\ka,\tau),
\label{as5eq2}\\[3pt]
\text{and}
\sum_{\substack{\text{p.o.s $\pr$ on $I$:}\\
\text{$\tl$ dominates $\pr$}}}
\CF^\stk\bigl(Q(I,\pr,\tl)\bigr)\dstb(X,I,\pr,\ka,\tau)&=
\dst(X,I,\tl,\ka,\tau).
\label{as5eq3}
\ea
\label{as5thm1}
\end{thm}

\begin{proof} Define $S=\{(i,j)\in I\t I:i\ne j$ and
$i\tl j\}$, and let $s=\md{S}$. Choose some arbitrary
{\it total order} $\le$ on $S$. Define a finite type
algebraic $\K$-space $\fG$ by
\begin{equation*}
\fG=\coprod_{\substack{\text{p.o.s $\pr$ on $I$:}\\
\text{$\tl$ dominates $\pr$}}}\fM(X,I,\pr,\ka)_\A.
\end{equation*}
Define a 1-morphism $\phi_r:\fG\ra\fG$ for $r=1,\ldots,s$ by
\begin{equation*}
\phi_r\vert_{\fM(X,I,\pr,\ka)_\A}=
\begin{cases} \id:\fM(X,I,\pr,\ka)_\A\ra\fM(X,I,\pr,\ka)_\A, & m\ne r, \\
Q(I,\pr,\ls):\fM(X,I,\pr,\ka)_\A\ra\fM(X,I,\ls,\ka)_\A, & m=r,
\end{cases}
\end{equation*}
if $\tl$ dominates $\pr$ by $m$ steps, where $\ls$ is defined as
follows: let $(i,j)\in S$ be $\le$-least such that (a) $i\npr j$,
(b) if $i\ne k\in I$ with $i\pr k$ then $j\pr k$, and (c) if
$j\ne k\in I$ with $k\pr i$, then $k\pr j$. Then set $a\ls b$ if
either $a\pr b$ or~$a=i,b=j$.

By \cite[Lem.~6.4]{Joyc3} and (a)--(c), $\ls$ is a partial order and
dominates $\pr$ by one step, and $(i,j)\in S$ gives $i\tl j$, so
that $\tl$ dominates $\ls$. Conversely, if $\tl$ dominates $\ls$
dominates $\pr$ by one step then it arises in this way for a unique
$(i,j)\in S$. As $r\ge 1$ there is at least one $\ls$ with $\tl$
dominates $\ls$ dominates $\pr$ by one step, by
\cite[Prop.~6.5]{Joyc3}. Thus the set of $(i,j)\in S$ which from
which we choose the $\le$-least element is nonempty, and $\phi_r$ is
well-defined.

If $\tl$ dominates $\pr$ by $m$ steps then $\phi_r$ fixes
$\pr$ if $m\ne r$, and takes $\pr$ to $\ls$ if $m=r$, where
$\tl$ dominates $\ls$ by $r-1$ steps. So by induction
$\phi_r\ci\phi_{r+1}\ci\cdots\ci\phi_s$ takes each $\pr$
to some $\ls$, where $\tl$ dominates $\ls$ by less than
$r$ steps. When $r=1$ we have $\ls=\tl$, as $\tl$ dominates
$\ls$ by 0 steps. It follows easily that
\e
\phi_1\ci\phi_2\ci\cdots\ci\phi_s\vert_{\fM(X,I,\pr,\ka)_\A}=
Q(I,\pr,\tl).
\label{as5eq4}
\e

Define
\e
{\mathcal C}_s=\coprod_{\substack{\text{p.o.s $\pr$ on $I$:}\\
\text{$\tl$ dominates $\pr$}}}\Mssb(X,I,\pr,\ka,\tau)_\A\subseteq\fG(\K).
\label{as5eq5}
\e
Then ${\mathcal C}_s$ is a constructible set in $\fG$ by Theorem
\ref{as4thm3}. For $r=s,s-1,\ldots,1$ define ${\mathcal
C}_{r-1}=(\phi_r)_*({\mathcal C}_r)$. As $\phi_r$ is a 1-morphism,
Proposition \ref{as2prop} shows that ${\mathcal C}_r$ is also
constructible for $r=s,s\!-\!1,\ldots,0$. Equation \eq{as5eq4} gives
\e
{\mathcal C}_0=\!\!\!\!\!
\coprod_{\substack{\text{p.o.s $\pr$ on $I$:}\\
\text{$\tl$ dominates $\pr$}}}
\!\!\!\!\!\!\!
Q(I,\pr,\tl)_*\bigl(\Mssb(X,I,
\pr,\ka,\tau)_\A\bigr)\!=\!\Mss(X,I,\tl,\ka,\tau)_\A,
\label{as5eq6}
\e
as every $(I,\tl,\ka)$-configuration admits a best improvement
by~\cite[Lem.~6.2]{Joyc3}.

Suppose $[(\si,\io,\pi)]\in{\mathcal C}_{r-1}$ for $r\le s$,
with $(\si,\io,\pi)$ an $(I,\ls)$-configuration. We shall
determine $(\phi_r)_*^{-1}\bigl([(\si,\io,\pi)]\bigr)$ in
${\mathcal C}_r$. If $(\si,\io,\pi)$ is not best
then by Theorem \ref{as3thm1} there are $i\ne j\in I$ with
$i\ls j$ but there exists no $k\in I$ with $i\ne k\ne j$ and
$i\ls k\ls j$, such that \eq{as3eq1} is split.

Now $i\tl j$ as $\tl$ dominates $\ls$, so $(i,j)\in S$. Let
$(i,j)$ be greatest in the total order $\le$ on $S$ satisfying
these conditions. Define $\pr$ by $a\pr b$ if $a\ls b$ and
$a\ne i$ or $b\ne j$. Then $\pr$ is a partial order on $I$
and $\ls$ dominates $\pr$ by one step. Furthermore, Theorem
\ref{as3thm1} and the construction of the ${\mathcal C}_r,\phi_r$
imply that $(\phi_r)_*^{-1}\bigl([(\si,\io,\pi)]\bigr)$ is
exactly the set of isomorphism classes $[(\si',\io',\pi')]$
of $(I,\pr)$-improvements $(\si',\io',\pi')$ of $(\si,\io,\pi)$,
which are in 1-1 correspondence with~$\Hom\bigl(\si(\{j\}),
\si(\{i\})\bigr)$.

Regard $\Hom\bigl(\si(\{j\}),\si(\{i\})\bigr)\cong\K^l$ as an affine
$\K$-variety. Using \cite[\S 6.2 \& Ass.~7.1(iv)]{Joyc3} one can
construct a $\K$-subvariety $V$ of $\fG$ isomorphic to $\K^l$, such
that $V(\K)=(\phi_r)_*^{-1}\bigl([(\si,\io,\pi)]\bigr)$. Hence
$\chi^\na\bigl((\phi_r)_*^{-1}\bigl([(\si,\io,\pi)]\bigr)\bigr)
\!=\!\chi(V)\!=\!\chi(\K^l)\!=\!1$, by \eq{as2eq2}. If
$(\si,\io,\pi)$ is best then $(\phi_r)_*^{-1}\bigl([(\si,\io,
\pi)]\bigr)\!=\!\bigl\{[(\si,\io,\pi)]\bigr\}$, so again~$\chi^\na
\bigl((\phi_r)_*^{-1}\bigl([(\si,\io,\pi)]\bigr)\bigr)=1$.

Write $\de_{{\mathcal C}_r}$ for the characteristic function of
${\mathcal C}_r$. Then $\de_{{\mathcal C}_r}$ is a constructible
function, as ${\mathcal C}_r$ is a constructible set. Since
$\chi^\na\bigl((\phi_r)_*^{-1}(x)\bigr)=1$ for all $x\in{\mathcal
C}_r$ and $m_{\phi_r}\equiv 1$ in Definition \ref{as2def6} as
$\fM(X,I,\pr,\ka)_\A$ is an algebraic $\K$-space with trivial
stabilizer groups, we see that $\CF^\stk(\phi_r)\de_{{\mathcal
C}_r}=\de_{{\mathcal C}_{r-1}}$ for all $r$. Hence
$\CF^\stk(\phi_1\ci\cdots\ci\phi_s)\de_{{\mathcal C}_s}
=\de_{{\mathcal C}_0}$, by \eq{as2eq3}. Equation \eq{as5eq1} then
follows from \eq{as5eq4}--\eq{as5eq6}. To prove \eq{as5eq2} and
\eq{as5eq3} we proceed in the same way, but define ${\mathcal C}_s$
in \eq{as5eq5} using $\Msib,\Mstb(X,I,\pr,\ka,\tau)_\A$ rather
than~$\Mssb(X,I,\pr,\ka,\tau)_\A$.
\end{proof}

Here are analogues of \eq{as5eq1}--\eq{as5eq3} for
$\fM(I,\pr,\ka)_\A$ rather than~$\fM(X,I,\pr,\ka)_\A$.

\begin{thm} Let Assumption \ref{as3ass} hold, $(\tau,T,\le)$ be
a permissible weak stability condition on $\A$, and\/ $(I,\tl,\ka)$
be $\A$-data. Then
\ea
\sum_{\substack{\text{p.o.s $\pr$ on $I$:}\\
\text{$\tl$ dominates $\pr$}}}
\CF^\stk\bigl(Q(I,\pr,\tl)\bigr)\dssb(I,\pr,\ka,\tau)&=
\dss(I,\tl,\ka,\tau),
\label{as5eq7}\\[3pt]
\sum_{\substack{\text{p.o.s $\pr$ on $I$:}\\
\text{$\tl$ dominates $\pr$}}}
\CF^\stk\bigl(Q(I,\pr,\tl)\bigr)\dsib(I,\pr,\ka,\tau)&=
\dsi(I,\tl,\ka,\tau),
\label{as5eq8}\\[3pt]
\text{and}
\sum_{\substack{\text{p.o.s $\pr$ on $I$:}\\
\text{$\tl$ dominates $\pr$}}}
\CF^\stk\bigl(Q(I,\pr,\tl)\bigr)\dstb(I,\pr,\ka,\tau)&=
\dst(I,\tl,\ka,\tau).
\label{as5eq9}
\ea
\label{as5thm2}
\end{thm}

\begin{proof} Let $X\in\A$ with $[X]=\ka(I)$ in $K(\A)$, and
$\pr$ be a partial order on $I$ dominated by $\tl$. Consider
the Cartesian square
\begin{equation*}
\xymatrix@C=60pt@R=8pt{
\fM(X,I,\pr,\ka)_\A \ar[r]_{Q(I,\pr,\tl)} \ar[d]^{\,\Pi_X} &
\fM(X,I,\tl,\ka)_\A \ar[d]_{\Pi_X\,} \\
\fM(I,\pr,\ka)_\A \ar[r]^{Q(I,\pr,\tl)} & \fM(I,\tl,\ka)_\A.
}
\end{equation*}
The $Q(I,\pr,\tl)$ are representable by Theorem \ref{as3thm2}(a),
and the $\Pi_X$ of finite type by Definition \ref{as3def6}. Thus
Theorem \ref{as2thm1} shows the following commutes:
\e
\begin{gathered}
\xymatrix@C=70pt@R=15pt{
\CF\bigl(\fM(X,I,\pr,\ka)_\A\bigr)
\ar[r]_{\CF^\stk(Q(I,\pr,\tl))} &
\CF\bigl(\fM(X,I,\tl,\ka)_\A\bigr) \\
\CF\bigl(\fM(I,\pr,\ka)_\A\bigr)
\ar[r]^{\CF^\stk(Q(I,\pr,\tl))} \ar[u]_{\Pi_X^*}
& \CF\bigl(\fM(I,\tl,\ka)_\A\bigr). \ar[u]^{\Pi_X^*}
}
\end{gathered}
\label{as5eq10}
\e

Using \eq{as5eq1}, commutativity of \eq{as5eq10},
$\Pi_X^*\bigl(\dss(I,\pr,\ka,\tau)\bigr)=\dss(X,I,\pr,\ka,\tau)$ and
$\Pi_X^*\bigl(\dss(I,\tl,\ka,\tau)\bigr)=\dss(X,I,\tl,\ka,\tau)$
shows that
\begin{align*}
\Pi_X^*
\raisebox{-10pt}{\begin{Large}$\displaystyle\biggl[$\end{Large}}
\sum_{\substack{\text{p.o.s $\pr$ on $I$:}\\
\text{$\tl$ dominates $\pr$}}}
\CF^\stk\bigl(Q(I,\pr,\tl)\bigr)\dssb(I,\pr,\ka,\tau)
\raisebox{-10pt}{\begin{Large}$\displaystyle\biggr]$\end{Large}}&=\\
\sum_{\substack{\text{p.o.s $\pr$ on $I$:}\\
\text{$\tl$ dominates $\pr$}}}
\CF^\stk\bigl(Q(I,\pr,\tl)\bigr)\ci\Pi_X^*
\bigl(\dssb(I,\pr,\ka,\tau)\bigr)&=\\
\sum_{\substack{\text{p.o.s $\pr$ on $I$:}\\
\text{$\tl$ dominates $\pr$}}}
\CF^\stk\bigl(Q(I,\pr,\tl)\bigr)\dssb(X,I,\pr,\ka,\tau)
&\begin{aligned}[t]&=\dss(X,I,\tl,\ka,\tau)\\
&=\Pi_X^*\bigl(\dss(I,\tl,\ka,\tau)\bigr).
\end{aligned}
\end{align*}
This implies that \eq{as5eq7} holds at all $[(\si,\io,\pi)]$ in
$\M(I,\tl,\ka)_\A$ with $\si(I)=X$. Since this is true for all
$X\in\A$ with $[X]=\ka(I)$, we have proved \eq{as5eq7}. Equations
\eq{as5eq8} and \eq{as5eq9} follow from \eq{as5eq2} and \eq{as5eq3}
in the same way.
\end{proof}

\subsection{Relating semistables and semistable-indecomposables}
\label{as52}

Next we shall write $\dss^\al(\tau)$ and $\dss(K,\tl,\mu,\tau)$ in
terms of the $\dsi(I,\pr,\ka,\tau)$.

\begin{thm} Let Assumption \ref{as3ass} hold, $(\tau,T,\le)$ be a
permissible weak stability condition on $\A$, and\/ $\al\in C(\A)$.
Then
\e
\sum_{n=1}^\iy\frac{1}{n!}\cdot
\sum_{\substack{\text{$\ka:\{1,\ldots,n\}\ra C(\A):$}\\
\text{$\ka(\{1,\ldots,n\})\!=\!\al$,
$\tau\!\ci\!\ka\!\equiv\!\tau(\al)$}}}\!\!\!\!\!\!\!\!\!\!\!\!\!\!\!
\begin{aligned}[t]
\CF^\stk\bigl(\bs\si(\{1,\ldots,n\})\bigr)\dsi(\{1,\ldots,n\},
\bu,\ka,\tau)&\\
=\dss^\al(\tau)&,
\end{aligned}
\label{as5eq11}
\e
where $\bu$ is the partial order on $\{1,\ldots,n\}$ with\/ $i\bu j$
if and only if\/ $i=j$. Only finitely many functions
$\dsi(\{1,\ldots,n\},\bu,\ka,\tau)$ in this sum are nonzero.
\label{as5thm3}
\end{thm}

\begin{proof} Suppose $n,\ka$ are as in \eq{as5eq11} with
$\dsi(\{1,\ldots,n\},\bu,\ka,\tau)\ne 0$. If $n=1$ the only
possibility is $\ka(1)=\al$, so let $n>1$. Pick $1\le i<n$, and set
$\be=\ka(\{1,\ldots,i\})$ and $\ga=\ka(\{i+1,\ldots,n\})$. Then
$\be,\ga\in C(\A)$ with $\al=\be+\ga$, and $\tau\!\ci\!\ka\!
\equiv\!\tau(\al)$ implies $\tau(\al)=\tau(\be)=\tau(\ga)$, and
$\dsi(\{1,\ldots,n\},\bu,\ka,\tau)\ne 0$ implies that
$\Oss^\be(\tau)\ne\emptyset\ne\Oss^\ga(\tau)$. Hence there are only
finitely many possibilities for $\be,\ga$, by Proposition
\ref{as4prop1}, and it quickly follows that there are only finitely
many nonzero terms in~\eq{as5eq11}.

Fix $0\not\cong X\in\A$ with $[X]=\al\in C(\A)$. Let the
pairwise-nonisomorphic indecomposable factors of $X$ be
$S_1,\ldots,S_k$, with multiplicities $m_1,\ldots,m_k\ge 1$, so that
$X\cong\bigop_{a=1}^k\bigop^{m_a}S_a$. It is easy to see that $X$ is
$\tau$-semistable if and only if each $S_a$ is also
$\tau$-semistable with $\tau([S_a])=\tau(\al)$.

Let $[(\si,\io,\pi)]\in \Msi(\{1,\ldots,n\},\bu,\ka,\tau)_\A$ with
$\si(\{1,\ldots,n\})\cong X$, for $n,\ka$ as in \eq{as5eq11}. Then
by definition and properties of configurations, $\si(\{i\})$ is
$\tau$-semistable and indecomposable for all $i=1,\ldots,n$ with
$X\cong\bigop_{i=1}^n\si(\{i\})$. So by definition of the $S_a,m_a$,
there must exist a unique, surjective map
$\phi:\{1,\ldots,n\}\ra\{1,\ldots,k\}$ with
$\md{\phi^{-1}(\{a\})}=m_a$, such that $\si(\{i\})\cong S_{\phi(i)}$
for all $i=1,\ldots,n$. This forces $n=\sum_{a=1}^km_a$. It also
implies each $S_a$ is also $\tau$-semistable with
$\tau([S_a])=\tau(\al)$, so $X$ is $\tau$-semistable from above.

Thus, if $X$ is {\it not\/} $\tau$-semistable, there exist no such
$n,\ka$ and $[(\si,\io,\pi)]$, so both sides of \eq{as5eq11} are
zero at $[X]$. Suppose $X$ is $\tau$-semistable, and consider the
set of possible choices $n,\ka$ and $[(\si,\io,\pi)]$ with
$\si(\{1,\ldots,n\})\cong X$. From above $n=\sum_{a=1}^km_a$, and it
is easy to show the possible choices of $\ka$ and $[(\si,\io,\pi)]$
are in 1-1 correspondence with maps
$\phi:\{1,\ldots,n\}\ra\{1,\ldots,k\}$ with
$\md{\phi^{-1}(\{a\})}=m_a$ for all $a$. There are exactly
$n!/m_1!\cdots m_k!$ such maps $\phi$, and in each case we have
$\Aut(\si,\io,\pi)\cong\bigot_{a=1}^k\Aut(S_a)^{m_a}$. So by
definition of $\CF^\stk(\bs\si(\{1,\ldots,n\}))$, we see that the
left hand side of \eq{as5eq11} at $[X]$ is
\e
\frac{1}{n!}\,\cdot\,\frac{n!}{m_1!\cdots m_k!}\,\cdot\,
\chi\left(\frac{\Aut\bigl(\bigop_{a=1}^k\bigop^{m_a}S_a\bigr)}{
\prod_{a=1}^k\Aut(S_a)^{m_a}}\right).
\label{as5eq12}
\e

Using elementary facts about finite-dimensional algebras taken from
Benson \cite[\S 1]{Bens} applied to the $\K$-algebras $\End(S_a)$
and $\End(X)$, we find that
\begin{equation*}
\Aut(S_a)\cong\K^\t\lt J_a\;\>\text{and}\;\>
\Aut\bigl(\ts\bigop_{a=1}^k\bigop^{m_a}S_a\bigr)\cong
\bigl(\ts\prod_{a=1}^k\GL(m_a,\K)\bigr)\lt J_X,
\end{equation*}
where $J_a$ and $J_X$ are the {\it Jacobson radicals} of $\End(S_a)$
and $\End(X)$, which are nilpotent $\K$-groups isomorphic as
$\K$-varieties to finite-dimensional vector spaces
$\K^{l_a},\K^{l_X}$. Using these isomorphisms we construct a natural
fibration
\begin{equation*}
\ts\Pi:\Aut\bigl(\bigop_{a=1}^k\bigop^{m_a}S_a\bigr)/
\prod_{a=1}^k\Aut(S_a)^{m_a}\longra\prod_{a=1}^k
\GL(m_a,\K)/(\K^\t)^{m_a},
\end{equation*}
where $(\K^\t)^{m_a}$ is the maximal torus of diagonal matrices
in~$\GL(m_a,\K)$.

The fibre of $\Pi$ is the quotient of nilpotent groups
$J_X/\bigl(\prod_{a=1}^kJ_a^{m_a}\bigr)$, which is isomorphic as a
$\K$-variety to $\K^{l_X-m_1l_1-\cdots-m_kl_k}$. Therefore every
fibre of $\Pi$ has Euler characteristic 1, so by properties of
$\chi$ we have
\e
\chi\left(\frac{\Aut\bigl(\bigop_{a=1}^k\bigop^{m_a}S_a\bigr)}{
\prod_{a=1}^k\Aut(S_a)^{m_a}}\right)= \chi\left(\prod_{a=1}^k
\frac{\GL(m_a,\K)}{(\K^\t)^{m_a}}\right)=\prod_{a=1}^km_a!.
\label{as5eq13}
\e
Combining \eq{as5eq12} and \eq{as5eq13} shows the left hand side of
\eq{as5eq11} at $[X]$ is 1, the same as the right hand side. This
proves \eq{as5eq11}, and the theorem.
\end{proof}

\begin{thm} Let Assumption \ref{as3ass} hold, and\/
$(\tau,T,\le)$ be a permissible weak stability condition on $\A$.
Then for all\/ $\A$-data $(K,\tl,\mu)$ we have
\e
\begin{aligned}
\sum_{\substack{\text{iso. classes}\\ \text{of finite }\\
\text{sets $I$}}} \frac{1}{\md{I}!}\cdot
\sum_{\substack{
\text{$\ka\!:\!I\!\ra\!C(\A)$, surjective $\phi\!:\!I\!\ra\!K$:}\\
\text{$\ka(\phi^{-1}(k))=\mu(k)$ for $k\in K$,}\\
\text{$\tau\ci\mu\ci\phi\equiv\tau\ci\ka:I\ra T$.}\\
\text{Define $\pr$ on $I$ by $i\pr j$ if $i=j$}\\
\text{or $\phi(i)\ne\phi(j)$ and $\phi(i)\tl\phi(j)$}}}
\!\!\!\!\!\!\!
\begin{aligned}[t]
\CF^\stk\bigl(Q(I,\pr,K,\tl,\phi)\bigr)&\\
\dsi(I,\pr,\ka,\tau)&=\\
\dss(K,\tl,\mu&,\tau).
\end{aligned}
\end{aligned}
\label{as5eq14}
\e
Only finitely many functions $\dsi(I,\pr,\ka,\tau)$ in this sum are
nonzero.
\label{as5thm4}
\end{thm}

\begin{proof} First we prove only finitely many $\dsi(I,\pr,\ka,\tau)$
in \eq{as5eq14} are nonzero. Let $I,\ka,\phi$ be as in \eq{as5eq14},
fix $k\in K$, set $I_k=\phi^{-1}(\{k\})$, $\al=\mu(k)$ and
$n=\md{I_k}$, choose a bijection $\imath:\{1,\ldots,n\}\ra I_k$, and
write $\ka'=\ka\ci\imath$. Then $\al,n,\ka'$ are as in \eq{as5eq11},
so Theorem \ref{as5thm3} shows there are only finitely many $n,\ka'$
with $\dsi(\{1,\ldots,n\},\bu,\ka',\tau)\not\equiv 0$. But
$\dsi(I,\pr,\ka,\tau)\not\equiv 0$ in \eq{as5eq14} implies
$\dsi(\{1,\ldots,n\},\bu,\ka',\tau)\not\equiv 0$. So there are
finitely many possibilities for $I_k,\ka\vert_{I_k}$ up to
isomorphism for each $k\in K$, and thus only finitely many
for~$I,\ka,\phi$.

For each $k\in K$, let $I_k$ be a finite set and $\ka_k:I_k\ra
C(\A)$ a map with $\ka_k(I_k)=\mu(k)$ and $\tau\ci\ka_k\equiv
\tau\ci\mu(k)$. Define $I=\coprod_{k\in K}I_k$ and $\phi:I\ra K$ by
$\phi(i)=k$ if $i\in I_k$. Define a partial order $\pr$ on $I$ using
$K,\mu,\phi$ as in \eq{as5eq14}. Now by applying the proof of
\cite[Th.~7.10]{Joyc3} $\md{K}$ times, we can show that the
following commutative diagram of 1-morphisms of stacks is a
Cartesian square:
\e
\begin{gathered}
\xymatrix@C=95pt@R=9pt{ \fM(I,\pr,\ka)_\A \ar[d]^{\,\prod_{k\in
K}S(I,\pr,I_k)} \ar[r]_{Q(I,\pr,K,\tl,\phi)} &
\fM(K,\tl,\mu)_\A \ar[d]_{\prod_{k\in K}\bs\si(\{k\})\,} \\
\prod_{k\in K}\fM(I_k,\bu,\ka_k)_\A \ar[r]^{\prod_{k\in
K}\bs\si(I_k)} & \prod_{k\in K}\fObj_\A^{\mu(k)}.}
\end{gathered}
\label{as5eq15}
\e
Theorem \ref{as3thm2} shows the rows are representable, and the
right 1-morphism finite type. As \eq{as5eq15} is Cartesian the left
1-morphism is finite type. Applying Theorem \ref{as2thm1} to
\eq{as5eq15} and $\prod_{k\in
K}\dsi(I_k,\bu,\ka_k,\tau)\!\in\!\CF\bigl(\prod_{k\in
K}\fM(I_k,\bu,\ka_k)_\A\bigr)$ yields
\e
\begin{aligned}
\CF^\stk\bigl(Q(I,\pr,K,\tl,\phi)\bigr)\dsi(I,\pr,\ka,\tau)&=\\
\CF^\stk\bigl(Q(I,\pr,K,\tl,\phi)\bigr)\bigl[\ts\prod_{k\in K}
S(I,\pr,I_k)^*\dsi(I_k,\bu,\ka_k,\tau)\bigr]&=\\
\ts\prod_{k\in K}\bs\si(\{k\})^*\bigl[\CF^\stk(\bs\si(I_k))
\dsi(I_k,\bu,\ka_k,\tau)\bigr]&.
\end{aligned}
\label{as5eq16}
\e

Now $\dss(K,\tl,\mu,\tau)=\prod_{k\in K}\bs\si(\{k\})^*[
\dss^{\mu(k)}(\tau)]$. Use \eq{as5eq11} with $I_k,\mu(k)$ in place
of $\{1,\ldots,n\},\al$ to substitute for $\dss^{\mu(k)}(\tau)$
here, taking the product in $\CF(\fM(K,\tl,\mu)_\A)$ of $\md{K}$
copies of \eq{as5eq11} pulled back by $\bs\si(\{k\})^*$. Using
\eq{as5eq16} then yields \eq{as5eq14}, {\it except that\/} rather
than summing over isomorphism classes of sets $I$ and maps $\phi$ we
sum over isomorphism classes of sets $I_k$ for $k\in K$ (here the
sum over sets $I_k$ replaces the sum over $n$ in \eq{as5eq11}, with
$\md{I_k}=n$), and instead of the factor $1/\md{I}!$ we
have~$1/\prod_{k\in K}\md{I_k}!$.

The sums over $I,\phi$ and over $I_k$, $k\in K$ are related as
follows: given $I,\phi$ we set $I_k=\phi^{-1}(\{k\})$ for $k\in K$,
and given $I_k$ for $k\in K$ we define $I=\coprod_{k\in K}I_k$ and
$\phi:I\ra K$ by $\phi\vert_{I_k}\equiv k$. But this is not a 1-1
correspondence: fixing $I_k$ for $k\in K$ up to isomorphism forces
$\md{I}=\sum_{k\in K}\md{I_k}$, which fixes $I$ up to isomorphism;
but there are $\md{I}!/\prod_{k\in K}\md{I_k}!$ choices of
$\phi:I\ra K$ with $\md{\phi^{-1}(\{k\})}=\md{I_k}$ for $k\in K$.
This exactly cancels the difference in the combinatorial factors
$1/\md{I}!$ and $\prod_{k\in K}\md{I_k}!$, proving~\eq{as5eq14}.
\end{proof}

\subsection{Counting best $\tau$-stable configurations}
\label{as53}

Now let $X\in\A$ be $\tau$-semistable. If $[(\si,\io,\pi)]
\in\Mst(X,I,\pr,\ka,\tau)_\A$ with $\tau\ci\ka(i)=\tau([X])$ for all
$i\in I$ then $\si(\{i\})$ is $\tau$-stable for all $i\in I$,
and we call $(\si,\io,\pi)$ a $\tau$-{\it stable configuration}.
From Theorem \ref{as4thm2} we find that $\si(\{i\})$ for $i\in I$
are the $\tau$-stable factors of $X$, and up to isomorphism
depend only on $X$. So $\md{I}$ also depends only on $X$. We
shall calculate the Euler characteristic of the family of all
best $\tau$-stable configurations for $X$ up to isomorphism,
the union of $\Mstb(X,I,\pr,\ka,\tau)_\A$ over isomorphism classes
of $\A$-data $(I,\pr,\ka)$ with $\tau\ci\ka\equiv\tau([X])$.
Consider the following situation.

\begin{dfn} Let Assumption \ref{as3ass} hold, $(\tau,T,\le)$
be a permissible stability condition on $\A$, and $X\in\A$ be
$\tau$-semistable. Then Theorem \ref{as4thm2} decomposes $X$
into $\tau$-stable factors with the same $\tau$-value as $X$,
uniquely up to isomorphism and order. Let $X$ have nonisomorphic
$\tau$-stable factors $S_1,\ldots,S_n$ with multiplicities
$l_1,\ldots,l_n>0$.

For any $\tau$-stable $(I,\pr,\ka)$-configuration $(\si,\io,\pi)$
with $\si(I)=X$ and $\tau\ci\ka\equiv\tau([X])$, the $\si(\{i\})$
for $i\in I$ are isomorphic to $S_m$ with multiplicities $l_m$ for
$m=1,\ldots,n$. Thus $\md{I}=\sum_{m=1}^nl_m$. Fix an indexing
set $I$ with $\md{I}=\sum_{m=1}^nl_m$. For $m=1,\ldots,n$ define
$k_m=\dim\Hom(S_m,X)$. Then $\bigop^{k_m}S_m\cong S_m\ot\Hom(S_m,X)
\subset X$, so~$0\le k_m\le l_m$.
\label{as5def}
\end{dfn}

Fix $a\in I$, and set $J=I\sm\{a\}$. Let $\ls$ be a partial
order on $J$, and define $\tl$ on $I$ by $i\tl j$ for
$i,j\in I$ if either $i,j\in J$ and $i\ls j$, or $i=a$.
Define $\phi:I\ra\{1,2\}$ by $\phi(a)=1$ and $\phi(j)=2$
for $j\in J$. Let $(\ti\si,\ti\io,\ti\pi)$ be a
$(\{1,2\},\le)$-configuration with $\ti\si(\{1,2\})=X$ and
$\ti\si(\{1\})$ $\tau$-stable with $\tau\bigl([\ti\si
(\{1\})]\bigr)=\tau\bigl([X]\bigr)$. Then $\ti\si(\{1\})$ is
(isomorphic to) one of the $\tau$-stable factors of $X$.
Define~$Y=\ti\si(\{2\})$.

Choose $\ka:I\ra K(\A)$ such that $(I,\tl,\ka)$ is $\A$-data,
$\tau\ci\ka\equiv\tau([X])$, $\ka(a)=[\ti\si(\{1\})]$, and
$[X]=\ka(I)$. Then $(J,\ls,\ka\vert_J)$ is also $\A$-data,
and $[Y]=\ka(J)$. Define $\mu:\{1,2\}\ra K(\A)$ by $\mu(1)=
\ka(a)$, $\mu(2)=\ka(J)$. Consider the diagram of 1-morphisms
\begin{equation*}
\xymatrix@C=35pt@R=-2pt{
\fM(I,\tl,\ka)_\A
\ar[rr]_{Q(I,\tl,\{1,2\},\le,\phi)} \ar[dd]^{\,S(I,\tl,J)} &&
\fM(\{1,2\},\le,\mu)_\A \ar[dd]_{\bs\si(\{2\})\,} \\
&&& \Spec\K. \ar[ul]_(0.4){(\ti\si,\ti\io,\ti\pi)} \ar[dl]^Y \\
\fM(J,\ls,\ka)_\A
\ar[rr]^{\bs\si(J)} && \fObj_\A^{\ka(J)}
}
\end{equation*}

By \cite[Th.~7.10]{Joyc3}, the left hand side is a Cartesian square.
And as $\ti\si(\{2\})=Y$, the right hand side commutes. Therefore
$S(I,\tl,J)$ induces a 1-isomorphism
\e
\begin{split}
S(I,\tl,J)_*:\,&\fM(I,\tl,\ka)_\A\t_{Q(I,\tl,\{1,2\},\le,\phi),
\fM(\{1,2\},\le,\mu)_\A,(\ti\si,\ti\io,\ti\pi)}\Spec\K
\longra\\
&\fM(J,\ls,\ka)_\A\t_{\bs\si(J),\fObj_\A^{\ka(J)},Y}
\Spec\K=\fM(Y,J,\ls,\ka)_\A.
\end{split}
\label{as5eq17}
\e
But as $\ti\si(\{1,2\})=X$ we have a commutative diagram
\begin{equation*}
\xymatrix@C=75pt@R=15pt{
\fM(I,\tl,\ka)_\A
\ar[r]_(0.6){Q(I,\tl,\{1,2\},\le,\phi)} \ar[dr]_{\bs\si(I)} &
\fM(\{1,2\},\le,\mu)_\A \ar[d]_{\bs\si(\{1,2\})} &
\Spec\K \ar[l]^{(\ti\si,\ti\io,\ti\pi)}
\ar[dl]^X \\
& \fObj_\A^{\ka(I)}.
}
\end{equation*}
Therefore $\bs\si(\{1,2\})$ induces a 1-morphism
\e
\begin{split}
\bs\si(\{1,2\})_*:\,&\fM(I,\tl,\ka)_\A\t_{Q(I,\tl,\{1,2\},\le,\phi),
\fM(\{1,2\},\le,\mu)_\A,(\ti\si,\ti\io,\ti\pi)}\Spec\K
\longra\\
&\fM(I,\tl,\ka)_\A\t_{\bs\si(I),\fObj_\A^{\ka(J)},X}
\Spec\K=\fM(X,I,\tl,\ka)_\A.
\end{split}
\label{as5eq18}
\e

As \eq{as5eq17} is a 1-isomorphism it is invertible, so \eq{as5eq17}
and \eq{as5eq18} give a 1-morphism
\begin{equation*}
\bs\si(\{1,2\})_*\ci S(I,\tl,J)_*^{-1}:\fM(Y,J,\ls,\ka)_\A\longra
\fM(X,I,\tl,\ka)_\A.
\end{equation*}
On the underlying geometric spaces the 1-isomorphism \eq{as5eq17}
gives a bijection, and \eq{as5eq18} an injective map with image
$Q(I,\tl,\{1,2\},\le,\phi)_*^{-1}
\bigl(\{[(\ti\si,\ti\io,\ti\pi)]\}\bigr)$. Thus we have a 1-1
correspondence
\e
\begin{gathered}
\bigl(\bs\si(\{1,2\})_*\ci S(I,\tl,J)_*^{-1}\bigr)_*:\M(Y,J,\ls,\ka)_\A
\longra \\
Q(I,\tl,\{1,2\},\le,\phi)_*^{-1}\bigl(\{[(\ti\si,\ti\io,\ti\pi)]\}\bigr)
\subseteq\M(X,I,\tl,\ka)_\A.
\end{gathered}
\label{as5eq19}
\e
Here is how to understand \eq{as5eq19}: it maps $[(\si',\io',\pi')]
\mapsto[(\si,\io,\pi)]$, for $(\si',\io',\pi')$ a
$(J,\ls,\ka)$-configuration with $\si'(J)=Y$, and $(\si,\io,\pi)$
the $(I,\tl,\ka)$-configuration constructed by {\it substituting}
$(\si',\io',\pi')$ into $(\ti\si,\ti\io,\ti\pi)$ at 2,
using~\cite[Def.~5.7]{Joyc3}.

As $\ti\si(\{1\})$ is $\tau$-stable, \eq{as5eq19} is a 1-1
correspondence between $\tau$-stable configurations in its domain
and range. Hence
\e
\begin{gathered}
\CF^\stk\bigl(\bs\si(\{1,2\})_*\ci S(I,\tl,J)_*^{-1}\bigr)
\dst(Y,J,\ls,\ka,\tau)=\\
\dst(X,I,\tl,\ka,\tau)\cdot\de_{Q^{-1}([(\ti\si,\ti\io,\ti\pi)])},
\end{gathered}
\label{as5eq20}
\e
writing $Q^{-1}([(\ti\si,\ti\io,\ti\pi)])$ as a shorthand for
$Q(I,\tl,\{1,2\},\le,\phi)^{-1}\bigl([(\ti\si,\ti\io,\ti\pi)]\bigr)$,
and $\de_{Q^{-1}([(\ti\si,\ti\io,\ti\pi)])}$ for its characteristic
function.

We now apply \eq{as5eq3} to rewrite $\dst(X,I,\tl,\ka,\tau)\cdot
\de_{Q^{-1}([(\ti\si,\ti\io,\ti\pi)])}$ as a sum over partial orders
$\pr$ on $I$ dominated by $\tl$. The operators
$\CF^\stk\bigl(Q(I,\pr,\tl)\bigr)$ commute with multiplication by
$\de_{Q^{-1}([\si,\io,\pi)])}$. Substituting this into \eq{as5eq20}
gives
\e
\begin{gathered}
\sum_{\substack{\text{p.o.s $\pr$ on $I$:}\\
\text{$\tl$ dominates $\pr$}}}
\begin{aligned}[t]
\CF^\stk\bigl(&Q(I,\pr,\tl)\bigr)\bigl(\dstb(X,I,\pr,\ka,\tau)
\cdot\de_{Q^{-1}([(\ti\si,\ti\io,\ti\pi)])}\bigr)=\\
&\CF^\stk\bigl(\bs\si(\{1,2\})_*\ci S(I,\tl,J)_*^{-1}\bigr)
\dst(Y,J,\ls,\ka,\tau).
\end{aligned}
\end{gathered}
\label{as5eq21}
\e

One can show using Theorem \ref{as3thm1} that the image under
$Q(I,\pr,\tl)_*$ of a best $(I,\pr)$-configuration in
$Q^{-1}([(\ti\si,\ti\io,\ti\pi)])$ is the image under
$\bigl(\bs\si(\{1,2\})_*\ci S(I,\tl,J)_*^{-1}\bigr)_*$ of a best
$(J,\ls)$-configuration if and only if $\pr\vert_{J}=\ls$. So
restricting \eq{as5eq21} to $\pr$ with $\pr\vert_J=\ls$ gives
\begin{equation*}
\smash{\sum_{\substack{\text{p.o.s $\pr$ on $I:\pr\vert_J=\ls$}\\
\text{and $j\npr a$ for all $j\in J$}}}}\!\!\!\!
\begin{aligned}[t]
\CF^\stk\bigl(&Q(I,\pr,\tl)\bigr)\bigl(\dstb(X,I,\pr,\ka,\tau)
\cdot\de_{Q^{-1}([(\ti\si,\ti\io,\ti\pi)])}\bigr)=\\
&\CF^\stk\bigl(\bs\si(\{1,2\})_*\ci S(I,\tl,J)_*^{-1}\bigr)
\dstb(Y,J,\ls,\ka,\tau).
\end{aligned}
\end{equation*}

Taking weighted Euler characteristics of both sides, using
\eq{as2eq3}, and summing over all $\ls,\ka$ with
$\tau\ci\ka\equiv\tau([X])$ proves:

\begin{prop} Let\/ $X,I$ be as above, $a\!\in\!I$ and\/
$J\!=\!I\sm\{a\}$. Define $\phi:I\!\ra\!\{1,2\}$ by
$\phi(a)\!=\!1$ and\/ $\phi(j)\!=\!2$ for $j\in J$. Let\/
$(\ti\si,\ti\io,\ti\pi)$ be a $(\{1,2\},\le)$-configuration
with\/ $\ti\si(\{1,2\})=X$ and\/ $\ti\si(\{1\})$ $\tau$-stable
with\/ $\tau\bigl([\ti\si(\{1\})]\bigr)\!=\!\tau\bigl([X]\bigr)$.
Define $Y=\ti\si(\{2\})$. Then
\e
\begin{gathered}
\sum_{\substack{\text{$\pr,\ka$: $(I,\pr,\ka)$ is $\A$-data,}\\
\text{$a$ is $\pr$-minimal, $[\ka(a)]\!=\![\ti\si(\{1\})]$,}\\
\text{$\ka(I)=[X]$, $\tau\ci\ka\equiv\tau([X])$}}}
\!\!\!\!\!\!\!\!\!
\begin{aligned}[t]
\chi^\na&\Bigl(\Mstb(X,I,\pr,\ka,\tau)_\A\,\cap\\
&Q(I,\pr,\{1,2\},\le,\phi)_*^{-1}\bigl(\{[(\ti\si,
\ti\io,\ti\pi)]\}\bigr)\Bigr)\!=
\end{aligned}\\
\sum_{\substack{\text{$\ls,\la$: $(J,\ls,\la)$ is $\A$-data,}\\
\text{$\la(J)=[Y]$, $\tau\ci\la\equiv\tau([X])$}}}
\chi^\na\bigl(\Mstb(Y,J,\ls,\la,\tau)_\A\bigr).
\end{gathered}
\label{as5eq22}
\e
Only finitely many terms in each sum are nonzero.
\label{as5prop1}
\end{prop}

We have not yet verified only finitely many terms in \eq{as5eq22}
are nonzero. Set $[X]=\al$, and suppose $\pr,\ka$ are as on the
l.h.s.\ of \eq{as5eq22} with $\Mstb(X,I,\pr,
\ka,\tau)_\A\ne\emptyset$. Let $K\ne\emptyset,I$ be an $(I,\pr)$
s-set, and set $\be=\ka(K)$, $\ga=\ka(I\sm K)$. Then $\be,\ga\in
C(\A)$ with $\al=\be+\ga$, and $\tau(\al)= \tau(\be)=\tau(\ga)$ as
$\tau\ci\ka\equiv\tau([X])$. If
$[(\si,\io,\pi)]\in\Mstb(X,I,\pr,\ka,\tau)_\A$ then $[\si(K)]
\in\Oss^\be(\tau)$ and $[\si(I\sm K)]\in\Oss^\ga(\tau)$,
so~$\Oss^\be(\tau)\ne\emptyset\ne\Oss^\ga(\tau)$.

Hence Proposition \ref{as4prop1} implies there are only finitely
many possibilities for $\ka(K),\ka(I\sm K)$. As this holds for all
$(I,\pr)$ s-sets $K\ne\emptyset,I$, and there are only finitely many
choices for $\pr$, there are also only finitely many choices for
$\ka$. So only finitely many terms on the l.h.s.\ of \eq{as5eq22}
are nonzero. The proof for the r.h.s.\ is the same. We can easily
extend this proof to fix not just one $\pr$-minimal element $a\in
I$, but a {\it minimal subset\/}~$A\subseteq I$.

\begin{prop} Let\/ $X,I$ be as above, $A\subseteq I$ and\/
$J=I\sm A$. Let\/ $b\notin A$ and set\/ $B=A\cup\{b\}$. Define
a partial order $\ps$ on $B$ by $r\ps s$ if either $r=s$ or $s=b$.
Define $\phi:I\ra B$ by $\phi(a)=a$ for $a\in A$ and\/ $\phi(i)=b$
for $i\in I\sm A$. Let\/ $(\ti\si,\ti\io,\ti\pi)$ be a
$(B,\ps)$-configuration with\/ $\ti\si(B)=X$ and\/ $\ti\si(\{a\})$
$\tau$-stable for all $a\in A$ with\/ $\tau\bigl([\ti\si(\{a\})]\bigr)=
\tau\bigl([X]\bigr)$. Define $Y=\ti\si(\{b\})$. Then
\e
\begin{gathered}
\sum_{\substack{\text{$\pr,\ka$: $(I,\pr,\ka)$ is $\A$-data,}\\
\text{each\/ $a\in A$ is $\pr$-minimal,}\\
\text{$[\ka(a)]\!=\![\ti\si(\{a\})]$, $\ka(I)\!=\![X]$,}\\
\text{$\tau\ci\ka\!\equiv\!\tau([X])$}}}
\begin{aligned}[t]
\chi^\na&\Bigl(\Mstb(X,I,\pr,\ka,\tau)_\A\,\cap\\
&Q(I,\pr,B,\ps,\phi)_*^{-1}\bigl(\{[(\ti\si,
\ti\io,\ti\pi)]\}\bigr)\Bigr)=
\end{aligned}\\
\sum_{\substack{\text{$\ls,\la$: $(J,\ls,\la)$ is $\A$-data,}\\
\text{$\la(J)=[Y]$, $\tau\ci\la\equiv\tau([X])$}}}
\chi^\na\bigl(\Mstb(Y,J,\ls,\la,\tau)_\A\bigr).
\end{gathered}
\label{as5eq23}
\e
Only finitely many terms in each sum are nonzero.
\label{as5prop2}
\end{prop}

We now calculate the Euler characteristic of the set of all
$[(\ti\si,\ti\io,\ti\pi)]$ satisfying the conditions in
Proposition~\ref{as5prop2}.

\begin{prop} Let\/ $X,I$ and\/ $S_m,k_m,l_m$ for $m=1,\ldots,n$
be as in Definition \ref{as5def}, and set\/ $k=\sum_{m=1}^nk_m$. For
$A\subseteq I$ with\/ $\md{A}\le k$, define $(B,\ps)$ as in
Proposition \ref{as5prop2}, and define
\e
\begin{gathered}
\M_A=\Bigl\{\bigl[(\ti\si,\ti\io,\ti\pi)\bigr]\in
\!\!\!\!\!\!\!\!\!\!\!
\coprod_{\substack{\text{$\mu$:$(B,\ps,\mu)$ is $\A$-data,}\\
\text{$\mu(B)=[X]$,}\\
\text{$\tau\ci\mu\equiv\tau([X])$}}}
\!\!\!\!\!\!\!\!\!\!\!\!\!\!\!\!
\M(X,B,\ps,\mu)_\A:
\begin{aligned}[t]
&\text{$\ti\si(\{a\})$ is $\tau$-stable}\\
&\text{for all $a\in A$}\Bigr\}.
\end{aligned}
\end{gathered}
\label{as5eq24}
\e
Then $\M_A$ is constructible with\/~$\chi^\na(\M_A)=k!/
\bigl(k-\md{A}\bigr)!$.
\label{as5prop3}
\end{prop}

\begin{proof} Write $P\bigl(\Hom(S_m,X)\bigr)$ for the
projective space of $\Hom(S_m,X)$. Then $P\bigl(\Hom
(S_m,X)\bigr)\cong\KP^{k_m-1}$ as $\Hom(S_m,X)\cong\K^{k_m}$.
Regard $P\bigl(\Hom(S_m,X)\bigr)$ as (the set of geometric
points of) a projective $\K$-variety. Define
\e
\begin{split}
{\mathcal N}_A=\Bigl\{&\psi:A\ra\ts\coprod_{m=1}^n
P\bigl(\Hom(S_m,X)\bigr):\text{$\psi$ is injective and}\\
&\quad\text{$\psi(A)\cap P\bigl(\Hom(S_m,X)\bigr)$ is
linearly independent for all $m$}\Bigr\}.
\end{split}
\label{as5eq25}
\e
Here a finite subset $S$ of a projective space $P(V)$ is
{\it linearly independent} if there exists no linear subspace
$U\subseteq V$ with $S\subseteq P(U)$ and $\dim U<\md{S}$.
Then ${\mathcal N}_A$ is an open set in the projective
$\K$-scheme $\prod_{a\in A}\coprod_{m=1}^nP\bigl(\Hom(S_m,X)\bigr)$,
so it is (the set of geometric points of) a quasiprojective
$\K$-scheme.

Define a map $\Phi:\M_A\ra{\mathcal N}_A$ as follows. If
$[(\ti\si,\ti\io,\ti\pi)]\in\M_A$ and $a\in A$ then
$\ti\si(\{a\})$ is $\tau$-stable with $\tau\bigl([\ti\si(
\{a\})]\bigr)=\tau\bigl([X]\bigr)$, so it follows that
$\ti\si(\{a\})$ is isomorphic to one of the $\tau$-{\it
stable factors} of $X$. Thus there exists an isomorphism
$i:S_m\ra\ti\si(\{a\})$ for some unique $m=1,\ldots,n$.
As $S_m$ is $\tau$-stable $\End(S_m)=\K$, so $i$ is
unique up to multiplication by a nonzero element of~$\K$.

As $i,\ti\io(\{a\},B)$ are injective we have $0\ne\ti\io(\{a\},B)
\ci i\in\Hom(S_m,X)$, and the class $\bigl[\ti\io(\{a\},B)\ci
i\bigr]\in P\bigl(\Hom(S_m,X)\bigr)$ is independent of choice
of $i$. Define $\psi(a)=\bigl[\ti\io(\{a\},B)\ci i\bigr]$.
This defines a map $\psi:A\ra\coprod_{m=1}^nP\bigl(\Hom(S_m,X)
\bigr)$. Define~$\Phi\bigl([(\ti\si,\ti\io,\ti\pi)]\bigr)=\psi$.

Now $\Phi$ essentially maps $[(\ti\si,\ti\io,\ti\pi)]$ in $\M_A$ to
a set of stable subobjects in $X$ parametrized by $A$. Using
\cite[Th.s 4.2 \& 4.5]{Joyc3} we deduce necessary and sufficient
conditions for such a set of subobjects to come from a
$(B,\ps)$-configuration $(\ti\si,\ti\io,\ti\pi)$ with $\ti\si(B)=X$,
and they turn out to be that $\psi$ is injective and $\psi(A)\cap
P\bigl(\Hom(S_m,X)\bigr)$ is linearly independent for all $m$. It
follows that $\Phi$ maps to ${\mathcal N}_A$, and is a 1-1
correspondence.

Using \cite[Ass.~7.1(iv)]{Joyc3} and general facts from \cite{Joyc1}
and \cite{Joyc3}, it is not difficult to see that $\Phi$ is a {\it
pseudoisomorphism}, in the sense of \cite[\S 4.2]{Joyc1}. The point
of invoking \cite[Ass.~7.1(iv)]{Joyc3} is that it gives us a {\it
tautological morphism} $\th_{S_m,X}$, a family of morphisms $S_m\ra
X$ parametrized by the base $\K$-scheme $\Hom(S_m,X)$. Using this it
is easy, for instance, to construct a $\K$-substack $P_m$ of
$\fM(X,\{1,2\},\le,\mu)_\A$ isomorphic to $P\bigl(\Hom(S_m,X)
\bigr)$, where $\mu(1)=[S_m]$ and $\mu(2)=[X]-[S_m]$, with
\begin{equation*}
P_m(\K)=\bigl\{[(\ti\si,\ti\io,\ti\pi)]\in\M(X,\{1,2\},
\le,\mu)_\A:\ti\si(\{1\})\cong S_m\bigr\}.
\end{equation*}
When $\md{A}=1$ we have $\M_A\cong\coprod_{m=1}^nP_m(\K)$ and
${\mathcal N}_A\cong\coprod_{m=1}^nP\bigl(\Hom(S_m,X)\bigr)$,
and the result follows. The case $\md{A}>1$ is a
straightforward generalization.

From \cite[Def.~4.8]{Joyc1} we see that $\chi^\na(\M_A)=
\chi({\mathcal N}_A)$. Thus the proposition follows from
$\chi({\mathcal N}_A)=k!/\bigl(k-\md{A}\bigr)!$. One can prove this
using $P\bigl(\Hom(S_m,X)\bigr)\cong\KP^{k_m-1}$,
$k=\sum_{m=1}^nk_m$, \eq{as5eq25}, and properties of $\chi$
including \eq{as2eq2}, by a long but elementary calculation that we
leave as an exercise.
\end{proof}

In the next theorem, note that the set of $\pr$-minimal elements in
$I$ {\it contains} $A$ in \eq{as5eq26}, and is {\it equal to} $A$
in~\eq{as5eq27}.

\begin{thm} Let\/ $X,I$ and\/ $S_m,k_m,l_m$ for $m=1,\ldots,n$
be as in Definition \ref{as5def}, and set\/ $k=\sum_{m=1}^nk_m$.
Then for each\/ $A\subseteq I$ with\/ $\md{A}\le k$ we have
\ea
\sum_{\substack{\text{$\pr,\ka$: $(I,\pr,\ka)$ is $\A$-data,}\\
\text{each $a\in A$ is $\pr$-minimal,}\\
\text{$\ka(I)=[X]$, $\tau\ci\ka\equiv\tau([X])$}}}\!\!\!\!\!\!\!\!\!\!\!
\chi^\na\bigl(\Mstb(X,I,\pr,\ka,\tau)_\A\bigr)&=
\frac{\bigl(\md{I}-\md{A}\bigr)!k!}{\bigl(k-\md{A}\bigr)!},
\;\>\text{and}
\label{as5eq26}\\
\sum_{\substack{\text{$\pr,\ka$: $(I,\pr,\ka)$ is $\A$-data,}\\
\text{$A$ is the $\pr$-minimal set,}\\
\text{$\ka(I)=[X]$, $\tau\ci\ka\equiv\tau([X])$}}}\!\!\!\!\!\!\!\!\!\!\!
\chi^\na\bigl(\Mstb(X,I,\pr,\ka,\tau)_\A\bigr)&=
\begin{cases} 0, & \md{A}<k,\\
\bigl(\md{I}-k\bigr)!k!, & \md{A}=k.\end{cases}
\label{as5eq27}
\ea
Only finitely many terms in each sum are nonzero.
\label{as5thm5}
\end{thm}

\begin{proof} The argument after Proposition \ref{as5prop1}
shows only finitely many terms in \eq{as5eq26} and \eq{as5eq27} are
nonzero. Note that the $\pr$-minimal set $A$ in $I$ always has
$1\le\md{A}\le k$ by definition of $k,k_m$, as $X$ has only $k$
linearly independent $\tau$-stable subobjects. First we show
\eq{as5eq26} and \eq{as5eq27} are equivalent.

Suppose \eq{as5eq27} holds. Then letting the $\pr$-minimal set in
\eq{as5eq26} be $A'$, for $A\!\subseteq\!A$ with $\md{A}\!\le\!k$
summing \eq{as5eq27} with $A'$ in place of $A$ over all
$A'\!\subseteq\!A$ and using a simple combinatorial argument proves
\eq{as5eq26}. Now \eq{as5eq27} holds trivially when $\md{A}=0$ as
both sides are zero. Hence, \eq{as5eq27} for
$1\!\le\!\md{A}\!\le\!k$ implies \eq{as5eq26}. By a more complicated
argument we find \eq{as5eq26} for $1\!\le\!\md{A}\!\le\!k$ implies
\eq{as5eq27}. Hence, if \eq{as5eq26} holds when
$1\!\le\!\md{A}\!\le\!k$, then both \eq{as5eq26} and \eq{as5eq27}
hold for~$\md{A}\!\le\!k$.

We can now prove the theorem by induction on $\md{I}$. The result is
trivial when $\md{I}=1$, giving the first step. Suppose by induction
that \eq{as5eq26} and \eq{as5eq27} hold whenever $\md{I}\le m$, and
let $\md{I}=m+1$. Let $A\subseteq I$ with $1\le\md{A}\le k$, and set
$J=I\sm A$. Let $\M_A$ be as in Proposition \ref{as5prop3}. Define
$\fG$ and $T\subseteq\fG(\K)$ by
\e
\fG=\!\!\!\!\!\!\!\!\!\!\!\!\!\!
\coprod_{\substack{\text{$\pr,\ka$: $(I,\pr,\ka)$ is $\A$-data,}\\
\text{each $a\in A$ is $\pr$-minimal,}\\
\text{$\ka(I)=[X]$, $\tau\ci\ka\equiv\tau([X])$}}}
\!\!\!\!\!\!\!\!\!\!\!\!\!\!
\fM(X,I,\pr,\ka)_\A,
\quad
T=\!\!\!\!\!\!\!\!\!\!\!\!\!\!
\coprod_{\substack{\text{$\pr,\ka$: $(I,\pr,\ka)$ is $\A$-data,}\\
\text{each $a\in A$ is $\pr$-minimal,}\\
\text{$\ka(I)=[X]$, $\tau\ci\ka\equiv\tau([X])$}}}
\!\!\!\!\!\!\!\!\!\!\!\!\!\!
\Mstb(X,I,\pr,\ka,\tau)_\A.
\label{as5eq28}
\e

Let $B,\ps,\phi$ be as in Proposition \ref{as5prop2}. For each
$\pr,\ka$ in the definition of $\fG$ in \eq{as5eq28}, define
$\mu:B\ra C(\A)$ by $\mu(c)=\ka(\phi^{-1}(\{c\}))$. Then we have a
1-morphism $Q(I,\pr,B,\ps,\phi):\fM(X,I,
\pr,\ka)_\A\ra\fM(X,B,\ps,\mu)_\A$. Define a 1-morphism
\e
\psi:\fG\ra
\!\!\!\!\!\!\!\!\!\!\!\!\!\!\!\!
\coprod_{\substack{\text{$\mu$:$(B,\ps,\mu)$ is $\A$-data,}\\
\text{$\mu(B)=[X]$,}\\
\text{$\tau\ci\mu\equiv\tau([X])$}}}
\!\!\!\!\!\!\!\!\!\!\!\!\!\!\!\!
\fM(X,B,\ps,\mu)_\A\quad\text{by}\quad
\psi=\!\!\!\!\!\!\!\!\!\!\!\!\!
\coprod_{\substack{\text{$\pr,\ka$: $(I,\pr,\ka)$ is $\A$-data,}\\
\text{each $a\in A$ is $\pr$-minimal,}\\
\text{$\ka(I)=[X]$, $\tau\ci\ka\equiv\tau([X])$}}}
\!\!\!\!\!\!\!\!\!\!\!\!\!
Q(I,\pr,B,\ps,\phi).
\label{as5eq29}
\e
Comparing \eq{as5eq24}, \eq{as5eq29} and the definition of $\phi$
shows that $\psi_*$ maps $T\ra\M_A$. Let $[(\ti\si,
\ti\io,\ti\pi)]\in\M_A$, and define $Y=\ti\si(\{b\})$. Then
\eq{as5eq23} gives an expression for $\chi^\na\bigl(
T\cap\psi_*^{-1}([(\ti\si,\ti\io,\ti\pi)])\bigr)$. Now the r.h.s.\
of \eq{as5eq23} is the l.h.s.\ of \eq{as5eq26} with $Y$ in place of
$X$, $J$ in place of $I$, and $\emptyset$ in place of~$A$.

Since $\md{I}=m+1$, $\md{A}\ge 1$ and $J=I\sm A$ we have $\md{J}\le
m$. Hence by the inductive hypothesis, \eq{as5eq26} holds for
$Y,J,\emptyset$. So for all $[(\ti\si,\ti\io,\ti\pi)]\in\M_A$ we
have
\begin{equation*}
\chi^\na\bigl(T\cap\psi_*^{-1}([(\ti\si,\ti\io,\ti\pi)])\bigr)
=\md{J}!=\bigl(\md{I}-\md{A}\bigr)!.
\end{equation*}
Proposition \ref{as5prop3} and general properties of $\chi^\na$ now
imply that
\begin{equation*}
\sum_{\substack{\text{$\pr,\ka$: $(I,\pr,\ka)$ is $\A$-data,}\\
\text{each $a\in A$ is $\pr$-minimal,}\\
\text{$\ka(I)=[X]$, $\tau\ci\ka\equiv\tau([X])$}}}
\!\!\!\!\!\!\!\!\!\!\!\!\!\!\!\!
\begin{aligned}[t]
\chi^\na\bigl(\Mstb(X,I,\pr,\ka,\tau)_\A\bigr)&=\!\chi^\na(T)
\!=\!\bigl(\md{I}\!-\!\md{A}\bigr)!\chi^\na(\M_A)\\
&=\bigl(\md{I}-\md{A}\bigr)!k!/\bigl(k-\md{A}\bigr)!.
\end{aligned}
\end{equation*}
Hence \eq{as5eq26} holds for $1\le\md{A}\le k$ with this fixed $I$,
and so \eq{as5eq26} and \eq{as5eq27} hold for $\md{A}\le k$ with
this $I$ from above. This completes the inductive step.
\end{proof}

Equation \eq{as5eq26} with $A=\emptyset$ calculates the Euler
characteristic of the family of all best $\tau$-stable
configurations for~$X$.

\begin{cor} Let Assumption \ref{as3ass} hold, $(\tau,T,\le)$
be a permissible stability condition on $\A$, and\/ $X\in\A$
be $\tau$-semistable. Fix a finite set\/ $I$ such that\/ $X$
has $\md{I}$ $\tau$-stable factors in Theorem \ref{as4thm2},
counted with multiplicity. Then
\e
\sum_{\substack{\text{$\pr,\ka$: $(I,\pr,\ka)$ is $\A$-data,}\\
\text{$\ka(I)=[X]$, $\tau\ci\ka\equiv\tau([X])$}}}
\chi^\na\bigl(\Mstb(X,I,\pr,\ka,\tau)_\A\bigr)=\md{I}!.
\label{as5eq30}
\e
Only finitely many $\Mstb(X,I,\pr,\ka,\tau)_\A$ in this sum
are nonempty.
\label{as5cor1}
\end{cor}

We turn this into an identity on constructible functions:

\begin{thm} Let Assumption \ref{as3ass} hold, $(\tau,T,\le)$ be a
permissible stability condition on $\A$, and\/ $\al\in C(\A)$. Then
\e
\sum_{\substack{\text{iso. classes}\\ \text{of finite }\\
\text{sets $I$}}} \frac{1}{\md{I}!}\cdot
\sum_{\substack{
\text{$\pr,\ka$: $(I,\pr,\ka)$ is $\A$-data,}\\
\text{$\ka(I)=\al$, $\tau\ci\ka\equiv\tau(\al)$}}}\!\!\!\!\!\!\!
\CF^\stk\bigl(\bs\si(I)\bigr)\dstb(I,\pr,\ka,\tau)=\dss^\al(\tau).
\label{as5eq31}
\e
Only finitely many functions $\dstb(I,\pr,\ka,\tau)$ in this sum are
nonzero.
\label{as5thm6}
\end{thm}

\begin{proof} A similar proof to that in \S\ref{as52} showing
\eq{as5eq11} has only finitely many nonzero terms proves that only
finitely many $\dstb(I,\pr, \ka,\tau)$ are nonzero in \eq{as5eq31}.
Let $(I,\pr,\ka)$ be as in \eq{as5eq31} and $[(\si,\io,\pi)]\in
\Mstb(I,\pr,\ka,\tau)_\A$. Then $\tau\ci\ka\equiv\tau(\al)$ implies
$\si(I)$ is $\tau$-semistable, so $[\si(I)]\in\Oss^\al(\tau)$. Hence
both sides of \eq{as5eq31} are zero outside $\Oss^\al(\tau)$. But if
$X\in\A$ is $\tau$-semistable with $[X]=\al$ in $K(\A)$ then
\eq{as5eq30} and the definitions of $\dstb(I,\pr,\ka,\tau)$ and
$\CF^\stk$ imply that both sides of \eq{as5eq31} are equal at
$[X]\in\Oss^\al(\tau)$, by an argument similar to
Theorem~\ref{as5thm2}.
\end{proof}

\subsection{Counting best $\tau$-stable refinements}
\label{as54}

Our next result in effect computes the Euler characteristic
of the family of all {\it best\/ $\tau$-stable refinements} of
a $\tau$-semistable $(K,\tl)$-configuration~$(\si,\io,\pi)$.

\begin{thm} Let Assumption \ref{as3ass} hold, and\/
$(\tau,T,\le)$ be a permissible stability condition on $\A$. Then
for all\/ $\A$-data $(K,\tl,\mu)$ we have
\e
\begin{aligned}
\sum_{\substack{\text{iso. classes}\\ \text{of finite }\\
\text{sets $I$}}} \frac{1}{\md{I}!}\cdot
\sum_{\substack{
\text{$\pr,\ka,\phi$: $(I,\pr,\ka)$ is $\A$-data,}\\
\text{$\phi:I\ra K$ is surjective,}\\
\text{$i\pr j$ implies $\phi(i)\tl\phi(j)$,}\\
\text{$\ka(\phi^{-1}(k))=\mu(k)$ for $k\in K$,}\\
\text{$\tau\ci\mu\ci\phi\equiv\tau\ci\ka:I\ra T$}}}
\begin{aligned}[t]
\CF^\stk\bigl(Q(I,\pr,K,\tl,\phi)\bigr)&\\
\dstb(I,\pr,\ka,\tau)&=\\
\dss(K,\tl,\mu&,\tau).
\end{aligned}
\end{aligned}
\label{as5eq32}
\e
Only finitely many functions $\dstb(I,\pr,\ka,\tau)$ in this sum are
nonzero.
\label{as5thm7}
\end{thm}

\begin{proof} A similar proof to that in \S\ref{as52} showing
\eq{as5eq14} has only finitely many nonzero terms proves that only
finitely many $\dstb(I,\pr,\ka,\tau)$ are nonzero in \eq{as5eq32},
as \eq{as5eq31} has only finitely many nonzero terms. For each $k\in
K$, let $(I_k,\ls_k,\ka_k)$ be $\A$-data with $\ka_k(I_k)=\mu(k)$
and $\tau\ci\ka_k\equiv\tau\ci\mu(k)$. Define $I=\coprod_{k\in
K}I_k$ and $\phi:I\ra K$ by $\phi(i)=k$ if $i\in I_k$. Define $\ls$
on $I$ by $i\ls j$ for $i,j\in I$ if either (a) $\phi(i)\tl\phi(j)$
and $\phi(i) \ne\phi(j)$, or (b) $\phi(i)=\phi(j)=k$ and $i\ls_k j$.
Then $i\ls j$ implies $\phi(i)\tl\phi(j)$. Define $\ka:I\ra K(\A)$
by $\ka\vert_{I_k}=\ka_k$. Then $(I,\ls,\ka)$ is $\A$-data with
$\ka(\phi^{-1}(k))=\mu(k)$ for $k\in K$ and $\tau\ci\mu
\ci\phi\equiv\tau\ci\ka$, as in~\eq{as5eq32}.

As for \eq{as5eq15}, the following commutative diagram is a
Cartesian square
\e
\begin{gathered}
\xymatrix@C=95pt@R=9pt{
\fM(I,\ls,\ka)_\A \ar[d]^{\,\prod_{k\in K}S(I,\ls,I_k)}
\ar[r]_{Q(I,\ls,K,\tl,\phi)} &
\fM(K,\tl,\mu)_\A \ar[d]_{\prod_{k\in K}\bs\si(\{k\})\,} \\
\prod_{k\in K}\fM(I_k,\ls_k,\ka_k)_\A \ar[r]^{\prod_{k\in
K}\bs\si(I_k)} & \prod_{k\in K}\fObj_\A^{\mu(k)},
}
\end{gathered}
\label{as5eq33}
\e
with representable rows and finite type columns. Since
$\dstb(I_k,\ls_k,\ka_k,\tau)\in\CF(\fM(I_k,\ls_k,\ka_k)_\A)$ by
Theorem \ref{as4thm3}, we may apply Theorem \ref{as2thm1} to
\eq{as5eq33} and the function $\prod_{k\in K}
\dstb(I_k,\ls_k,\ka_k,\tau)\in\CF\bigl(\prod_{k\in K}
\fM(I_k,\ls_k,\ka_k)_\A\bigr)$. This yields
\e
\begin{gathered}
\CF^\stk\bigl(Q(I,\ls,K,\tl,\phi)\bigr)\bigl[\ts\prod_{k\in K}
S(I,\ls,I_k)^*\dstb(I_k,\ls_k,\ka_k,\tau)\bigr]=\\
\ts\prod_{k\in K}\bs\si(\{k\})^*\bigl[\CF^\stk(\bs\si(I_k))
\dstb(I_k,\ls_k,\ka_k,\tau)\bigr].
\end{gathered}
\label{as5eq34}
\e
Using Theorem \ref{as5thm2}, the definition of $\dst(I,\ls,
\ka,\tau)$ and $I=\coprod_{k\in K}I_k$ gives
\e
\begin{gathered}
\sum_{\substack{\text{p.o.s $\pr$ on $I$:}\\
\text{$\ls$ dominates $\pr$}}}
\begin{aligned}[t]\CF^\stk\bigl(&Q(I,\pr,\ls)\bigr)\dstb(I,\pr,\ka,\tau)
=\dst(I,\ls,\ka,\tau)\\
&=\ts\prod_{k\in K}S(I,\ls,I_k)^*
\dst(I_k,\ls_k,\ka_k,\tau).
\end{aligned}
\end{gathered}
\label{as5eq35}
\e

One can show using Theorem \ref{as3thm1} that the image of a best
configuration under $S(I,\ls,I_k)_*\ci Q(I,\ls,\pr)_*$ is best if
and only if $\pr\vert_{I_k}=\ls_k$. So restricting \eq{as5eq35} to
$\pr$ with $\pr\vert_{I_k}=\ls_k$ for all $k$ proves that
\begin{equation*}
\sum_{\substack{\text{p.o.s $\pr$ on $I$:}\\
\text{$\pr\vert_{I_k}=\ls_k$, $k\in K$,}\\
\text{$i\pr j$ implies $\phi(i)\!\tl\!\phi(j)$}}}
\begin{aligned}[t]
\CF^\stk\bigl(&Q(I,\pr,\ls)\bigr)\dstb(I,\pr,\ka,\tau)\\
&=\ts\prod_{k\in K}S(I,\ls,I_k)^*\dstb(I_k,\ls_k,\ka_k,\tau).
\end{aligned}
\end{equation*}
Applying $\CF^\stk\bigl(Q(I,\ls,K,\tl,\phi)\bigr)$ to this equation,
noting that $Q(I,\ls,K,\tl,\phi)\ci Q(I,\pr,\ls)=
Q(I,\pr,K,\tl,\phi)$ and using \eq{as2eq3} and \eq{as5eq34} gives
\e
\begin{gathered}
\sum_{\substack{\text{p.o.s $\pr$ on $I$:}\\
\text{$\pr\vert_{I_k}=\ls_k$, $k\in K$,}\\
\text{$i\pr j$ implies $\phi(i)\!\tl\!\phi(j)$}}}\!\!\!\!\!\!\!\!
\begin{aligned}[t]
\CF^\stk&\bigl(Q(I,\pr,K,\tl,\phi)\bigr)\dstb(I,\pr,\ka,\tau)\\
=&\,\ts\prod_{k\in K}\bs\si(\{k\})^*\bigl[\CF^\stk(\bs\si(I_k))
\dstb(I_k,\ls_k,\ka_k,\tau)\bigr].
\end{aligned}
\end{gathered}
\label{as5eq36}
\e

Now using Theorem \ref{as5thm6} to rewrite $\dss^{\mu(k)}(\tau)$ for
each $k\in K$ yields
\e
\begin{aligned}
\dss&(K,\tl,\mu,\tau)=\ts\prod_{k\in K}\bs\si(\{k\})^*
\bigl(\dss^{\mu(k)}(\tau)\bigr)\\
&=\prod_{k\in K}\,
\sum_{\substack{\text{iso. classes}\\ \text{of finite }\\
\text{sets $I_k$}}} \frac{1}{\md{I_k}!}\cdot
\!\!\!\!
\sum_{\substack{
\text{$\ls_k,\ka_k$: $(I_k,\ls_k,\ka_k)$ is $\A$-data,}\\
\text{$\ka_k(I_k)=\mu(k)$, $\tau\ci\ka_k\equiv\tau(\mu(k))$}}}
\!\!\!\!\!\!\!\!\!\!\!\!\!\!\!\!\!\!
\begin{aligned}[t]
\bs\si(\{k\})^*\bigl[&\CF^\stk\bigl(\bs\si(I_k)\bigr)\\
&\dstb(I_k,\ls_k,\ka_k,\tau)\bigr]
\end{aligned}
\\
&=
\sum_{\substack{\text{iso. classes}\\ \text{of finite sets}\\
\text{$I_k$, all $k\in K$}}}
\Bigl[\,\prod_{k\in K}\frac{1}{\md{I_k}!}\,\Bigr]\cdot
\sum_{\substack{
\text{$\ls_k,\ka_k$ for all $k\in K$:}\\
\text{$(I_k,\ls_k,\ka_k)$ is $\A$-data,}\\
\text{$\ka_k(I_k)=\mu(k)$, $\tau\ci\ka_k\equiv\tau(\mu(k))$}}}\\
&\quad
\sum_{\substack{\text{p.o.s $\pr$ on $I=\coprod_{k\in K}I_k$:}\\
\text{$\pr\vert_{I_k}=\ls_k$, $k\in K$,}\\
\text{$i\pr j$ implies $\phi(i)\!\tl\!\phi(j)$}}}\!\!\!\!\!\!\!\!
\CF^\stk\bigl(Q(I,\pr,K,\tl,\phi)\bigr)\dstb(I,\pr,\ka,\tau),
\end{aligned}
\label{as5eq37}
\e
substituting in \eq{as5eq36} at the last step.

The sums over $I_k,\ls_k,\ka_k$ for all $k\in K$ and $\pr$ in
\eq{as5eq37} are equivalent to the sums over $I,\pr,\ka$ and $\phi$
in \eq{as5eq32}, with the following proviso. If we choose sets $I_k$
for $k\in K$ in \eq{as5eq37}, then in \eq{as5eq32} the first sum
fixes a unique set $I$ with $\md{I}=\sum_{k\in K}\md{I_k}$, and
there are then $\md{I}!/ \prod_{k\in K}\md{I_k}!$ possible
surjective maps $\phi:I\ra K$ with $\bmd{\phi^{-1}(\{k\})}=\md{I_k}$
for all $k\in K$. Thus, for each choice of data $I_k$ in
\eq{as5eq37}, there are $\md{I}!/\prod_{k\in K}\md{I_k}!$
corresponding choices of data $I,\phi$ in \eq{as5eq32}. This exactly
cancels the difference between the factors $\prod_{k\in
K}1/\md{I_k}!$ in \eq{as5eq37} and $1/\md{I}!$ in \eq{as5eq32}. So
\eq{as5eq37} and \eq{as5eq32} are equivalent, completing the proof.
\end{proof}

\section{Combinatorial inversion of the identities of \S\ref{as5}}
\label{as6}

Next we prove some more identities involving pushforwards of the
characteristic functions $\dss,\dsi,\dst,\dssb,\dsib,\dstb
(I,\pr,\ka,\tau)$ under 1-morphisms $Q(I,\pr,K,\tl,\phi)$. Equations
\eq{as5eq7}, \eq{as5eq8}, \eq{as5eq9}, \eq{as5eq14} and \eq{as5eq32}
above are of this type. By inverting these explicitly we find six
further identities, \eq{as6eq3}, \eq{as6eq4}, \eq{as6eq5},
\eq{as6eq10}, \eq{as6eq11} and \eq{as6eq14} below. These mean that
given the $\fM(I,\pr,\ka)_\A$ and $Q(I,\pr,K,\tl,\phi)$, any one of
the six families $\dss,\dsi,\dst,\dssb,\dsib,\dstb(*,\tau)$
determines the other five.

In contrast to \S\ref{as5}, the arguments of this section are all
{\it combinatorial\/} in nature. Our principal techniques are
substituting one complicated sum inside another, and rearranging the
order of summation. We continue to suppose $\K$ has characteristic
zero.

\subsection{Inverting identities \eq{as5eq7}--\eq{as5eq9}}
\label{as61}

In \eq{as5eq7}--\eq{as5eq9} we wrote
$\dss,\dsi,\dst(I,\tl,\ka,\tau)$ in terms of
$\dssb,\dsib,\dstb(I,\pr,\ka,\tau)$. We now invert these. We shall
need some integers~$n(I,\pr,\tl)$.

\begin{dfn} Let $I$ be a finite set, and $\ls,\tl$ partial
orders on $I$, where $\tl$ dominates $\ls$. Define an integer
\e
n(I,\ls,\tl)=\sum_{\substack{\text{$n\ge 0$,
$\ls=\pr_0,\pr_1,\ldots,\pr_n=\tl$:}\\
\text{$\pr_m$ is a partial order on $I$, $0\le m\le n$,}\\
\text{$\pr_m$ strictly dominates $\pr_{m-1}$, $1\le m\le n$}}}(-1)^n.
\label{as6eq1}
\e
\label{as6def1}
\end{dfn}

If $\tl$ dominates $\ls$ by $l$ steps, as in Definition
\ref{as3def1}, then $0\le n\le l$ in \eq{as6eq1}, so the sum
\eq{as6eq1} is finite. The $n(I,\ls,\tl)$ satisfy the following
equation:

\begin{prop} Let\/ $I$ be a finite set and\/ $\ls,\tl$ partial
orders on $I$, where $\tl$ dominates $\ls$. Then
\e
\sum_{\substack{\text{partial orders $\pr$ on $I$:}\\
\text{$\tl$ dominates $\pr$ dominates $\ls$}}}n(I,\pr,\tl)=
\begin{cases} 1, & \ls=\tl, \\ 0, & \ls\ne\tl.\end{cases}
\label{as6eq2}
\e
Also, the same equation holds with\/ $n(I,\pr,\tl)$
replaced by~$n(I,\ls,\pr)$.
\label{as6prop1}
\end{prop}

\begin{proof} If $\ls=\tl$ then in \eq{as6eq1} there is only
one possibility, $n=0$ and $\ls=\pr_0=\tl$, so $n(I,\ls,\tl)=1$.
Also in \eq{as6eq2} we have $\pr=\ls=\tl$, so the top line of
\eq{as6eq2} is immediate. Suppose $\ls\ne\tl$. Then every term in
\eq{as6eq1} has $n\ge 1$, and by setting $\pr=\pr_1$, replacing $n$
by $n-1$ and $\pr_m$ by $\pr_{m+1}$ we rewrite \eq{as6eq1} as
\begin{align*}
n(I,\ls,\tl)&=\sum_{\substack{\text{p.o.s $\pr$ on $I$:}\\
\text{$\tl$ dominates $\pr$}\\
\text{$\pr$ strictly dominates $\ls$}}}\,\,\,\,
\sum_{\substack{\text{$n\ge 0$, $\pr=\pr_0,\ldots,\pr_n=\tl$:}\\
\text{$\pr_m$ is a p.o. on $I$,}\\
\text{$\pr_m$ strictly dominates $\pr_{m-1}$}}}(-1)^{n+1}\\
&=-\sum_{\substack{\text{p.o.s $\pr$ on $I$:}\\
\text{$\tl$ dominates $\pr$}\\
\text{$\pr$ strictly dominates $\ls$}}}n(I,\pr,\tl),
\end{align*}
The bottom line of \eq{as6eq2} follows immediately. We prove
\eq{as6eq2} with $n(I,\pr,\tl)$ replaced by $n(I,\ls,\pr)$ in a
similar way, writing $\pr$ for $\pr_{n-1}$ in~\eq{as6eq1}.
\end{proof}

Here are the inverses of the identities of Theorem~\ref{as5thm2}.

\begin{thm} Let Assumption \ref{as3ass} hold, $(\tau,T,\le)$
be a permissible weak stability condition on $\A$, and\/
$(I,\tl,\ka)$ be $\A$-data, as in Definition \ref{as3def5}. Then
\ea
\sum_{\substack{\text{p.o.s $\pr$ on $I$:}\\
\text{$\tl$ dominates $\pr$}}}\!\!\!\!\!\!\!\!\!
n(I,\pr,\tl)\CF^\stk\bigl(Q(I,\pr,\tl)\bigr)
\dss(I,\pr,\ka,\tau)&=\dssb(I,\tl,\ka,\tau),
\label{as6eq3}\\[3pt]
\sum_{\substack{\text{p.o.s $\pr$ on $I$:}\\
\text{$\tl$ dominates $\pr$}}}\!\!\!\!\!\!\!\!\!
n(I,\pr,\tl)\CF^\stk\bigl(Q(I,\pr,\tl)\bigr)
\dsi(I,\pr,\ka,\tau)&=\dsib(I,\tl,\ka,\tau),
\label{as6eq4}\\[3pt]
\sum_{\substack{\text{p.o.s $\pr$ on $I$:}\\
\text{$\tl$ dominates $\pr$}}}\!\!\!\!\!\!\!\!\!
n(I,\pr,\tl)\CF^\stk\bigl(Q(I,\pr,\tl)\bigr)
\dst(I,\pr,\ka,\tau)&=\dstb(I,\tl,\ka,\tau).
\label{as6eq5}
\ea
\label{as6thm1}
\end{thm}

\begin{proof} Substituting \eq{as5eq7} into the left hand side
of \eq{as6eq3} gives
\begin{align*}
\sum_{\substack{\text{p.o.s $\pr$ on $I$:}\\
\text{$\tl$ dominates $\pr$}}}
n(I,\pr,\tl)\CF^\stk\bigl(Q(I,\pr,\tl)\bigr)\dss(I,\pr,\ka,\tau)=&\\
\!\!\!\sum_{\substack{\text{p.o.s $\pr$ on $I$:}\\
\text{$\tl$ dominates $\pr$}}}\!\!\!\!\!\!\! n(I,\pr,\tl)
\sum_{\substack{\text{p.o.s $\ls$ on $I$:}\\
\text{$\pr$ dominates $\ls$}}}
\begin{aligned}[t]
\CF^\stk\bigl(Q(I,\pr,\tl)\bigr)\Bigl[
\CF^\stk\bigl(Q(I,\ls,\pr)\bigr)\quad\,\,{}&\\
\dssb(I,\ls,\ka,\tau)\Bigr]=&
\end{aligned}\\
\sum_{\substack{\text{p.o.s $\ls$ on $I$:}\\
\text{$\tl$ dominates $\ls$}}}
\raisebox{-12pt}{\begin{Large}$\displaystyle\Biggl[$\end{Large}}
\sum_{\substack{\text{p.o.s $\pr$ on $I$:}\\
\text{$\tl$ dominates $\pr$,}\\
\text{$\pr$ dominates $\ls$}}}n(I,\pr,\tl)
\raisebox{-12pt}{\begin{Large}$\displaystyle\Biggr]$\end{Large}}
\begin{aligned}[t]
\CF^\stk\bigl(Q(I,\ls,\tl)\bigr)\,{}&\\
\dssb(I,\ls,\ka,\tau),&
\end{aligned}
\end{align*}
exchanging sums over $\pr,\ls$ and using $Q(I,\pr,\tl)\ci
Q(I,\ls,\pr)=Q(I,\ls,\tl)$ and \eq{as2eq3}. By \eq{as6eq2} the
bracketed sum on the last line is 0 unless $\ls=\pr=\tl$, when it is
1. But then $Q(I,\ls,\tl)$ is the identity, so the final line
reduces to $\dssb(I,\tl,\ka,\tau)$, giving \eq{as6eq3}. The proofs
of \eq{as6eq4}--\eq{as6eq5} from \eq{as5eq8}--\eq{as5eq9} are the
same.
\end{proof}

\subsection{Inverting \eq{as5eq11} and \eq{as5eq14}}
\label{as62}

We invert \eq{as5eq11} to write $\dsi^\be(\tau)$ in terms
of~$\dss^{\ka(i)}(\tau)$.

\begin{thm} Let Assumption \ref{as3ass} hold, $(\tau,T,\le)$ be a
permissible weak stability condition on $\A$, and\/ $\be\in C(\A)$.
Then
\e
\sum_{n=1}^\iy\frac{(-1)^{n-1}}{n}\cdot \!\!\!\!
\sum_{\substack{\text{$\ka:\{1,\ldots,n\}\ra C(\A):$}\\
\text{$\ka(\{1,\ldots,n\})\!=\!\be$, $\tau\!\ci\!\ka\!\equiv
\!\tau(\be)$}}}\!\!\!\!\!\!\!\!\!\!\!\!\!\!\!\!\!\!\!
\begin{aligned}[t]
\CF^\stk\bigl(\bs\si(\{1,\ldots,n\})\bigr)\dss(\{1,\ldots,n\},
\bu,\ka,\tau)&\\
=\dsi^\be(\tau)&,
\end{aligned}
\label{as6eq6}
\e
where $\bu$ is the partial order on $\{1,\ldots,n\}$ with\/ $i\bu j$
if and only if\/ $i=j$. Only finitely many functions
$\dss(\{1,\ldots,n\},\bu,\ka,\tau)$ in this sum are nonzero.
\label{as6thm2}
\end{thm}

\begin{proof} We could give a straight combinatorial proof of
\eq{as6eq6}, but the author finds the following infinite series
proof more attractive, and we will also reuse the method in Theorem
\ref{as7thm1}. The motivation is that equation \eq{as5eq11} looks
like an exponential series, so its inverse \eq{as6eq6} should look
like a log. To make \eq{as5eq11} look more like an exponential,
define $\square:\CF(\fObj_\A)\t\CF(\fObj_\A)\ra \CF(\fObj_\A)$ by
$f\square g=P_{\sst(\{1,2\},\bu)}(f,g)$. Then as in \cite[\S
4.8]{Joyc4}, $\square$ is an {\it associative, commutative
multiplication} on $\CF(\fObj_\A)$. Also in \eq{as5eq11} we have
\begin{align*}
\CF^\stk\bigl(\bs\si(\{1,\ldots,n\})\bigr)\dsi(\{1,\ldots,n\},\bu,
\ka,\tau)&=\\
P_{\sst(\{1,\ldots,n\},\bu)}\bigl(\dsi^{\ka(i)}(\tau):i=1,\ldots,n\bigr)
&=\square_{i=1}^n\dsi^{\ka(i)}(\tau),
\end{align*}
where $\square_{i=1}^n$ is the product over $i=1,\ldots,n$ using
$\square$. So \eq{as5eq11} becomes
\e
\sum_{n=1}^\iy\frac{1}{n!}\cdot
\sum_{\substack{\text{$\ka:\{1,\ldots,n\}\ra C(\A):$}\\
\text{$\ka(\{1,\ldots,n\})\!=\!\al$,
$\tau\!\ci\!\ka\!\equiv\!\tau(\al)$}}}\!\!\!\!\!\!\!\!\!
\square_{i=1}^n\dsi^{\ka(i)}(\tau)=\dss^\al(\tau).
\label{as6eq7}
\e

To prove \eq{as6eq6}, fix $t\in T$, and consider the following
identity in~$\LCF(\fObj_\A)$:
\e
\de_{[0]}+\sum_{\al\in C(\A):\;\tau(\al)=t}\dss^\al(\tau)=
\de_{[0]}+\sum_{n\ge 1}\frac{1}{n!}
\raisebox{-5pt}{\begin{Large}$\displaystyle\Bigl[$\end{Large}}
\sum_{\be\in C(\A):\;\tau(\be)=t}\dsi^\be(\tau)
\raisebox{-5pt}{\begin{Large}$\displaystyle\Bigr]$
\end{Large}}^{\!\!\!\!\ts \square^n},
\label{as6eq8}
\e
where $f^{\square^n}$ for means $f\square f\square\cdots\square f$
with $f$ occurring $n$ times. All three sums in \eq{as6eq8} are
infinite, so we must explain what they mean.

One way to interpret \eq{as6eq8} is as a formal sum which packages
up finite identities in $\CF(\fObj_\A^\be)$ for each $\be\in C(\A)$.
It is easy to see that if $\tau(\be)\ne t$ then all terms in
\eq{as6eq8} are zero on $\fObj_\A^\be(\K)$, and if $\tau(\be)=t$
then the restriction of \eq{as6eq8} to $\fObj_\A^\be(\K)$ is exactly
\eq{as6eq7}, which has finitely many nonzero terms by Theorem
\ref{as5thm3}. This proves \eq{as6eq8} makes sense, and is true, as
such a finite formal sum. Another way to make sense of \eq{as6eq8}
is to use the ideas of~\cite[\S 4.2]{Joyc4}.

Now $\exp(x)=1+\sum_{n\ge 1}x^n/n!$, so \eq{as6eq8} may be rewritten
\begin{equation*}
\de_{[0]}+\sum_{\al\in C(\A):\;\tau(\al)=t}\dss^\al(\tau)=
\exp\raisebox{-5pt}{\begin{Large}$\displaystyle\Bigl[$\end{Large}}
\sum_{\be\in C(\A):\;\tau(\be)=t}\dsi^\be(\tau)
\raisebox{-5pt}{\begin{Large}$\displaystyle\Bigr]$\end{Large}}.
\end{equation*}
Formally taking logs and using $\log(1+x)=\sum_{n\ge
1}(-1)^{n-1}x^n/n$ gives
\e
\sum_{\be\in C(\A):\;\tau(\be)=t}\dsi^\be(\tau)=
\sum_{n\ge 1}\frac{(-1)^n}{n}
\raisebox{-5pt}{\begin{Large}$\displaystyle\Bigl[$\end{Large}}
\sum_{\al\in C(\A):\;\tau(\al)=t}\dss^\al(\tau)
\raisebox{-5pt}{\begin{Large}$\displaystyle\Bigr]$
\end{Large}}^{\!\!\!\!\square^n}.
\label{as6eq9}
\e

Here \eq{as6eq9} is interpreted in the same way as \eq{as6eq8}. It
follows from \eq{as6eq8} and $\log\ci\exp x=x$ as an identity in
formal power series. If $\tau(\be)\ne t$ then all terms in
\eq{as6eq9} are zero on $\fObj_\A^\be(\K)$, and if $\tau(\be)=t$
then the restriction of \eq{as6eq9} to $\fObj_\A^\be(\K)$ is
\eq{as6eq6}. So taking $t=\tau(\be)$ proves \eq{as6eq6}. The proof
for \eq{as5eq11} in Theorem \ref{as5thm3} shows there are only
finitely many nonzero terms in~\eq{as6eq6}.
\end{proof}

Following the proof of Theorem \ref{as5thm4}, but starting from
\eq{as6eq6} rather than \eq{as5eq11}, gives the following formula
for $\dsi(K,\tl,\mu,\tau)$. The only differences are in exchanging
$\dsi(\cdots),\dss(\cdots)$, and the combinatorial factors in the
last part.

\begin{thm} Let Assumption \ref{as3ass} hold, and\/
$(\tau,T,\le)$ be a permissible weak stability condition on $\A$.
Then for all\/ $\A$-data $(K,\tl,\mu)$ we have
\e
\begin{aligned}
\sum_{\substack{\text{iso.}\\ \text{classes}\\ \text{of finite }\\
\text{sets $I$}}}\!\! \frac{(-1)^{\md{I}-\md{K}}}{\md{I}!}\!\cdot
\!\!\!\!\!\!\!\!
\sum_{\substack{
\text{$\ka\!:\!I\!\ra\!C(\A)$, surjective $\phi\!:\!I\!\ra\!K$:}\\
\text{$\ka(\phi^{-1}(k))=\mu(k)$ for $k\in K$,}\\
\text{$\tau\ci\mu\ci\phi\equiv\tau\ci\ka:I\ra T$.}\\
\text{Define $\pr$ on $I$ by $i\pr j$ if $i=j$}\\
\text{or $\phi(i)\ne\phi(j)$ and $\phi(i)\tl\phi(j)$}}}
\begin{aligned}[t]
\ts\prod_{k\in K}\bigl(\md{\phi^{-1}(\{k\})}-1\bigr)!&\cdot\\
\CF^\stk\bigl(Q(I,\pr,K,\tl,\phi)\bigr)&\\
\dss(I,\pr,\ka,\tau)=&\\
\dsi(K,\tl,\mu,\tau)&.
\end{aligned}
\end{aligned}
\label{as6eq10}
\e
Only finitely many functions $\dss(I,\pr,\ka,\tau)$ in this sum are
nonzero.
\label{as6thm3}
\end{thm}

\subsection{Writing $\dssb(*,\tau)$ in terms of $\dstb(*,\tau)$}
\label{as63}

\begin{dfn} Let $(I,\pr)$ be a finite poset, $K$ a
finite set, and $\phi:I\ra K$ a surjective map. We call
$(I,\pr,K,\phi)$ {\it allowable} if there exists a partial order
$\tl$ on $K$ such that $i\pr j$ implies $\phi(i)\tl\phi(j)$. For
$(I,\pr,K,\phi)$ allowable, define a partial order $\ls$ on $K$ by
$k\ls l$ for $k,l\in K$ if there exist $b\ge 0$ and
$i_0,\ldots,i_b$, $j_0,\ldots,j_b$ in $I$ with $\phi(i_0)=k$,
$\phi(j_b)=l$, and $i_a\pr j_a$ for $a=0,\ldots,b$, and
$\phi(i_a)=\phi (j_{a-1})$ for $a=1,\ldots,b$. Write
$\cP(I,\pr,K,\phi)=\ls$. It has the property that if $\tl$ is a
partial order on $K$, then $i\pr j$ implies $\phi(i)\tl\phi(j)$ if
and only if $\tl$ dominates~$\cP(I,\pr,K,\phi)$.
\label{as6def2}
\end{dfn}

Here is a transitivity property of allowable quadruples.
The proof is elementary, and left as an exercise.

\begin{lem} Suppose $(I,\pr,J,\psi)$ is allowable with
$\ls=\cP(I,\pr,J,\psi)$, and\/ $\xi:J\ra K$ is a surjective map.
Then $(J,\ls,K,\xi)$ is allowable if and only if\/
$(I,\pr,K,\xi\!\ci\!\psi)$ is allowable, and when they
are~$\cP(J,\ls,K,\xi)\!=\!\cP(I,\pr,K,\xi\!\ci\!\psi)$.
\label{as6lem}
\end{lem}

We can now write $\dssb(*,\tau)$ in terms of $\dstb(*,\tau)$.

\begin{thm} Let Assumption \ref{as3ass} hold, $(\tau,T,\le)$
be a permissible stability condition on $\A$, and\/
$(J,\ls,\la)$ be $\A$-data. Then
\e
\begin{aligned}
\sum_{\substack{\text{iso. classes}\\ \text{of finite }\\
\text{sets $I$}}} \frac{1}{\md{I}!}\cdot
\sum_{\substack{
\text{$\pr,\ka,\psi$: $(I,\pr,\ka)$ is $\A$-data,}\\
\text{$(I,\pr,J,\psi)$ is allowable,}\\
\text{$\ls=\cP(I,\pr,J,\psi)$,}\\
\text{$\ka(\psi^{-1}(j))=\la(j)$ for $j\in J$,}\\
\text{$\tau\ci\la\ci\psi\equiv\tau\ci\ka:I\ra T$}}}
\begin{aligned}[t]
\CF^\stk\bigl(Q(I,\pr,J,\ls,\psi)\bigr)&\\
\dstb(I,\pr,\ka,\tau)&=\\
\dssb(J,\ls,\la&,\tau).
\end{aligned}
\end{aligned}
\label{as6eq11}
\e
Only finitely many functions $\dstb(I,\pr,\ka,\tau)$ in
this sum are nonzero.
\label{as6thm4}
\end{thm}

\begin{proof} Substituting \eq{as5eq32} with $J,\la$ in place
of $K,\mu$ into \eq{as6eq3} with $J,\la,\ls,\tl$ in place of
$I,\ka,\tl,\pr$ and using $Q(J,\tl,\ls)\ci Q(I,\pr,J,\tl,\psi)
=Q(I,\pr,J,\ls,\psi)$, Theorems \ref{as2thm1}, \ref{as5thm7} and
\ref{as6thm1} and Definition \ref{as6def2} gives
\begin{gather*}
\sum_{\substack{\text{iso. classes}\\ \text{of finite }\\
\text{sets $I$}}} \frac{1}{\md{I}!}\cdot
\sum_{\substack{
\text{$\pr,\ka,\psi$: $(I,\pr,\ka)$ is $\A$-data,}\\
\text{$(I,\pr,J,\psi)$ is allowable,}\\
\text{$\ka(\psi^{-1}(j))=\la(j)$ for $j\in J$,}\\
\text{$\tau\ci\la\ci\psi\equiv\tau\ci\ka:I\ra T$}}}
\raisebox{-12pt}{\begin{Large}$\displaystyle\Biggl[$\end{Large}}
\sum_{\substack{\text{partial orders $\tl$ on $J$:}\\
\text{$\ls$ dominates $\tl$,}\\
\text{$\tl$ dominates $\cP(I,\pr,J,\psi)$}}} \!\!\!\!\!\!
n(J,\tl,\ls)
\raisebox{-12pt}{\begin{Large}$\displaystyle\Biggr]$\end{Large}}
\cdot\\
\CF^\stk\bigl(Q(I,\pr,J,\ls,\psi)\bigr)\dstb(I,\pr,\ka,\tau)=
\dssb(J,\ls,\la,\tau),
\end{gather*}
with only finitely many $\dstb(I,\pr,\ka,\tau)$ nonzero. By
\eq{as6eq2} the bracketed sum is 1 if $\ls=\cP(I,\pr,J,\psi)$ and 0
otherwise, and \eq{as6eq11} follows.
\end{proof}

\subsection{Inverting \eq{as6eq11}}
\label{as64}

Next we invert \eq{as6eq11}. We will need the following notation.

\begin{dfn} Let $I$ be a finite set. Then equivalence
relations $\sim$ on $I$ are in 1-1 correspondence with
subsets $S=\bigl\{(i,j)\in I\t I:i\sim j\bigr\}$ of $I\t I$
satisfying the properties (i) $(i,i)\in S$ for all $i\in I$,
\begin{itemize}
\setlength{\itemsep}{0pt}
\setlength{\parsep}{0pt}
\item[(ii)] $(i,j)\in S$ implies $(j,i)\in S$, and
\item[(iii)] $(i,j)\in S$ and $(j,k)\in S$ imply~$(i,k)\in S$.
\end{itemize}
Given $S\subseteq I\t I$ satisfying (i)--(iii), define an
equivalence relation $\sim_S$ on $I$ by $i\,\sim_S\,j$ if
$(i,j)\in S$. Write $[i]_S$ for the $\sim_S$-equivalence
class of $i$, set $I_S=\{[i]_S:i\in I\}$, and define $\psi_S:
I\ra I_S$ by~$\psi_S(i)=[i]_S$.

Now let $(I,\pr)$ be a finite poset, and define
\begin{equation*}
\U(I,\pr)\!=\!\bigl\{S\!\subseteq\!I\!\t\!I:\text{$S$ satisfies
(i)--(iii), $(I,\pr,I_S,\psi_S)$ is allowable}\bigr\}.
\end{equation*}
Suppose $(I,\pr,K,\phi)$ is allowable, and define
$S=\bigl\{(i,j)\in I\t I:\phi(i)=\phi(j)\bigr\}$.
Then it is easy to see that $S\in\U(I,\pr)$, and
there is a unique 1-1 correspondence
$\io:I_S\ra K$ with $\io\bigl([i]_S\bigr)=\phi(i)$
for $i\in I$ such that $\phi=\io\ci\psi_S$. So
$\U(I,\pr)$ classifies isomorphism classes of $K,\phi$
such that $(I,\pr,K,\phi)$ is allowable. Define
\begin{equation*}
N(I,\pr)=\sum_{\substack{\text{$n\ge 0$,
$S_0,\ldots,S_n\in\U(I,\pr)$:}\\
\text{$S_{m-1}\subset S_m$, $S_{m-1}\ne S_m$, $1\le m\le n$}\\
\text{$S_0=\{(i,i):i\in I\}$, $S_n=I\t I$}}}(-1)^n.
\end{equation*}

Now let $(I,\pr,K,\phi)$ be allowable, and define
\begin{equation*}
N(I,\pr,K,\phi)=\sum_{\substack{\text{$n\ge 0$,
$S_0,\ldots,S_n\in\U(I,\pr)$:}\\
\text{$S_{m-1}\subset S_m$, $S_{m-1}\ne S_m$, $1\le m\le n$}\\
\text{$S_0=\{(i,i):i\in I\}$,}\\
\text{$S_n=\{(i,j)\in I\t I:\phi(i)=\phi(j)\}$}}}(-1)^n.
\end{equation*}
\label{as6def3}
\end{dfn}

By a similar proof to Proposition \ref{as6prop1}, using Lemma
\ref{as6lem}, we can show:

\begin{prop} Let\/ $(I,\pr,K,\phi)$ be allowable. Then
\e
\begin{gathered}
\sum_{\substack{\text{iso. classes}\\ \text{of finite }\\
\text{sets $J$}}} \frac{1}{\md{J}!}\cdot
\sum_{\substack{\text{$\psi:I\!\ra\!J$, $\xi:J\!\ra\!K$}\\
\text{surjective, $\phi\!=\!\xi\!\ci\!\psi$:}\\
\text{$(I,\pr,J,\psi)$ allowable}}}
\!\!\!\!\!\!\!\!\! N(I,\pr,J,\psi)
=\begin{cases} 1, & \text{$\phi$ is a bijection,}
\\ 0, & \text{otherwise.}\end{cases}
\end{gathered}
\label{as6eq12}
\e
This also holds with\/ $N(I,\pr,J,\psi)$ replaced
by~$N\bigl(J,\cP(I,\pr,J,\psi),K,\xi\bigr)$.
\label{as6prop2}
\end{prop}

Here is a product formula for $N(I,\pr,K,\phi)$. We leave the proof
as an exercise; one possible starting point is to note that both
sides of \eq{as6eq13} satisfy~\eq{as6eq12}.

\begin{prop} Let\/ $(I,\pr,K,\phi)$ be allowable. Then
\e
N(I,\pr,K,\phi)=\ts\prod_{k\in K}N\bigl(\phi^{-1}(k),
\pr\vert_{\phi^{-1}(k)}\bigr).
\label{as6eq13}
\e
\label{as6prop3}
\end{prop}

Next we invert the identity of Theorem~\ref{as6thm4}.

\begin{thm} Let Assumption \ref{as3ass} hold, $(\tau,T,\le)$
be a permissible stability condition on $\A$, and\/
$(K,\tl,\mu)$ be $\A$-data. Then
\e
\begin{aligned}
\sum_{\substack{\text{iso. classes}\\ \text{of finite }\\
\text{sets $J$}}} \frac{1}{\md{J}!}\cdot \!\!\!\!\!
\sum_{\substack{
\text{$\ls,\la,\chi$: $(J,\ls,\la)$ is $\A$-data,}\\
\text{$(J,\ls,K,\chi)$ is allowable,}\\
\text{$\tl=\cP(J,\ls,K,\chi)$,}\\
\text{$\la(\chi^{-1}(k))=\mu(k)$ for $k\in K$,}\\
\text{$\tau\ci\mu\ci\chi\equiv\tau\ci\la:J\ra T$}}}
\!\!\!\!\!\!\!\!\!\!\!\!\!\!\!\!\!\!\!\!\!\!\!\!
\begin{aligned}[t]
N(J,\ls,K,\chi)\CF^\stk\bigl(Q(J,\ls,K,\tl,\chi&)\bigr)\\
\dssb(J,\ls,\la,\tau)&=\\
\dstb(K,\tl,\mu&,\tau).
\end{aligned}
\end{aligned}
\label{as6eq14}
\e
Only finitely many functions $\dssb(J,\ls,\la,\tau)$ in
this sum are nonzero.
\label{as6thm5}
\end{thm}

\begin{proof} Using the proof in \S\ref{as52} that only
finitely many $\dsi(I,\pr,\ka,\tau)$ in \eq{as5eq14} are nonzero, we
find that only finitely many $\dssb(J,\ls,\la,\tau)$ in \eq{as6eq14}
are nonzero. Substituting \eq{as6eq11} into the left hand side of
\eq{as6eq14} gives
\begin{gather}
\sum_{\substack{\text{iso. classes}\\ \text{of finite }\\
\text{sets $J$}}} \frac{1}{\md{J}!}\cdot
\sum_{\substack{
\text{$\ls,\la,\chi$: $(J,\ls,\la)$ is $\A$-data,}\\
\text{$(J,\ls,K,\chi)$ is allowable,}\\
\text{$\tl=\cP(J,\ls,K,\chi)$,}\\
\text{$\la(\chi^{-1}(k))=\mu(k)$ for $k\in K$,}\\
\text{$\tau\ci\mu\ci\chi\equiv\tau\ci\la:J\ra T$}}}
\!\!\!\!\!\!
N(J,\ls,K,\chi)\cdot
\nonumber \\
\sum_{\substack{\text{iso. classes}\\ \text{of finite }\\
\text{sets $I$}}} \frac{1}{\md{I}!}\cdot
\sum_{\substack{
\text{$\pr,\ka,\psi$: $(I,\pr,\ka)$ is $\A$-data,}\\
\text{$(I,\pr,J,\psi)$ is allowable,}\\
\text{$\ls=\cP(I,\pr,J,\psi)$,}\\
\text{$\ka(\psi^{-1}(j))=\la(j)$ for $j\in J$,}\\
\text{$\tau\ci\la\ci\psi\equiv\tau\ci\ka:I\ra T$}}}
\begin{aligned}[t]
\CF^\stk\bigl(Q(J,\ls,K,\tl,\chi)\bigr)&\\
\Bigl[\CF^\stk\bigl(Q(I,\pr,J,\ls,\psi)\bigr)&\\
\dstb(I,\pr,\ka,\tau)\Bigr]&=
\end{aligned}
\nonumber
\allowdisplaybreaks \\
\sum_{\substack{\text{iso. classes}\\ \text{of finite }\\
\text{sets $I$}}} \frac{1}{\md{I}!}\cdot
\sum_{\substack{
\text{$\pr,\ka,\phi$: $(I,\pr,\ka)$ is $\A$-data,}\\
\text{$(I,\pr,K,\phi)$ is allowable,}\\
\text{$\tl=\cP(I,\pr,K,\phi)$,}\\
\text{$\ka(\phi^{-1}(k))=\mu(k)$ for $k\in K$,}\\
\text{$\tau\ci\mu\ci\phi\equiv\tau\ci\ka:I\ra T$}}}
\begin{aligned}[t]
\CF^\stk\bigl(Q(I,\pr,K,\tl,\phi)\bigr)&\\
\dstb(I,\pr,\ka,\tau)&\cdot
\end{aligned}
\nonumber\\
\raisebox{-12pt}{\begin{Large}$\displaystyle\Biggl[$\end{Large}}
\sum_{\substack{\text{iso. classes}\\ \text{of finite }\\
\text{sets $J$}}} \frac{1}{\md{J}!}\cdot
\sum_{\substack{\text{$\psi:I\!\ra\!J$, $\chi:J\!\ra\!K$}\\
\text{surjective, $\phi\!=\!\chi\!\ci\!\psi$:}\\
\text{$(I,\pr,J,\psi)$ allowable,}\\
\text{$\ls=\cP(I,\pr,J,\psi)$}}} \!\!\!\!\!\!\!\!\! N(J,\ls,K,\chi)
\raisebox{-12pt}{\begin{Large}$\displaystyle\Biggr]$\end{Large}},
\label{as6eq15}
\end{gather}
setting $\phi=\chi\ci\psi$ and using $\CF^\stk(Q(I,\pr,K,\tl,\phi))=\CF^\stk
(Q(J,\ls,K,\tl,\chi))\ci\CF^\stk(Q(I,\pr,J,\ls,\psi))$ in the third line.

Here, given $(I,\pr,J,\psi)$ allowable and $\ls=\cP(I,\pr,J,\psi)$,
Lemma \ref{as6lem} shows that $(J,\ls,K,\chi)$ allowable and
$\tl=\cP(J,\ls,K,\chi)$ in the first line of \eq{as6eq15} is
equivalent to $(I,\pr,K,\phi)$ allowable and $\tl=\cP(I,\pr,K,\phi)$
in the third line. Also, $\tau\ci\mu\ci\chi\equiv\tau\ci\la$ and
$\tau\ci\la\ci\psi\equiv \tau\ci\ka$ in the first and second lines
of \eq{as6eq15} are equivalent to
$\tau\ci\mu\ci\phi\equiv\tau\ci\ka$ in the third,
as~$\phi=\chi\ci\psi$.

Now Proposition \ref{as6prop2} shows that the bracketed term on the
last line of \eq{as6eq15} is 1 if $\phi$ is a bijection, and 0
otherwise. When $\phi$ is a bijection $\md{I}=\md{K}$. The first sum
on the third line in \eq{as6eq15} fixes a unique $I$ with
$\md{I}=\md{K}$. Then in the second sum there are $\md{I}!$
bijections $\phi:I\ra K$. So by dropping the factor $1/\md{I}!$ on
the third line we may take $I=K$ and $\phi=\id_K$. Then $\tl=\pr$,
$\ka=\mu$, and $\CF^\stk(Q(K,\tl,K,\tl,\id_K))$ is the identity.
Thus, the last two lines of \eq{as6eq15} reduce to
$\dstb(K,\tl,\mu,\tau)$, the right hand side of \eq{as6eq14}. This
completes the proof.
\end{proof}

\section{(Lie) algebras of constructible functions}
\label{as7}

We now define and study some interesting subalgebras $\Hp_\tau,
\Ht_\tau$ of $\CF(\fObj_\A)$, for $(\tau,T,\le)$ a permissible weak
stability condition. These encode information about the moduli
spaces $\Oss^\al,\Osi^\al,\Ost^\al(\tau)$ for all $\al\in C(\A)$. We
will see in \cite{Joyc5} that these subalgebras are essentially {\it
independent of choice of\/} $(\tau,T,\le)$, and that changing weak
stability conditions amounts to changing bases in
$\Hp_\tau,\Ht_\tau$. We suppose $\K$ has characteristic zero
throughout this section.

\subsection{The algebras $\Hp_\tau,\Ht_\tau$ and Lie algebras
$\Lp_\tau,\Lt_\tau$}
\label{as71}

From \S\ref{as33}, $\CF(\fObj_\A)$ is a $\Q$-algebra with
multiplication $*$ and identity $\de_{[0]}$. Given a permissible
weak stability condition $(\tau,T,\le)$ we define two interesting
subalgebras $\Hp_\tau,\Ht_\tau$ of~$\CF(\fObj_\A)$.

\begin{dfn} Let Assumption \ref{as3ass} hold, and $(\tau,T,\le)$ be
a permissible weak stability condition on $\A$. Define $\Q$-vector
subspaces $\Hp_\tau,\Ht_\tau$ in $\CF(\fObj_\A)$ by
\ea
\Hp_\tau&=\bigl\langle\CF^\stk(\bs\si(I))\dss(I,\pr,
\ka,\tau):\text{$(I,\pr,\ka)$ is $\A$-data}\bigr\rangle{}_\Q,
\label{as7eq1}\\
\Ht_\tau&=\bigl\langle\de_{[0]},
\dss^{\al_1}(\tau)*\cdots*\dss^{\al_n}(\tau):\al_1,\ldots,\al_n\in
C(\A)\bigr\rangle{}_\Q.
\label{as7eq2}
\ea
Here $\langle\cdots\rangle_\Q$ is the set of all finite $\Q$-linear
combinations of the elements `$\,\cdots$', and $\A$-data is defined
in Definition \ref{as3def5}. Define
$\Lp_\tau=\Hp_\tau\cap\CFi(\fObj_\A)$
and~$\Lt_\tau=\Ht_\tau\cap\CFi(\fObj_\A)$.
\label{as7def1}
\end{dfn}

In \cite{Joyc5} we will show that if $(\tau,T,\le)$ and
$(\ti\tau,\ti T,\le)$ are permissible weak stability conditions on
$\A$, then (under some finiteness conditions)
$\Hp_\tau=\Hp_{\ti\tau}$ and $\Ht_\tau=\Ht_{\ti\tau}$, so that
$\Hp_\tau,\Ht_\tau$ are {\it independent of the choice of\/}
$(\tau,T,\le)$. To relate $\Ht_\tau$ and $\Hp_\tau$, let
$(\{1,\ldots,n\},\le,\ka)$ be $\A$-data. Then
\e
\dss(\{1,\ldots,n\},\le,\ka,\tau)=\ts\prod_{i=1}^n
\bigl(\bs\si(\{i\})^*(\dss^{\ka(i)}(\tau))\bigr).
\label{as7eq3}
\e
Generalizing the argument of \cite[Th.~4.3]{Joyc4} we then find that
\e
\dss^{\ka(1)}(\tau)*\!\cdots\!*\dss^{\ka(n)}(\tau)
=\CF^\stk(\bs\si(\{1,\ldots,n\}))\dss(\{1,\ldots,n\},\le,\ka,\tau).
\label{as7eq4}
\e
Thus $\Hp_\tau$ is the span of $\CF^\stk(\bs\si(I))
\dss(I,\pr,\ka,\tau)$ for $\A$-data $(I,\pr,\ka)$ with $\pr$ a {\it
partial order}, and $\Ht_\tau$ the span with $\pr$ a {\it total
order}. This explains the notation.

Now in \cite[\S 4.8]{Joyc4} we defined multilinear operations
$P_\sIp$ on $\CF(\fObj_\A)$ for $(I,\pr)$ a finite poset, and
generalizing \eq{as7eq3} shows that
\e
\CF^\stk(\bs\si(I))\dss(I,\pr,\ka,\tau)=P_\sIp\bigl(\dss^{\ka(i)}
(\tau):i\in I\bigr).
\label{as7eq5}
\e
Thus an alternative expression for $\Hp_\tau$ is
\e
\Hp_\tau=\bigl\langle P_\sIp\bigl(\dss^{\ka(i)}(\tau):i\in I\bigr):
\text{$(I,\pr,\ka)$ is $\A$-data}\bigr\rangle{}_\Q.
\label{as7eq6}
\e
It follows from \cite[Th.~4.22]{Joyc4} that $\Hp_\tau$ is closed
under the operations~$P_\sIp$.

\begin{prop} In Definition \ref{as7def1}, $\Hp_\tau,\Ht_\tau$ are
subalgebras of\/ $\CF(\fObj_\A)$ and\/ $\Lp_\tau,\Lt_\tau$ Lie
subalgebras of\/ $\CFi(\fObj_\A)$, with\/
$\Ht_\tau\!\subseteq\!\Hp_\tau$
and\/~$\Lt_\tau\!\subseteq\!\Lp_\tau$.
\label{as7prop1}
\end{prop}

\begin{proof} Clearly $\Ht_\tau$ is the subalgebra of $\CF(\fObj_\A)$
generated by the $\dss^\al(\tau)$ for all $\al\in C(\A)$. As
$\Hp_\tau$ is closed under the operations $P_\sIp$ from above, it is
closed under $*=P_{\sst(\{1,2\},\le)}$. Writing
$(\emptyset,\emptyset,\emptyset)$ for the trivial $\A$-data we have
$\CF^\stk(\bs\si(\emptyset))\dss (\emptyset,\emptyset,
\emptyset,\tau)=\de_{[0]}$, so $\Hp_\tau$ contains the identity
$\de_{[0]}$, and is a subalgebra of $\CF(\fObj_\A)$. Therefore
$\Lp_\tau,\Lt_\tau$ are intersections of Lie subalgebras
$\Hp_\tau,\Ht_\tau$ and $\CFi(\fObj_\A)$ of $\CF(\fObj_\A)$, so they
are Lie subalgebras. The inclusion $\Ht_\tau\subseteq \Hp_\tau$ is
obvious from \eq{as7eq4}, and this
implies~$\Lt_\tau\subseteq\Lp_\tau$.
\end{proof}

We now apply the work of \S\ref{as5}--\S\ref{as6} to study
$\Hp_\tau$. There we constructed eleven transformations between the
six families $\dss,\dsi,\dst,\dssb,\dsib,\dstb (*,\tau)$. Their
equation numbers are displayed below.
\e
\text{
\begin{footnotesize}
$\displaystyle \!\!\! \xymatrix@C=18pt{\dsib(*,\tau)
\ar@/^1.2pc/[r]_{\ts\eq{as5eq8}} & \dsi(*,\tau)
\ar@/^1.2pc/[r]_{\ts\eq{as5eq14}} \ar@/^1.2pc/[l]_{\ts\eq{as6eq4}} &
\dss(*,\tau) \ar@/^1.2pc/[r]_{\ts\eq{as6eq3}}
\ar@/^1.2pc/[l]_{\ts\eq{as6eq10}} & \dssb(*,\tau)
\ar@/^1.2pc/[l]_{\ts\eq{as5eq7}} \ar@/^1.2pc/[r]_{\ts\eq{as6eq14}} &
\dstb(*,\tau) \ar@/^1.2pc/[r]_{\ts\eq{as5eq9}}
\ar@<-1ex>@/_2pc/[ll]^{\ts\eq{as5eq32}}
\ar@/^1.2pc/[l]_{\ts\eq{as6eq11}} & \dst(*,\tau).
\ar@/^1.2pc/[l]_{\ts\eq{as6eq5}} } \!\!\! $
\end{footnotesize}
}
\label{as7eq7}
\e
Note that the identities involving $\dst,\dstb(*,\tau)$ require
$(\tau,T,\le)$ to be a stability condition, but the other identities
work for $(\tau,T,\le)$ a {\it weak\/} stability condition.
Combining these, we can write any of the six families
$\dss,\dsi,\dst,\dssb,\dsib,\dstb(*,\tau)$ in terms of any of the
others. Applying $\CF^\stk(\bs\si(I))$ to \eq{as5eq7}, noting that
$\bs\si(I)\ci Q(I,\pr,\tl)=\bs\si(I)$ and using \eq{as2eq3} yields
\e
\sum_{\substack{\text{p.o.s $\pr$ on $I$:}\\
\text{$\tl$ dominates $\pr$}}}
\CF^\stk\bigl(\bs\si(I)\bigr)\dssb(I,\pr,\ka,\tau)=
\CF^\stk\bigl(\bs\si(I)\bigr)\dss(I,\tl,\ka,\tau).
\label{as7eq8}
\e
Similarly, all eleven transformations \eq{as7eq7} imply
transformations between the six families $\CF^\stk
(\bs\si(I))\dss,\ldots,\dstb(*,\tau)$ in $\CF(\fObj_\A)$. Thus we
deduce:

\begin{cor} In Definition \ref{as7def1} we have
\e
\begin{split}
\Hp_\tau&=\bigl\langle\CF^\stk(\bs\si(I))\dssb(I,\pr,
\ka,\tau):\text{$(I,\pr,\ka)$ is $\A$-data}\bigr\rangle{}_\Q\\
&=\bigl\langle\CF^\stk(\bs\si(I))\dsi(I,\pr,
\ka,\tau):\text{$(I,\pr,\ka)$ is $\A$-data}\bigr\rangle{}_\Q\\
&=\bigl\langle\CF^\stk(\bs\si(I))\dsib(I,\pr,
\ka,\tau):\text{$(I,\pr,\ka)$ is $\A$-data}\bigr\rangle{}_\Q\\
&=\bigl\langle\CF^\stk(\bs\si(I))\dst(I,\pr,
\ka,\tau):\text{$(I,\pr,\ka)$ is $\A$-data}\bigr\rangle{}_\Q\\
&=\bigl\langle\CF^\stk(\bs\si(I))\dstb(I,\pr,
\ka,\tau):\text{$(I,\pr,\ka)$ is $\A$-data}\bigr\rangle{}_\Q,
\end{split}
\label{as7eq9}
\e
supposing $(\tau,T,\le)$ is a stability condition in the last two
lines.
\label{as7cor1}
\end{cor}

The material of \S\ref{as5}--\S\ref{as6}, and other identities in
\cite{Joyc5}, can therefore be interpreted as giving {\it basis
change formulae} in the infinite-dimensional algebra $\Hp_\tau$. In
particular, $\Hp_\tau$ contains $\dss^\al,\dsi^\al,\dst^\al(\tau)$
for all $\al\in C(\A)$. We can interpret this as saying that
$\Hp_\tau$ contains information about $\tau$-semistability,
$\tau$-semistable indecomposables, and $\tau$-stability, but
$\Ht_\tau$ only information about $\tau$-semistability.

We can write down the multiplication relations in $\Hp_\tau$
explicitly for the six spanning sets $\CF^\stk(\bs\si(I))\dss,
\ldots,\dstb(I,\pr, \ka,\tau)$. Let $(I,\pr,\ka)$ and $(J,\ls,\la)$
be $\A$-data with $I\cap J=\emptyset$. Define $\A$-data
$(K,\tl,\mu)$ by $K=I\amalg J$, $\mu\vert_I=\ka$, $\mu\vert_J=\la$,
and $k\tl l$ if either $k,l\in I$ and $k\pr l$, or $k,l\in J$ and
$k\ls l$, or $k\in I$ and $l\in J$. Then from \cite[Th.~4.22]{Joyc4}
we deduce that
\e
\begin{aligned}
\bigl(\CF^\stk(\bs\si(I))\dss(I,\pr,\ka,\tau)\bigr)\,*&\,
\bigl(\CF^\stk(\bs\si(J))\dss(J,\ls,\la,\tau)\bigr)\\
=&\,\CF^\stk(\bs\si(K))\dss(K,\tl,\mu,\tau).
\end{aligned}
\label{as7eq10}
\e
The same holds with $\dsi(*)$ or $\dst(*)$ in place of $\dss(*)$.
Using \eq{as5eq7}--\eq{as5eq9} and \eq{as6eq3}--\eq{as6eq5} we can
now deduce the multiplication relations for the $\CF^\stk(\bs\si(I))
\dssb,\ab\dsib,\ab\dstb(I,\pr,\ka,\tau)$, and the answer turns out
as follows. Let $(I,\pr,\ka),(J,\ls,\la),K$ and $\mu$ be as above,
but do not define $\tl$. Then
\e
\begin{aligned}
\bigl(\CF^\stk(\bs\si(I))\dssb(I,\pr,\ka,\tau)\bigr)\,*\,
\bigl(\CF^\stk(\bs\si(J))&\dssb(J,\ls,\la,\tau)\bigr)\\
=\sum_{\substack{\text{p.o.s $\tl$ on $K$: $\tl\vert_I=\pr$,
$\tl\vert_J=\ls$}\\ \text{and $i\in I$, $j\in J$ implies $j\ntl
i$}}} \CF^\stk(\bs\si(K))&\dssb(K,\tl,\mu,\tau).
\end{aligned}
\label{as7eq11}
\e
The same holds with $\dsib(*)$ or $\dstb(*)$ in place of~$\dssb(*)$.

\subsection{The structure of the Lie algebra $\Lp_\tau$}
\label{as72}

If $(I,\pr)$ is a finite poset, let $\approx$ be the equivalence
relation on $I$ generated by $i\approx j$ if $i\pr j$ or $j\pr i$,
and define the {\it connected components} of $(I,\pr)$ to be the
$\approx$-equivalence classes. Equivalently, if $\Ga$ is the
directed graph with vertices $I$ and edges
$\smash{\mathop{\bu}\limits^{\sst i} \ra\mathop{\bu}\limits^{\sst
j}}$ for $i,j\in I$ with $i\pr j$, then the connected components of
$(I,\pr)$ are the sets of vertices of connected components of $\Ga$.
We call $(I,\pr)$ {\it connected\/} if it has exactly one connected
component. Then we prove:

\begin{prop} Let Assumption \ref{as3ass} hold, $(\tau,T,\le)$
be a permissible weak stability condition on $\A$, and\/
$(I,\pr,\ka)$ be $\A$-data. If\/ $(I,\pr)$ has $k$ connected
components, then $\CF^\stk(\bs\si(I))\dsib(I,\pr,\ka,\tau)$ and
$\CF^\stk(\bs\si(I))\dstb(I,\pr,\ka,\tau)$ are supported on points
$[X_1\op\cdots\op X_k]\in\fObj_\A(\K)$, with all\/ $X_a$
indecomposable.
\label{as7prop2}
\end{prop}

\begin{proof} Let $I_1,\ldots,I_k\subseteq I$ be the connected
components of $(I,\pr)$, so that $I=I_1\amalg\cdots\amalg I_k$.
Suppose $\CF^\stk(\bs\si(I))\dsib(I,\pr,\ka,\tau)$ is nonzero on
$[X]\in\fObj_\A(\K)$. Then there exists $[(\si,\io,\pi)]\in
\Msib(I,\pr,\ka)_\A$ with $\si(I)=X$ making a nonzero contribution
to $\CF^\stk(\bs\si(I))\dsib(I,\pr,\ka,\tau)$ at $[X]$. The $I_1,
\ldots,I_k$ are s-sets in $I$, so setting $X_a=\si(I_a)$ we find
$\io(I_a,I):X_a\ra X$ defines a subobject $X_a\subset X$, with
$X=X_1\op\cdots\op X_k$. We shall prove $X_a$ is indecomposable for
$a=1,\ldots,k$. Write $(\si_a,\io_a,\pi_a)$ for the
$(I_a,\pr)$-subconfiguration of $(\si,\io,\pi)$. Then
$\si_a(I_a)=X_a$, and as $(\si,\io,\pi)$ is best Theorem
\ref{as3thm1} implies $(\si_a,\io_a,\pi_a)$ is best.

Let $T_a$ be a maximal torus in $\Aut(\si_a,\io_a,\pi_a)$ containing
$\bigl\{\la\id_{(\si_a,\io_a,\pi_a)}:\la\in\K^\t\bigr\}$. Then
$\bs\si(\{i\})_*(T_a)$ is a $\K$-subtorus of $\Aut(\si_a(\{i\}))$
containing $\bigl\{\la\id_{\si_a(\{i\})}: 0\ne\la\in\K\bigr\}$. But
since $\si_a(\{i\})=\si(\{i\})$ is {\it indecomposable} for $i\in
I_a$, $\Aut(\si_a(\{i\}))$ has {\it rank one}, so
$\bigl\{\la\id_{\si_a(\{i\})}:0\ne\la\in\K\bigr\}$ is a maximal
torus of $\Aut(\si_a(\{i\}))$. Since $\bs\si(\{i\})_*(T_a)$ must be
contained in a maximal torus of $\Aut(\si_a(\{i\}))$, we see that
\e
\bs\si(\{i\})_*(T_a)=\bigl\{\la\id_{\si_a(\{i\})}:0\ne\la\in\K\bigr\}.
\label{as7eq12}
\e

We claim that
\ea
&\bigl(\ts\prod_{i\in
I_a}\bs\si(\{i\})_*\bigr)(T_a)=\bigl\{\ts\prod_{i\in
I_a}\la\id_{\si(\{i\})}:0\ne\la\in\K\bigr\}\cong\K^\t,
\label{as7eq13}\\
\text{where}\quad &\ts\prod_{i\in
I_a}\bs\si(\{i\})_*:\Aut(\si_a,\io_a,\pi_a)\longra\ts\prod_{i\in
I_a}\Aut\bigl(\si_a(\{i\})\bigr).
\label{as7eq14}
\ea
The right hand side of \eq{as7eq13} is the image of
$\{\la\id_{(\si_a,\io_a,\pi_a)}:\la\in\K^\t\}\subseteq T_a$, so the
left hand side of \eq{as7eq13} contains the right. To prove the
opposite inclusion, let $\al\in T_a$. Then for each $i\in I_a$,
equation \eq{as7eq12} gives $\al(\{i\})=\la_i\id_{\si(\{i\})}$ for
some $\la_i\in\K^\t$. We must show $\la_i=\la$ for some
$\la\in\K^\t$ and all~$i\in I_a$.

Suppose $i\ne j\in I_a$ with $i\pr j$ but there is no $k\in I_a$
with $i\ne k\ne j$ and $i\pr k\pr j$. By Theorem \ref{as3thm1} the
short exact sequence \eq{as3eq1} is not split, and so corresponds to
a {\it nonzero} $\ga_{ij}\in\Ext^1(\si(\{j\}),\si(\{i\}))$. But
$\al$ induces an automorphism of \eq{as3eq1}, so
$\ga_{ij}\ci\la_j\id_{\si(\{j\})}
=\la_i\id_{\si(\{i\})}\ci\ga_{ij}$, giving $\la_i=\la_j$ as
$\ga_{ij}\ne 0$. Since $(I_a,\pr)$ is {\it connected\/} there are
enough such pairs $i,j$ to force $\la_i=\la$ for all $i\in I_a$.
This proves~\eq{as7eq13}.

If $(\si',\io',\pi')$ is a $(\{1,2\},\le)$-configuration then the
kernel of $\bs\si(\{1\})\t\bs\si(\{2\}):\Aut(\si',\io',\pi')\ra\Aut
(\si'(\{1\}))\t\Aut(\si'(\{2\})$ is $\Hom\bigl(
\si'(\{2\}),\si'(\{1\})\bigr)$. Generalizing this, one can show by
induction on $\md{I_a}$ that the kernel of \eq{as7eq14} is a {\it
nilpotent\/} $\K$-group. Thus, \eq{as7eq14} is injective on the
maximal torus $T_a$, and \eq{as7eq13} implies that $T_a=\bigl\{
\la\id_{(\si_a,\io_a,\pi_a)}:\la\in\K^\t\bigr\}$, so
$\Aut(\si_a,\io_a,\pi_a)$ has rank one.

The contribution of $[(\si,\io,\pi)]$ to $\CF^\stk(\bs\si(I))\dsib
(I,\pr,\ka,\tau)$ at $[X]$, which is nonzero by assumption, is
\e
\ts\prod_{a=1}^k\chi\bigl(\Aut(X_a)/\bs\si(I_a)_*
(\Aut(\si_a,\io_a,\pi_a))\bigr).
\label{as7eq15}
\e
Suppose $\Aut(X_a)$ has rank greater than one, and consider the
action of a maximal torus of $\Aut(X_a)$ on $\Aut(X_a)/\bs\si(I_a)_*
(\Aut(\si_a,\io_a,\pi_a))$. Since $\Aut(\si_a,\io_a,\pi_a)$ has rank
one, the orbits of this action are all of the form $(\K^\t)^l$ for
$l\ge 1$, which implies that the Euler characteristic in
\eq{as7eq15} is zero, a contradiction. Thus $\Aut(X_a)$ has rank
one, and $X_a$ is indecomposable for $a=1,\ldots,k$, as we have to
prove. Since $\tau$-stable objects are indecomposable, the same
proof works for~$\CF^\stk(\bs\si(I))\dstb(I,\pr,\ka,\tau)$.
\end{proof}

We can now deduce an alternative description of~$\Lp_\tau$.

\begin{prop} In Definition \ref{as7def1} we have
\e
\begin{split}
\Lp_\tau\!&=\!\bigl\langle\CF^\stk(\bs\si(I))\dsib(I,\pr,\ka,\tau):
\text{$(I,\pr,\ka)$ $\A$-data, $(I,\pr)$
connected\/}\bigr\rangle{}_\Q\\
&=\!\bigl\langle\CF^\stk(\bs\si(I))\dstb(I,\pr,\ka,\tau):
\text{$(I,\pr,\ka)$ $\A$-data, $(I,\pr)$
connected\/}\bigr\rangle{}_\Q,
\end{split}
\label{as7eq16}
\e
supposing $(\tau,T,\le)$ is a stability condition in the second
line. There is a natural\/ $\Q$-algebra isomorphism $\Phi^{\rm
pa}_\tau:U(\Lp_\tau)\ra\Hp_\tau$, where $U(\Lp_\tau)$ is the
universal enveloping algebra of\/~$\Lp_\tau$.
\label{as7prop3}
\end{prop}

\begin{proof} Equation \eq{as7eq16} follows from Definition
\ref{as7def1}, \eq{as7eq9} and Proposition \ref{as7prop2}. From
above, the multiplication relations for the
$\CF^\stk(\bs\si(I))\dsib(I,\pr,\ka,\tau)$ are given by \eq{as7eq11}
with $\dsib(*)$ in place of $\dssb(*)$. From this it is easy to see
that if $(I,\pr)$ has connected components $I_1,\ldots,I_k$ then
\begin{gather*}
\CF^\stk(\bs\si(I))\dsib(I,\pr,\ka,\tau)=
\bigl(\CF^\stk(\bs\si(I_1))\dsib(I_1,\pr,\ka,\tau)\bigr)*\cdots\\
*\bigl(\CF^\stk(\bs\si(I_k))\dsib(I_k,\pr,\ka,\tau)\bigr)
+\bigl(\text{$\Q$-linear combination of}\\
\text{$\CF^\stk(\bs\si(J))\dsib(J,\ls,\la,\tau)$ for $(J,\ls)$ with
$<k$ connected components}\bigr).
\end{gather*}
Then \eq{as7eq16} and induction on $k$ shows that
$\CF^\stk(\bs\si(I))\dsib(I,\pr,\ka,\tau)$ is contained in the
algebra generated by $\Lp_\tau$ for all $(I,\pr,\ka)$, so $\Hp_\tau$
is generated by $\Lp_\tau$ by \eq{as7eq9}. The isomorphism
$\Phi^{\rm pa}_\tau$ follows using Proposition~\ref{as3prop2}.
\end{proof}

\subsection{Functions $\ep^\al(\tau)$ and the Lie algebra $\Lt_\tau$}
\label{as73}

We would like to prove an analogue of Proposition \ref{as7prop3} for
the Lie algebra $\Lt_\tau$. The methods of
\S\ref{as71}--\S\ref{as72} do not really help, as the restriction to
total orders $(I,\pr)$ in the spanning set $\dss(I,\pr,\ka,\tau)$
does not translate to nice restrictions in the other spanning sets
such as $\dsib(*,\tau)$. Instead we introduce alternative generators
$\ep^\al(\tau)$, $\al\in C(\A)$, for the algebra $\Ht_\tau$. These
will be important in the author's paper \cite{Joyc6} on holomorphic
generating functions for invariants counting $\tau$-semistable
objects.

\begin{dfn} Let Assumption \ref{as3ass} hold, and $(\tau,T,\le)$
be a permissible weak stability condition on $\A$. For $\al\in
C(\A)$, define $\ep^\al(\tau)$ in $\CF(\fObj_\A)$ by
\e
\ep^\al(\tau)=
\!\!\!\!\!\!\!\!\!\!\!\!\!\!\!\!
\sum_{\substack{\text{$\A$-data $(\{1,\ldots,n\},\le,\ka):$}\\
\text{$\ka(\{1,\ldots,n\})=\al$, $\tau\ci\ka\equiv\tau(\al)$}}}
\!\!\!\!\!\!\!
\frac{(-1)^{n-1}}{n}\,\,\dss^{\ka(1)}(\tau)*\dss^{\ka(2)}(\tau)*
\cdots*\dss^{\ka(n)}(\tau).
\label{as7eq17}
\e
If $n,\ka$ give a nonzero term in \eq{as7eq17} and $1\!\le\!i\!<\!n$
then $\be\!=\!\ka(\{1,\ldots,i\})$, $\ga\!=\!
\ka(\{i\!+\!1,\ldots,n\})$ lie in $C(\A)$ with
$\al\!=\!\be\!+\!\ga$, $\tau(\al)\!=\!\tau(\be)\!=\!\tau(\ga)$ and
$\Oss^\be(\tau) \!\ne\!\emptyset\!\ne\!\Oss^\ga(\tau)$. There are
only finitely many such $\be,\ga$ by Proposition \ref{as4prop1}, and
so only finitely many nonzero terms in \eq{as7eq17}. Thus
$\ep^\al(\tau)$ is well-defined.
\label{as7def2}
\end{dfn}

Here is the {\it inverse} of \eq{as7eq17}. The proof follows that of
Theorem \ref{as6thm2} closely, but using the associative
multiplication $*$ on $\CF(\fObj_\A)$ rather than $\square$, and
exchanging the r\^oles of $\exp$ and~$\log$.

\begin{thm} Let Assumption \ref{as3ass} hold, $(\tau,T,\le)$ be a
permissible weak stability condition on $\A$, and\/ $\be\in C(\A)$.
Then
\e
\dss^\be(\tau)= \!\!\!\!\!\!\!\!\!\!\!\!\!\!\!\!
\sum_{\substack{\text{$\A$-data $(\{1,\ldots,m\},\le,\la):$}\\
\text{$\la(\{1,\ldots,m\})=\be$, $\tau\ci\la\equiv\tau(\be)$}}}
\!\!\!\! \frac{1}{m!}\,\,\ep^{\la(1)}(\tau)*\ep^{\la(2)}(\tau)*
\cdots*\ep^{\la(m)}(\tau).
\label{as7eq18}
\e
There are only finitely many nonzero terms in~\eq{as7eq18}.
\label{as7thm1}
\end{thm}

Equations \eq{as7eq17}--\eq{as7eq18} show that the $\ep^\al(\tau)$
lie in the subalgebra of $\CF(\fObj_\A)$ generated by the
$\dss^\al(\tau)$ and vice versa, so they generate the same
subalgebra, which is $\Ht_\tau$ by \eq{as7eq2}. Therefore
\e
\Ht_\tau=\bigl\langle\de_{[0]},
\ep^{\al_1}(\tau)*\cdots*\ep^{\al_n}(\tau):
\al_1,\ldots,\al_n\in C(\A)\bigr\rangle{}_\Q.
\label{as7eq19}
\e
Here is an important property of the $\ep^\al(\tau)$, which the
coefficient $(-1)^{n-1}/n$ in \eq{as7eq17} was chosen to achieve.

\begin{thm} In Definition \ref{as7def2} we
have~$\ep^\al(\tau)\in\CFi(\fObj_\A)$.
\label{as7thm2}
\end{thm}

\begin{proof} Let $\al\in C(\A)$, $X\in\A$ with $[X]=\al$, and
$(\{1,\ldots,n\},\le,\ka)$ be $\A$-data with
$\ka(\{1,\ldots,n\})=\al$. By Definition \ref{as3def6} we have a
Cartesian square
\begin{equation*}
\xymatrix@C=90pt@R=10pt{ \fM(X,\{1,\ldots,n\},\le,\ka)_\A
\ar[r]_(0.55){\bs\si(\{1,\ldots,n\})}
\ar@<-3ex>[d]^{\,\Pi_X} & \Spec\K \ar[d]_{X\,} \\
\fM(\{1,\ldots,n\},\le,\ka)_\A \ar[r]^(0.55){\bs\si(\{1,\ldots,n\})}
&\fObj_\A. }
\end{equation*}
Applying \eq{as2eq5} to this and using \eq{as7eq4} and
$\CF(\Spec\K)=\Q$ we have
\begin{align}
(\dss^{\ka(1)}(\tau)\!*\!\cdots\!*\!\dss^{\ka(n)}(\tau))&([X])\!=\!
X^*\!\ci\!\CF^\stk(\bs\si(\{1,\ldots,n\}))\dss(\{1,\ldots,n\},\le,
\ka,\tau)
\nonumber\\
&=\!\CF^\stk(\bs\si(\{1,\ldots,n\}))\!\ci\!\Pi_X^*
(\dss(\{1,\ldots,n\},\le,\ka,\tau))
\nonumber\\
&=\CF^\stk(\bs\si(\{1,\ldots,n\}))\dss(X,\{1,\ldots,n\},\le,\ka,\tau)
\nonumber\\
&=\chi^\na\bigl(\Mss(X,\{1,\ldots,n\},\le,\ka,\tau)_\A\bigr).
\label{as7eq20}
\end{align}

To prove Theorem \ref{as7thm2} it is sufficient to show that if
$X=Y\op Z$ with $0\not\cong Y,Z$ then $\ep^\al(\tau)([X])=0$. By
\eq{as7eq17} and \eq{as7eq20} this is equivalent to
\e
\sum_{\substack{\text{$\A$-data $(\{1,\ldots,n\},\le,\ka):$}\\
\text{$\ka(\{1,\ldots,n\})\!=\!\al$,
$\tau\!\ci\!\ka\!\equiv\!\tau(\al)$}}} \!\!\!\!\!\!\!\!
\frac{(-1)^{n-1}}{n} \chi^\na\bigl(\Mss(Y\op
Z,\{1,\ldots,n\},\le,\ka,\tau)_\A\bigr)=0.
\label{as7eq21}
\e
Now $\Aut(Y\op Z)$ acts naturally on $\Mss(Y\op Z,\{1,\ldots,n\},
\le,\ka,\tau)_\A$. Define $G$ to be the subgroup $\{\id_Y+\ga\id_Z:0
\ne\ga\in\K\}$ of $\Aut(Y\op Z)$, so that $G\cong\K^\t$. Then each
orbit of $G$ on $\Mss(Y\op Z,\{1,\ldots,n\},\le,\ka,\tau)_\A$ is
either a single point or free. Since $\chi(\K^\t)=0$, by properties
of the Euler characteristic we have
\e
\begin{split}
&\chi^\na\bigl(\Mss(Y\op Z,\{1,\ldots,n\},\le,\ka,\tau)_\A\bigr)=\\
&\chi^\na\bigl(\Mss(Y\op Z,\{1,\ldots,n\},\le,\ka,\tau)_\A^G\bigr),
\end{split}
\label{as7eq22}
\e
where $(\cdots)^G$ is the fixed points of $G$, as the free orbits
contribute zero.

By \cite[Cor.~4.4]{Joyc3} there is a 1-1 correspondence between
$[(\si,\io,\pi)]\in\Mss(Y\!\op\!Z,\{1,\ldots,n\},\le,\ka,\tau)_\A$
and filtrations $0\!=\!A_0\!\subset\!\cdots\!\subset\!A_n\!=\!
Y\!\op\!Z$ with $S_i=A_i/A_{i-1}$ $\tau$-semistable with
$\tau([S_i])\!=\!\ka(i)$. The condition for $[(\si,\io,\pi)]$ to be
$G$-invariant turns out to be $A_i=B_i\op C_i$ for all $i$ as
subobjects of $Y\!\op\!Z$, where $B_i\!=\!A_i\!\cap\!Y$ and
$C_i\!=\!A_i\!\cap\!Z$. Then $0\!=\!B_0\!\subset\!\cdots\!
\subset\!B_n\!=\!Y$ and~$0\!=\!C_0\!\subset\!\cdots\!\subset\!
C_n\!=\!Z$.

Let $0=B_0'\subset\cdots\subset B_l'=Y$ and $0=C_0'\subset
\cdots\subset C_m'=Z$ be the filtrations obtained by omitting
repetitions, that is, omit $B_i$ if $B_i=B_{i-1}$ and so on. There
are unique maps $\phi:\{0,\ldots,n\}\ra\{0,\ldots,l\}$ and
$\psi:\{0,\ldots,n\}\ra\{0,\ldots,m\}$ with $B_i=B_{\phi(i)}'$ and
$C_i=C_{\psi(i)}'$ for all $i$. They are surjective, with $i\le j$
implies $\phi(i)\le\phi(j)$ and $\psi(i)\le\psi(j)$. Also, the
condition that $A_i\ne A_{i-1}$ implies that either
$\phi(i-1)\ne\phi(i)$ or $\psi(i-1)\ne\psi(i)$ for
all~$i=1,\ldots,n$.

Conversely, if we fix $l,m>0$ and filtrations $0=B_0'\subset
\cdots\subset B_l'=Y$ and $0=C_0'\subset\cdots\subset C_m'=Z$ such
that $T_i'=B_i'/B_{i-1}'$ and $U_i'=C_i'/C_{i-1}'$ are
$\tau$-semistable with $\tau([T_i'])=\tau([U_i'])=\tau(\al)$, the
possible $n,\ka$ and $0=A_0\subset\cdots\subset A_n=X$ coming from
$[(\si,\io,\pi)]\in\Mss(Y\op Z,\{1,\ldots,n\},\le,\ka,\tau)_\A^G$
and yielding these $k,l,B_i',C_i'$ from the construction above are
classified by such $\phi,\psi$. Therefore the contribution to
\eq{as7eq21} from such $[(\si,\io,\pi)]$ is
\e
\sum_{n=\max(l,m)}^{l+m}
\sum_{\substack{
\text{surjective $\phi:\{0,\ldots,n\}\ra\{0,\ldots,l\}$}\\
\text{and $\psi:\{0,\ldots,n\}\ra\{0,\ldots,m\}$:}\\
\text{$i\!\le\!j$ implies $\phi(i)\!\le\!\phi(j)$ and
$\psi(i)\!\le\!\psi(j)$,}\\
\text{$\phi(i\!-\!1)\!\ne\!\phi(i)$ or $\psi(i\!-\!1)
\!\ne\!\psi(i)$ for $1\le i\le n$}}} \frac{(-1)^{n-1}}{n}\,.
\label{as7eq23}
\e
We shall show \eq{as7eq23} is zero. Integrating this over all
$l,m,B_i',C_i'$ and using \eq{as7eq22} and properties of Euler
characteristics proves \eq{as7eq21}, and Theorem~\ref{as7thm2}.

For $n,\phi,\psi$ as in \eq{as7eq23}, define $E=\bigl\{i\in
\{1,\ldots,n\}:\phi(i-1)=\phi(i)$ and $\psi(i-1)\ne\psi(i) \bigr\}$
and $F=\bigl\{i\in\{1,\ldots,n\}:\phi(i-1)\ne\phi(i)$ and
$\psi(i-1)=\psi(i)\bigr\}$. Then $E,F$ are disjoint subsets of
$\{1,\ldots,n\}$ with $\md{E}=n-l$, $\md{F}=n-m$, and any such $E,F$
determine unique $\phi,\psi$. Thus for fixed $n$ the number of
$\phi,\psi$ in \eq{as7eq23} is $n!/(n-l)!(n-m)! (m+l-n)!$, and
\eq{as7eq23} reduces to
\e
\sum_{n=\max(l,m)}^{l+m}\frac{(-1)^{n-1}}{n}\,\cdot\,
\frac{n!}{(n-l)!(n-m)!(m+l-n)!}\,.
\label{as7eq24}
\e
Fixing $l>0$, multiplying \eq{as7eq24} by $t^m$ and summing over
$m=0,1,2,\ldots$ gives
\begin{gather*}
\frac{(-1)^{l-1}}{l}\sum_{n=l}^\iy\frac{(-1)^{n-l}(n-1)!\,
t^{n-l}}{(l-1)!(n-l)!}\sum_{m=n-l}^n\frac{l!\,t^{m+l-n}
}{(n-m)!(m+l-n)!}=\\
\frac{(-1)^{l\!-\!1}}{l}\sum_{a=0}^\iy\frac{(-1)^a(l\!-\!1\!+\!a)!
t^a}{(l\!-\!1)!\,a!}\sum_{b=0}^l\frac{l!\,t^b}{(l\!-\!b)!\,b!}\!=\!
\frac{(-1)^{l\!-\!1}}{l}(1\!+\!t)^{-l}(1\!+\!t)^l\!=\!
\frac{(-1)^{l\!-\!1}}{l},
\end{gather*}
using $a=n-l$, $b=m+l-n$ and the binomial theorem. Equating
coefficients of $t^m$, \eq{as7eq24} is zero when $m>0$, so
\eq{as7eq23} is zero. This completes the proof.
\end{proof}

Let $[X]\in\fObj_\A^\al(\K)$. For $\tau$-stable $X$, the only
nonzero term at $[X]$ in \eq{as7eq17} is $n=1$, $\ka(1)=\al$. If any
term in \eq{as7eq17} is nonzero at $[X]$ then $X$ has a filtration
$0=A_0\subset\cdots\subset A_n=X$ with $S_i=A_i/A_{i-1}$
$\tau$-semistable and $\tau([S_i])=\tau([X])$ for all $i$, so $X$ is
$\tau$-semistable. Hence by Theorem \ref{as7thm2} we have
\begin{equation*}
\ep^\al(\tau)\bigl([X])=\begin{cases} 1, &
\text{$X$ is $\tau$-stable,} \\
\text{in $\Q$,} & \text{$X$ is strictly $\tau$-semistable
and indecomposable,} \\
0, & \text{$X$ is $\tau$-unstable or decomposable.} \end{cases}
\end{equation*}
Thus $\ep^\al(\tau)$ interpolates between $\dsi^\al(\tau)$
and~$\dst^\al(\tau)$.

We can now prove an analogue of Proposition \ref{as7prop3}
for~$\Lt_\tau,\Ht_\tau$.

\begin{cor} Let Assumption \ref{as3ass} hold and $(\tau,T,\le)$
be a permissible weak stability condition on $\A$, and use the
notation of\/ \S\ref{as33} and\/ \S\ref{as71}. Then $\Lt_\tau$ is
the Lie subalgebra of\/ $\CFi(\fObj_\A)$ generated by the
$\ep^\al(\tau)$ for $\al\in C(\A)$. There is a natural\/
$\Q$-algebra isomorphism $\Phi^{\rm to}_\tau:U(\Lt_\tau)\ra
\Ht_\tau,$ where $U(\Lt_\tau)$ is the universal enveloping algebra
of\/~$\Lt_\tau$.
\label{as7cor2}
\end{cor}

\begin{proof} Write $\L'$ for the Lie subalgebra of
$\CFi(\fObj_\A)$ generated by the $\ep^\al(\tau)$ for all $\al\in
C(\A)$; this makes sense by Theorem \ref{as7thm2}. By \eq{as7eq19}
$\Ht_\tau$ is generated by the $\ep^\al(\tau)$, and so by $\L'$. But
$\L'\subseteq\Lt_\tau=\Ht_\tau\cap\CFi(\fObj_\A)$, so $\Ht_\tau$ is
also generated by $\Lt_\tau$. Thus Proposition \ref{as3prop2} gives
$\Q$-algebra isomorphisms $\Phi':U(\L')\ra\Ht_\tau$ and $\Phi^{\rm
to}_\tau :U(\Lt_\tau)\ra\Ht_\tau$. As $\L'\subseteq\Lt_\tau$ we have
$U(\L')\subseteq U(\Lt_\tau)$, with $\Phi^{\rm
to}_\tau\vert_{U(\L')}=\Phi'$. Since $\Phi',\Phi^{\rm to}_\tau$ are
isomorphisms this forces $\L'=\Lt_\tau$, so $\Lt_\tau$ is the Lie
subalgebra of $\CFi(\fObj_\A)$ generated by the~$\ep^\al(\tau)$.
\end{proof}

\subsection{The $\dss(*,\tau),\ldots$ have no universal linear
relations}
\label{as74}

The identities of \S\ref{as5}--\S\ref{as6} given in \eq{as7eq7}, and
their projections to $\CF(\fObj_\A)$ as in \eq{as7eq8}, are {\it
universal linear relations} between the $\dss,\dsi,\dst,\dssb,\dsib,
\dstb(*,\tau)$. By this we mean that they hold for all choices of
$\A,\fF_\A,K(\A),(\tau,T,\le)$ and auxiliary $\A$-data
$(I,\pr,\ka),\ldots$. Note also that each of these relations
expresses one of the families $\dss,\ldots,\dstb(*,\tau)$ in terms
of another; they can be thought of as basis change formulae between
six different bases of some universal algebra.

We claim that, in contrast, there are {\it no nontrivial universal
linear relations} involving just one of the families $\dss,\ldots,
\dstb(*,\tau)$. That is, the $\dss(I,\pr,\ka,\tau)$ over all
isomorphism classes of $\A$-data $(I,\pr,\ka)$ should have a kind of
{\it universal linear independence}: there are no systematic
relations on them that hold for all $\A,\fF_\A,K(\A),(\tau,T,\le)$,
only particular relations in each example. Before proving our
general result Theorem \ref{as7thm3}, we study an example and prove
linear independence of some collections of functions
in~$\CF(\fObj_\A)$.

\begin{ex} Fix a nonempty finite set $I$. Define a {\it quiver}
$Q=(Q_0,Q_1,b,e)$ to have vertices $I$ and an edge
$\smash{\mathop{\bu}\limits^{\sst i}\ra \mathop{\bu}\limits^{\sst
j}}$ for all $i,j\in I$, including $i=j$. That is, take $Q_0=I$,
$Q_1=I\t I$, $b:(i,j)\mapsto i$ and $e:(i,j)\mapsto j$. Set $\K=\C$
and consider the abelian category $\A=\nilCQ$ of nilpotent
$\C$-representations of $Q$, with data $K(\A),\fF_\A$ satisfying
Assumption \ref{as3ass} as in \cite[Ex.~10.6]{Joyc3}. Then
$K(\A)=\Z^I$, with elements of $K(\A)$ written as maps $\al:I\ra\Z$,
and $C(\A)$ is $\N^I\sm\{0\}$. For $i\in I$ define $e_i\in C(\A)$ by
$e_i(j)=1$ if $j=i$ and $e_i(j)=0$ otherwise. Then $\sum_{i\in
I}e_i=1$, where $1\in C(\A)$ maps $i\mapsto 1$ for all $i\in I$. Let
$(\tau,T,\le)$ be any (weak) stability condition on $\A$, such as
one defined using a {\it slope function\/} in Example \ref{as4ex1}.
Then $(\tau,T,\le)$ is {\it permissible\/} by
Corollary~\ref{as4cor}.

For $i\in I$ define ${\mathbf V}^i=(V^i,\rho^i)$ in $\A$ by
$V^i_i=\C$, $V^i_j=0$ for $i\ne j\in I$ and $\rho(e)=0$ for all
edges $e$ in $Q$. Then $[{\mathbf V}^i]=e_i$ in $C(\A)$, and
$\fObj_\A^i(\C)=\{[{\mathbf V}^i]\}$. Also ${\mathbf V}^i$ is
simple, so it is automatically $\tau$-stable, and we see that
\e
\dss^{e_i}(\tau)=\dsi^{e_i}(\tau)=\dst^{e_i}(\tau)=\ep^{e_i}(\tau)
=\de_{[{\mathbf V}^i]}.
\label{as7eq25}
\e
Define $\al=1\in C(\A)$ and $\ka:I\ra C(\A)$ by~$\ka(i)=e_i$.
\label{as7ex}
\end{ex}

\begin{prop} In the situation of Example \ref{as7ex} we have:
\begin{itemize}
\setlength{\itemsep}{0pt}
\setlength{\parsep}{0pt}
\item[{\rm(a)}] There exists no $\A$-data $(J,\ls,\la)$ with\/
$\md{J}>\md{I}$ and\/ $\la(J)=\al$. If\/ $(J,\ls,\la)$ is $\A$-data
with\/ $\md{J}=\md{I}$ and\/ $\la(J)=\al$, then there is a unique
bijection $\imath:J\ra I$ with\/~$\la=\ka\ci\imath$.
\item[{\rm(b)}] The functions $\CF^\stk(\bs\si(I))\dss(I,\pr,\ka,\tau)$
for all partial orders $\pr$ on $I$ are linearly independent in
$\CF(\fObj_\A)$. The same applies with\/ $\dss(\cdots)$ replaced by
$\dsi,\dst,\dssb,\dsib(\cdots)$ or~$\dstb(\cdots)$.
\item[{\rm(c)}] The subalgebra of\/ $\CF(\fObj_\A)$ generated by
$\dss^{e_i}(\tau)$ for $i\in I$ is freely generated, that is, there
are no polynomial relations in $\CF(\fObj_\A)$ on the
$\dss^{e_i}(\tau)$ for $i\in I$. The same holds for the
$\dsi^{e_i}(\tau)$, $\dst^{e_i}(\tau)$, and\/~$\ep^{e_i}(\tau)$.
\end{itemize}
\label{as7prop4}
\end{prop}

\begin{proof} For (a), if $(J,\ls,\la)$ is $\A$-data with
$\la(J)=\al=1$ then there must exist $\md{J}$ elements $\la(j)$ of
$\N^I\sm\{0\}$ adding up to 1. This is clearly impossible if
$\md{J}>\md{I}$, and if $\md{J}=\md{I}$ the elements $\la(j)$ for
$j\in J$ must be the set of all $e_i$, so there is a unique
bijection $\imath:J\ra I$ with $\la(i)=e_i$ for all~$i\in I$.

For (b), let $\tl$ be a partial order on $I$, and define
$(V,\rho)\in\A$ by
\begin{equation*}
V_i=\C\quad\text{for all $i\in I$ and}\quad
\rho(\smash{\mathop{\bu}\limits^{\sst i}\ra
\mathop{\bu}\limits^{\sst j}})=\begin{cases} 0, & i\tl j, \\
1, & i\ntl j. \end{cases}
\end{equation*}
Now $\A=\nilCQ$ is an abelian category of finite length, so by the
{\it Jordan--H\"older Theorem} the object $(V,\rho)$ in $\A$ has a
composition series into {\it simple factors}, which are unique up to
order and isomorphism. By construction these simple factors are
exactly ${\mathbf V}^i$ for~$i\in I$.

As the simple factors of $(V,\rho)$ are {\it pairwise
nonisomorphic}, we can apply the work of \cite[\S 3--\S 4]{Joyc3}.
These construct a unique partial order $\ls$ on the set $I$ indexing
the simple factors of $(V,\rho)$, and a {\it best\/
$(I,\ls)$-configuration} $(\si,\io,\pi)$ with $\si(I)=(V,\rho)$ and
$\si(\{i\})\cong{\mathbf V}^i$ for $i\in I$, which is unique up to
canonical isomorphism. Furthermore the $(I,\le)$ s-sets $J$
correspond to subobjects $S^J$ of $(V,\rho)$ induced
by~$\io(J,I):\si(J)\ra(V,\rho)$.

Now it is not difficult to show that the subobjects of $(V,\rho)$
are given by vector subspaces $V^J\le V$ of the form
$V^J=\bigop_{j\in J}V_j$ for $J\subseteq I$ an $(I,\tl)$ s-set.
Hence $\ls=\tl$, and by \eq{as7eq25} we see that for partial orders
$\pr$ on $I$ we have
\begin{equation*}
\Mssb,\Msib,\Mstb\bigl((V,\rho),I,\pr,\ka,\tau\bigr)_\A=
\begin{cases}\{[\si,\io,\pi]\}, & \pr=\tl, \\
\emptyset, & \pr\ne\tl.\end{cases}
\end{equation*}
As $(V,\rho)$ determines $(\si,\io,\pi)$ up to canonical isomorphism
$\bs\si(I)_*:\Aut(\si,\io,\pi)\ra\Aut(V,\rho)$ is an isomorphism, so
this implies that
\e
\CF^\stk(\bs\si(I))\dssb,\dstb,\dsib(I,\pr,\ka,\tau)\bigl([(V,\rho)]\bigr)
=\begin{cases} 1, & \pr=\tl, \\ 0, & \pr\ne\tl.
\end{cases}
\label{as7eq26}
\e

Since we can find such $(V,\rho)$ for each partial order $\tl$ on
$I$, \eq{as7eq26} implies the $\CF^\stk(\bs\si(I))
\dssb(I,\pr,\ka,\tau)$ for all $\pr$ on $I$ are linearly
independent, and similarly for $\dsib,\dstb(\cdots)$. But applying
$\CF^\stk(\bs\si(I))$ to \eq{as5eq7} and \eq{as6eq3} show that the
$\CF^\stk(\bs\si(I))\dssb(I,\pr,\ka,\tau)$ and $\CF^\stk(\bs\si(I))
\dss(I,\pr,\ka,\tau)$ over all $\pr$ span the same subspace of
$\CF(\fObj_\A)$, with dimension the number of partial orders on $I$,
so the $\CF^\stk(\bs\si(I))\dss(I,\pr,\ka,\tau)$ over all $\pr$ must
also be linearly independent. The same holds for
$\dsi,\dst(\cdots)$, using \eq{as5eq8}--\eq{as5eq9} and
\eq{as6eq4}--\eq{as6eq5}. This proves~(b).

For (c), let $(i_1,\ldots,i_n)$ be an ordered sequence in $I$,
allowing repeated elements. Then using similar techniques we can
construct $(V',\rho')$ in $\A=\nilCQ$ and a
$(\{1,\ldots,n\},\le)$-configuration $(\si',\io',\pi')$, unique up
to canonical isomorphism, with $\si'(\{1,\ldots,n\})=(V',\rho')$ and
$\si'(\{a\})\cong{\mathbf V}^{i_a}$ for $a=1,\ldots,n$, such that
there exists no such configuration for any other sequence
$(j_1,\ldots,j_m)$ in $I$. It follows that $\de_{[{\mathbf
V}^{j_1}]}*\cdots*\de_{[{\mathbf V}^{j_m}]}([(V',\rho')])$ is 1 if
$(j_1,\ldots,j_m)=(i_1,\ldots,i_n)$ and 0 otherwise. So the
$\de_{[{\mathbf V}^{i_1}]}*\cdots*\de_{[{\mathbf V}^{i_n}]}$ for all
sequences $(i_1,\ldots,i_n)$ are linearly independent in
$\CF(\fObj_\A)$. Part (c) now follows from~\eq{as7eq25}.
\end{proof}

To prove no universal linear relations exist on the $\dss(*,\tau)$,
say, we need to explain just what we mean by a universal linear
relation, which is not very obvious. In our next result we adopt a
rather restrictive definition \eq{as7eq27}, which includes the
identities of \S\ref{as5}--\S\ref{as6} and is sufficient for the
applications below. But the author expects the same principle to
hold for other universal forms.

\begin{thm} There exist no universal linear relations of the form
\e
\begin{aligned}
\sum_{\substack{\text{iso. classes of\/ $\A$-data $(J,\ls,\la)$
and surjective $\psi\!:\!J\!\ra\!K$:}\\
\text{$i\ls j$ $\Rightarrow$ $\psi(i)\!\tl\!\psi(j)$,
$\la(\psi^{-1}(k))\!=\!\mu(k)$ for $k\!\in\!K$,
$\tau\!\ci\!\mu\!\ci\!\psi\!\equiv\!\tau\!\ci\!\la$}}
\!\!\!\!\!\!\!\!\!\!\!\!\!\!\!\!\!\!\!\!\!\! }
\!\!\!\!\!\!\!\!\!\!\!\!\!\!\!\!\!\!\!\!\!\!\!\!\!\!\!\!\!\!\!\!\!
\!\!\!\!\!\!\!\!\!\!\!
C_{J,\ls,K,\tl,\psi}\CF^\stk(\bs\si(J))\dss(J,\ls,\la,\tau)=0
\end{aligned}
\label{as7eq27}
\e
in $\CF(\fObj_\A)$, which hold for all choices of\/
$\A,\fF_\A,K(\A)$ satisfying Assumption \ref{as3ass}, permissible
stability conditions or weak stability conditions $(\tau,T,\le)$ on
$\A$, and\/ $\A$-data $(K,\tl,\mu)$, where
$C_{J,\ls,K,\tl,\psi}\in\Q$ depends only on $J,\ls,K,\tl,\psi$ up to
isomorphism and is nonzero for at least one choice of\/
$J,\ldots,\psi$. The same applies with\/ $\dss(\cdots)$ replaced by
$\dsi,\dst,\dssb,\dsib(\cdots)$ or~$\dstb(\cdots)$.
\label{as7thm3}
\end{thm}

\begin{proof} Suppose for a contradiction that some such universal
linear relation exists. Choose $I,\pr,K,\tl,\phi$ with $\md{I}$ {\it
minimal\/} such that $C_{I,\pr,K,\tl,\phi}\ne 0$. Apply Example
\ref{as7ex} with this $I$, to get $\A,\fF_\A,K(\A)$ and $\ka:I\ra
C(\A)$. Let $(\tau,T,\le)$ be the trivial stability condition
$T=\{0\}$, $\tau\equiv 0$. Define $\mu:K\ra C(\A)$ by
$\ka(\phi^{-1}(k))=\mu(k)$ for $k\in K$. Then
$\tau\ci\mu\ci\phi\equiv\tau\ci\ka$ by choice of~$T$.

Consider equation \eq{as7eq27} with this data. Suppose
$(J,\ls,\la),\psi$ gives a nonzero term. We cannot have
$\md{J}<\md{I}$, since then $C_{J,\ls,K,\tl,\psi}=0$ by choice of
$I$. We cannot have $\md{J}>\md{I}$, as
$\tau\ci\mu\ci\psi\equiv\tau\ci\la$ implies $\la(J)=\mu(K)=\al$,
contradicting Proposition \ref{as7prop4}(a). Thus $\md{J}=\md{I}$,
and Proposition \ref{as7prop4}(a) gives a bijection $\imath:J\ra I$
with $\la=\ka\ci\imath$. Since the $\ka(i)$ for $i\in I$ are
linearly independent in $C(\A)$, and $\tau\ci\mu\ci\psi\equiv\tau
\ci\la$ we see that $\psi=\phi\ci\imath$. Thus $(J,\ls,\la),\psi$
are isomorphic to $(I,\ps,\ka),\phi$ for $\ps=\imath_*(\ls)$, and
\eq{as7eq27} reduces to
\begin{equation*}
\ts\sum_{\text{partial orders $\ps$ on $I$}}
C_{I,\ps,K,\tl,\phi}\CF^\stk(\bs\si(I))\dss(I,\ps,\ka,\tau)=0.
\end{equation*}
But the $\CF^\stk(\bs\si(I))\dss(I,\ps,\ka,\tau)$ for all $\ps$ are
linearly independent by Proposition \ref{as7prop4}(b), and
$C_{I,\pr,K,\tl,\phi}\ne 0$, a contradiction. The proof for
$\dsi,\ldots,\dstb(\cdots)$ is the same.
\end{proof}

Here are some remarks on this:
\begin{itemize}
\setlength{\itemsep}{0pt}
\setlength{\parsep}{0pt}
\item This implies a second result on nonexistence of universal
linear relations in $\CF(\fM(K,\tl,\mu)_\A)$ with
$\CF^\stk(\bs\si(J))\dss(J,\ls,\la,\tau)$ in \eq{as7eq27} replaced
by $\CF^\stk(Q(J,\ls,K,\tl,\phi))\dss(J,\ls,\la,\tau)$, since
applying $\CF^\stk(\bs\si(K))$ to such a relation would yield one of
the form \eq{as7eq27}. The identities of \eq{as7eq7} are of this
form, though mixing different families~$\dss,\ldots,\dstb(\cdots)$.
\item This second result shows that the identities of \eq{as7eq7}
are {\it unique} as universal linear relations. So, for instance,
\eq{as5eq7} is the {\it only} universal way to write
$\dss(I,\pr,\ka,\tau)$ in terms of the $\dssb(*,\tau)$, at least in
the form \eq{as7eq27}, since if there was another way we could take
the difference with \eq{as5eq7} to get a universal relation on the
$\dssb(*,\tau)$, contradicting Theorem~\ref{as7thm3}.
\item We can use a similar method of proof with Proposition
\ref{as7prop4}(c) to show there are {\it no universal polynomial
relations} in the $\dss^\al(\tau)$ for $\al\in C(\A)$, and similarly
for the $\dsi^\al,\dst^\al(\tau)$, and $\ep^\al(\tau)$. Effectively
this shows that the universal model for $\Ht_\tau$ is the free
associative $\Q$-algebra generated by $\dss^\al(\tau)$ for $\al\in
C(\A)$, or equivalently by $\ep^\al(\tau)$ for~$\al\in C(\A)$.
\item The theorem is evidence that the configurations framework is a
good one, and in particular, that partial orders are a good choice
of combinatorial data to keep track of collections of objects and
morphisms. For we know by closure of $\Ht_\tau$ under $*$ and other
operations that there are not too few partial orders to do
everything we want, and the theorem tells us there is no redundant
information, and so not too many partial orders.
\end{itemize}

\section{Generalization to stack functions}
\label{as8}

Finally we discuss the best way to generalize the constructible
functions material of \S\ref{as5}--\S\ref{as7} to stack functions.
We would like to define stack function versions
$\bdss^\al,\bdsi^\al,\bdst^\al(\tau)$ of $\dss^\al,\dsi^\al,
\dst^\al(\tau)$, and $\bdss,\ldots,\bdstb(I,\pr,\ka,\tau)$ of
$\dss,\ldots,\dstb(I,\pr,\ka,\tau)$, that satisfy analogues of the
identities of \S\ref{as5}--\S\ref{as6} and the (Lie) algebra ideas
of \S\ref{as7}; also, we want the transformation laws between
stability conditions $(\tau,T,\le)$ and $(\ti\tau,\ti T,\le)$
studied in \cite{Joyc5} for these stack functions to be
well-behaved.

The most obvious way to define these stack functions is
$\bdss^\al(\tau)=\bde_{\Oss^\al(\tau)}$,
$\bdsi^\al(\tau)=\bde_{\Osi^\al(\tau)},\ldots,
\bdstb(I,\pr,\ka,\tau)=\bde_{\Mstb(I,\pr,\ka,\tau)}$,
following Definition \ref{as4def4}. However, investigation shows
that this is {\it not\/} a helpful definition: none of the
identities of \S\ref{as5}--\S\ref{as6} would then hold, and much of
the (Lie) algebra material of \S\ref{as7} would not generalize
either.

There are two main reasons for this. The first is that constructible
function pushforwards $\CF^\stk(\cdots)$ use Euler characteristics,
and many of the identities of \S\ref{as5}--\S\ref{as6} make
essential use of $\chi(\K^m)=1$, and so will not work for general
stack function pushforwards. We could get round this by using the
stack function spaces $\oSF(\fF,\Th,\Om)$ of \cite[\S 6]{Joyc2},
which also set~$[\K^m]=1$.

The second is the idea of {\it virtual rank\/} introduced in
\cite[\S 5]{Joyc2}, and the corresponding idea of {\it virtual
indecomposable} in \cite[\S 5]{Joyc4}. The point here is that
experience shows that the best analogue of constructible functions
$\CFi(\fObj_\A)$ supported on indecomposables is not stack functions
supported on indecomposables, but stack functions $\SFai(\fObj_\A)$
`supported on virtual indecomposables', which can have nontrivial
components over decomposable objects.

Unfortunately these notions of virtual rank and virtual
indecomposable are technical and difficult to explain, but here is
the rough idea. On the stack functions $\SF(\fObj_\A)$ (or
$\SFa(\fObj_\A),\ldots$) we define linear maps $\Pi_n^{\rm
vi}:\SF(\fObj_\A)\ra\SF(\fObj_\A)$ for $n=0,1,2,\ldots$, the
projections to stack functions of `virtual rank $n$'. These satisfy
$(\Pi_n^{\rm vi})^2=\Pi_n^{\rm vi}$ and $\Pi_m^{\rm vi}\Pi_n^{\rm
vi}=0$ for $m\ne n$. If $[(\fR,\rho)]\in \SF(\fObj_\A)$ and $\fR$ is
a $\K$-stack whose stabilizer groups are all {\it abelian\/}
algebraic $\K$-groups, then $\Pi_n^{\rm
vi}\bigl([(\fR,\rho)]\bigr)=[(\fR_n,\rho)]$, where $\fR_n$ is the
locally closed $\K$-substack of $\fR$ of points whose stabilizer
groups have rank exactly~$n$.

If $[(\fR,\rho)]\in\SF(\fObj_\A)$ and $\fR$ is a $\K$-stack whose
stabilizer groups are nonabelian, then $\Pi_n^{\rm
vi}\bigl([(\fR,\rho)]\bigr)$ replaces each point $x\in\fR(\K)$ with
stabilizer group $\Aut_\K(x)=G$ by a finite $\Q$-linear combination
of points with stabilizer groups $C_G(T)$, the centralizer of $T$ in
$G$, for certain subgroups $T$ of the maximal torus $T^G$ of $G$. It
is like regarding a nonabelian stabilizer group $G$ as a formal
$\Q$-linear combination of torus stabilizer groups $(\C^\t)^k$ for
$\rk\,Z(G)\le k\le\rk\,G$, where $Z(G)$ is the centre of $G$, and
then $\Pi_n^{\rm vi}$ selects the $(\C^\t)^n$ components.

An object $X\in\A$ is indecomposable if and only if $\Aut(X)$ has
rank 1. By analogy, a stack function $f\in\SF(\fObj_\A)$ is said to
be {\it supported on virtual indecomposables\/} if it has virtual
rank 1, that is, $\Pi_1^{\rm vi}(f)=f$. We write $\SFai(\fObj_\A)$
for the subspace of $f\in\SFa(\fObj_\A)$ supported on virtual
indecomposables. The importance of these ideas for us is that there
is a deep compatibility between the projections $\Pi_n^{\rm vi}$ and
multiplication $*$ in $\SF(\fObj_\A),\SFa(\fObj_\A),\ldots$,
explored in \cite[\S 5]{Joyc4}. This implies, for instance, that
$\SFai(\fObj_\A)$ is a Lie algebra, that is, it is closed under the
Lie bracket $[f,g]=f*g-g*f$. In contrast, the subspace of
$f\in\SFa(\fObj_\A)$ supported on (actual, non-virtual)
indecomposable objects is not closed under~$[\,,\,]$.

These ideas suggest that the best definition for $\bdsi^\al(\tau)$
is not $\bde_{\Osi^\al(\tau)}$, but rather a `characteristic
function' of `$\tau$-semistable virtual indecomposables', perhaps
$\Pi^\vi_1(\bde_{\Oss^\al(\tau)})$ in the notation of \cite[\S
5]{Joyc2}, as in Theorem \ref{as8thm3} below. Following similar
reasoning, one can argue there should be stack function ideas of
`virtual $\tau$-stables' and `virtual best configurations', which
can have nonzero components over strictly $\tau$-semistable objects
and non-best configurations. However, there does not seem to be a
stack function idea of `virtual $\tau$-semistable': the appropriate
notion is just $\tau$-semistable in the usual sense.

Therefore the approach we choose is to first set
$\bdss^\al(\tau)=\bde_{\Oss^\al(\tau)}$ and
$\bdss(I,\pr,\ka,\tau)\!=\!\bde_{\Mss(I,\pr,\ka,\tau)}$, and then
define $\bdsi^\al,\bdst^\al(\tau)$ and
$\bdsi,\ldots,\bdstb(I,\pr,\ka, \tau)$ uniquely such that the
analogues of the identities of \S\ref{as5}--\S\ref{as6} hold. Of
course, the meaning of $\bdsi^\al,\bdst^\al(\tau)$ and
$\bdsi,\ldots,\bdstb(I,\pr,\ka,\tau)$ is then unclear, and we
discuss this after Theorem \ref{as8thm3}. The justification for this
approach is that nearly all of the (Lie) algebra material of
\S\ref{as7} generalizes very neatly, as we shall see below, and it
fits nicely with the ideas on changing stability conditions
in~\cite{Joyc5}.

For simplicity we work throughout with the spaces $\SF(\fF)$, but
the material below works equally well in the spaces
$\oSF(\fF,\Up,\La), \oSF(\fF,\Up,\La^\ci)$ or $\oSF(\fF,\Th,\Om)$ of
\cite{Joyc2}, and much of it also in~$\uSF(\fF,\Up,\La)$.

\begin{dfn} Let Assumption \ref{as3ass} hold, $(\tau,T,\le)$ be
a permissible weak stability condition on $\A$, $\al\in C(\A)$, and
$(I,\tl,\ka)$ be $\A$-data, as in Definition \ref{as3def5}. Define
\e
\begin{gathered}
\text{$\bdss^\al(\tau)=\bde_{\Oss^\al(\tau)}\in\SFa(\fObj_\A)$ or
$\SF(\fObj_\A^\al)$}\\
\text{and $\bdss(I,\tl,\ka,\tau)=\bde_{\Mss(I,\tl,\ka,\tau)}
\in\SF(\fM(I,\tl,\ka)_\A)$.}
\end{gathered}
\label{as8eq1}
\e
Since $\Mss(I,\tl,\ka,\tau)=\bigl(\prod_{i\in
I}\bs\si(\{i\})\bigr)_*^{-1}\bigl(\prod_{i\in
I}\Oss^{\ka(i)}(\tau)\bigr)$ we see that
\e
\bdss(I,\tl,\ka,\tau)=\ts\bigl(\prod_{i\in I}\bs\si(\{i\})\bigr)^*
\bigl(\bigot_{i\in I}\bdss^{\ka(i)}(\tau)\bigr).
\label{as8eq2}
\e
By analogy with \eq{as6eq3}, for $\A$-data $(I,\tl,\ka)$ define
\e
\bdssb(I,\tl,\ka,\tau)=
\sum_{\substack{\text{p.o.s $\pr$ on $I$:}\\
\text{$\tl$ dominates $\pr$}}}\!\!\!\!\!\!
n(I,\pr,\tl)Q(I,\pr,\tl)_*\,\bdss(I,\pr,\ka,\tau).
\label{as8eq3}
\e

By analogy with \eq{as6eq6}, setting $i\bu j$ if and only if $i=j$,
for $\al\in C(\A)$ define
\e
\bdsi^\al(\tau)=
\sum_{n=1}^\iy\frac{(-1)^{n-1}}{n}\cdot \!\!
\sum_{\begin{subarray}{l}\text{$\ka:\{1,\ldots,n\}\ra C(\A):$}\\
\text{$\ka(\{1,\ldots,n\})\!=\!\al$, $\tau\!\ci\!\ka\!\equiv
\!\tau(\al)$}\end{subarray}
\!\!\!\!\!\!\!\!\!\!\!\!\!\!\!\!\!\!\!\!\!\!\!\!\!\!\!\!\!\!
\!\!\!\!\!\!\!\!\!\!\!\! }
\bs\si(\{1,\ldots,n\})_*\,\bdss(\{1,\ldots,n\},\bu,\ka,\tau).
\label{as8eq4}
\e
By analogy with \eq{as8eq2} and \eq{as6eq4}, define
\ea
\bdsi(I,\tl,\ka,\tau)&=\ts\bigl(\prod_{i\in I}\bs\si(\{i\})\bigr)^*
\bigl(\bigot_{i\in I}\bdsi^{\ka(i)}(\tau)\bigr),
\nonumber\\
\bdsib(I,\tl,\ka,\tau)&=
\sum_{\substack{\text{p.o.s $\pr$ on $I$:}\\
\text{$\tl$ dominates $\pr$}}}\!\!\!\!\!\!
n(I,\pr,\tl)Q(I,\pr,\tl)_*\,\bdsi(I,\pr,\ka,\tau).
\label{as8eq5}
\ea

Now let $(\tau,T,\le)$ be a stability condition (not just a weak
one). By analogy with the case $K=\{k\}$ in \eq{as6eq14}, using
\eq{as6eq13} to simplify the $N(J,\ls,K,\chi)$, define
\e
\begin{aligned}
\bdst^\al(\tau)=
\sum_{\substack{\text{iso. classes}\\ \text{of finite sets $I$}}}
\frac{1}{\md{I}!}\cdot
\sum_{\begin{subarray}{l}
\text{$\pr,\ka$: $(I,\pr,\ka)$ is $\A$-data,}\\
\text{$\ka(I)=\al$, $\tau\ci\ka\equiv\tau(\al)$} \end{subarray}
\!\!\!\!\!\!\!\!\!\!\!\!\!\!\!\!\!\!\!\!\!\!\!\!\!\!\!\!\!\!\!\!\!\!\!}
N(I,\pr)\bs\si(I)_*\,\bdssb(I,\pr,\ka,\tau).
\end{aligned}
\label{as8eq6}
\e
By analogy with \eq{as8eq2} and \eq{as6eq5}, define
\ea
\bdst(I,\tl,\ka,\tau)&=\ts\bigl(\prod_{i\in I}\bs\si(\{i\})\bigr)^*
\bigl(\bigot_{i\in I}\bdst^{\ka(i)}(\tau)\bigr),
\label{as8eq7}\\
\bdstb(I,\tl,\ka,\tau)&=
\sum_{\substack{\text{p.o.s $\pr$ on $I$:}\\
\text{$\tl$ dominates $\pr$}}}\!\!\!\!\!\!
n(I,\pr,\tl)Q(I,\pr,\tl)_*\,\bdst(I,\pr,\ka,\tau).
\label{as8eq8}
\ea
By analogy with \eq{as7eq17}, for $\al\in C(\A)$ define
\e
\bep^\al(\tau)= \!\!\!\!\!\!\!\!\!\!\!\!\!\!\!\!
\sum_{\substack{\text{$\A$-data $(\{1,\ldots,n\},\le,\ka):$}\\
\text{$\ka(\{1,\ldots,n\})=\al$, $\tau\ci\ka\equiv\tau(\al)$}}}
\!\!\!\!\!\!\!
\frac{(-1)^{n-1}}{n}\,\,\bdss^{\ka(1)}(\tau)*\bdss^{\ka(2)}(\tau)*
\cdots*\bdss^{\ka(n)}(\tau).
\label{as8eq9}
\e
By the proofs in \S\ref{as5}--\S\ref{as7} there are only finitely
many nonzero terms in each equation, so they are all well-defined.
It is easy to show $\bdss^\al(\tau),\bdsi^\al(\tau),
\bdst^\al(\tau),\bep^\al(\tau)$ are supported on $\Oss^\al(\tau)$,
and $\bdss,\ldots,\bdstb(I,\tl,\ka,\tau)$
on~$\Mss(I,\tl,\ka,\tau)_\A$.
\label{as8def1}
\end{dfn}

Here are the analogues of the remaining eight identities in
\eq{as7eq7}, that is, \eq{as5eq7}, \eq{as5eq8}, \eq{as5eq9},
\eq{as5eq14}, \eq{as5eq32}, \eq{as6eq10}, \eq{as6eq11} and
\eq{as6eq14} respectively.

\begin{thm} For all\/ $\A$-data $(K,\tl,\mu)$ and\/
$\al\in C(\A)$ we have
\begin{gather}
\sum_{\substack{\text{p.o.s $\pr$ on $K$:}\\
\text{$\tl$ dominates $\pr$}}}Q(K,\pr,\tl)_*\,
\bdssb(K,\pr,\mu,\tau)=\bdss(K,\tl,\mu,\tau),
\label{as8eq10}
\allowdisplaybreaks\\
\sum_{\substack{\text{p.o.s $\pr$ on $K$:}\\
\text{$\tl$ dominates $\pr$}}}Q(K,\pr,\tl)_*\,
\bdsib(K,\pr,\mu,\tau)=\bdsi(K,\tl,\mu,\tau),
\label{as8eq11}
\allowdisplaybreaks\\
\sum_{\substack{\text{p.o.s $\pr$ on $K$:}\\
\text{$\tl$ dominates $\pr$}}}Q(K,\pr,\tl)_*\,
\bdstb(K,\pr,\mu,\tau)=\bdst(K,\tl,\mu,\tau),
\label{as8eq12}
\allowdisplaybreaks\\
\begin{aligned}
\sum_{\substack{\text{iso. classes}\\ \text{of finite }\\
\text{sets $I$}}} \frac{1}{\md{I}!}\cdot
\sum_{\substack{
\text{$\ka\!:\!I\!\ra\!C(\A)$, surjective $\phi\!:\!I\!\ra\!K$:}\\
\text{$\ka(\phi^{-1}(k))=\mu(k)$ for $k\in K$,}\\
\text{$\tau\ci\mu\ci\phi\equiv\tau\ci\ka:I\ra T$.}\\
\text{Define $\pr$ on $I$ by $i\pr j$ if $i=j$}\\
\text{or $\phi(i)\ne\phi(j)$ and $\phi(i)\tl\phi(j)$}}} \!\!\!
\begin{aligned}[t]
Q(I,\pr,K,\tl,\phi)_*&\\
\bdsi(I,\pr,\ka,\tau)&=\\
\bdss(K,\tl,\mu&,\tau),
\end{aligned}
\end{aligned}
\label{as8eq13}
\allowdisplaybreaks\\
\begin{aligned}
\sum_{\substack{\text{iso. classes}\\ \text{of finite }\\
\text{sets $I$}}} \frac{1}{\md{I}!}\cdot
\sum_{\substack{
\text{$\pr,\ka,\phi$: $(I,\pr,\ka)$ is $\A$-data,}\\
\text{$\phi:I\ra K$ is surjective,}\\
\text{$i\pr j$ implies $\phi(i)\tl\phi(j)$,}\\
\text{$\ka(\phi^{-1}(k))=\mu(k)$ for $k\in K$,}\\
\text{$\tau\ci\mu\ci\phi\equiv\tau\ci\ka:I\ra T$}}}
\begin{aligned}[t]
Q(I,\pr,K,\tl,\phi)_*&\\
\bdstb(I,\pr,\ka,\tau)&=\\
\bdss(K,\tl,\mu&,\tau),
\end{aligned}
\end{aligned}
\label{as8eq14}
\allowdisplaybreaks\\
\begin{aligned}
\sum_{\substack{\text{iso.}\\ \text{classes}\\ \text{of finite }\\
\text{sets $I$}}}\!\! \frac{(-1)^{\md{I}-\md{K}}}{\md{I}!}\!\cdot
\!\!\!\!
\sum_{\substack{
\text{$\ka\!:\!I\!\ra\!C(\A)$, surjective $\phi\!:\!I\!\ra\!K$:}\\
\text{$\ka(\phi^{-1}(k))=\mu(k)$ for $k\in K$,}\\
\text{$\tau\ci\mu\ci\phi\equiv\tau\ci\ka:I\ra T$.}\\
\text{Define $\pr$ on $I$ by $i\pr j$ if $i=j$}\\
\text{or $\phi(i)\ne\phi(j)$ and $\phi(i)\tl\phi(j)$}}} \!\!\!\!\!
\begin{aligned}[t]
\ts\prod_{k\in K}\bigl(\md{\phi^{-1}(\{k\})}-1\bigr)!&\cdot\\
Q(I,\pr,K,\tl,\phi)_*&\\
\bdss(I,\pr,\ka,\tau)&=\\
\bdsi(K,\tl,\mu&,\tau),
\end{aligned}
\end{aligned}
\label{as8eq15}
\allowdisplaybreaks\\
\begin{aligned}
\sum_{\substack{\text{iso. classes}\\ \text{of finite }\\
\text{sets $I$}}} \frac{1}{\md{I}!}\cdot
\sum_{\substack{
\text{$\pr,\ka,\phi$: $(I,\pr,\ka)$ is $\A$-data,}\\
\text{$(I,\pr,K,\phi)$ is allowable,}\\
\text{$\tl=\cP(I,\pr,K,\phi)$,}\\
\text{$\ka(\phi^{-1}(k))=\mu(k)$ for $k\in K$,}\\
\text{$\tau\ci\mu\ci\phi\equiv\tau\ci\ka:I\ra T$}}}
\begin{aligned}[t]
Q(I,\pr,K,\tl,\phi)_*&\\
\bdstb(I,\pr,\ka,\tau)&=\\
\bdssb(K,\tl,\mu&,\tau),
\end{aligned}
\end{aligned}
\label{as8eq16}
\allowdisplaybreaks\\
\begin{aligned}
\sum_{\substack{\text{iso. classes}\\ \text{of finite }\\
\text{sets $I$}}} \frac{1}{\md{I}!}\cdot \!\!\!\!\!
\sum_{\substack{
\text{$\pr,\ka,\phi$: $(I,\pr,\ka)$ is $\A$-data,}\\
\text{$(I,\pr,K,\phi)$ is allowable,}\\
\text{$\tl=\cP(I,\pr,K,\phi)$,}\\
\text{$\ka(\phi^{-1}(k))=\mu(k)$ for $k\in K$,}\\
\text{$\tau\ci\mu\ci\phi\equiv\tau\ci\ka:I\ra T$}}}
\!\!\!\!\!\!\!\!\!\!\!\!\!\!\!\!\!\!\!\!\!\!\!\!
\begin{aligned}[t]
N(I,\pr,K,\phi)Q(I,\pr,K,\tl,\phi)_*&\\
\bdssb(I,\pr,\ka,\tau)&=\\
\bdstb(K,\tl,\mu&,\tau),
\end{aligned}
\end{aligned}
\label{as8eq17}
\allowdisplaybreaks\\
\bdss^\al(\tau)= \!\!\!\!\!\!\!\!\!\!\!\!\!\!\!\!
\sum_{\substack{\text{$\A$-data $(\{1,\ldots,n\},\le,\ka):$}\\
\text{$\ka(\{1,\ldots,n\})=\al$, $\tau\ci\ka\equiv\tau(\al)$}}}
\!\!\!\! \frac{1}{n!}\,\,\bep^{\ka(1)}(\tau)*\bep^{\ka(2)}(\tau)*
\cdots*\bep^{\ka(n)}(\tau),
\label{as8eq18}
\end{gather}
supposing $(\tau,T,\le)$ is a stability condition in \eq{as8eq12},
\eq{as8eq14}, \eq{as8eq16} and\/ \eq{as8eq17}. There are only
finitely many nonzero terms in each equation.
\label{as8thm1}
\end{thm}

\begin{proof} The proofs in \S\ref{as5}--\S\ref{as6} imply there are
only finitely many nonzero terms in each equation. Equations
\eq{as8eq10}--\eq{as8eq12} are the inverses of \eq{as8eq3},
\eq{as8eq5}, \eq{as8eq8} respectively, and follow from them by the
reverse of the argument in \S\ref{as61}. The argument used to prove
\eq{as6eq10} from \eq{as6eq6} proves \eq{as8eq15} from \eq{as8eq4},
using \eq{as8eq2} along the way. Equation \eq{as8eq13} then follows
from \eq{as8eq15} as it is its combinatorial inverse, reversing the
argument in \S\ref{as62} that \eq{as6eq10} is the inverse
of~\eq{as5eq14}.

Combining \eq{as8eq3}, \eq{as8eq12} and \eq{as8eq17} gives an
identity writing $\bdst(K,\tl,\mu,\tau)$ as a linear combination of
$Q(I,\pr,K,\tl,\phi)_*\,\bdss(I,\pr,\ka,\tau)$. When $K=\{k\}$ this
is equivalent to the combination of \eq{as8eq3} and \eq{as8eq6}, and
so holds. The general case of the identity follows from the case
$K=\{k\}$ by \eq{as8eq7}, using \eq{as8eq2} along the way. We can
then recover \eq{as8eq17} from this identity as we already know
\eq{as8eq3}, \eq{as8eq12} and their inverses \eq{as8eq10},
\eq{as8eq8}. Equation \eq{as8eq16} follows from \eq{as8eq17} as it
is its combinatorial inverse, reversing the argument in \S\ref{as64}
that \eq{as6eq14} is the inverse of \eq{as6eq11}. We obtain
\eq{as8eq14} by substituting \eq{as8eq16} into \eq{as8eq10}.
Finally, \eq{as8eq18} is proved from \eq{as8eq9} in the same way as
\eq{as7eq18} from~\eq{as7eq17}.
\end{proof}

\begin{cor} If\/ $\cha\K=0$, $\pi^\stk_{\fObj_\A}$ takes
$\bdss^\al,\bdsi^\al,\bdst^\al,\bep^\al(\tau)$ to
$\dss^\al,\dsi^\al,\dst^\al,\ep^\al(\tau)$, and\/
$\pi^\stk_{\fM(I,\tl,\ka)_\A}$ takes $\bdss,\ldots,\bdstb
(I,\tl,\ka,\tau)$ to $\dss,\ldots,\dstb(I,\tl, \ka,\tau)$, and\/
$\pi^\stk_{\fObj_\A}$ takes $\bs\si(I)_*\,\bdss,\ldots,\bdstb
(I,\tl,\ka,\tau)$ to $\CF^\stk(\bs\si(I))\dss,\ldots,
\dstb(I,\tl,\ka,\tau),$ supposing $(\tau,T,\le)$ is a stability
condition for the $\bdst^\al(\tau)$ and\/ $\bdst,\bdstb(\cdots)$
cases.
\label{as8cor1}
\end{cor}

\begin{proof} By definition $\bdss^\al(\tau)\!=\!\io_{\fObj_\A}(
\dss^\al(\tau))$, so $\pi^\stk_{\fObj_\A}(\bdss^\al(\tau))=
\dss^\al(\tau)$ as $\pi^\stk_{\fObj_\A}\ci\io_{\fObj_\A}$ is the
identity. Similarly $\pi^\stk_{\fM(I,\tl,\ka)_\A}
(\bdss(I,\tl,\ka,\tau))=\dss(I,\tl,\ka,\tau)$. Now the identities of
Definition \ref{as8def1} and Theorem \ref{as8thm1} are all analogues
of identities on $\dss^\al,\ldots,\ep^\al(\tau)$ or
$\dss,\ldots,\dstb(I,\tl,\ka,\tau)$ in \S\ref{as5}--\S\ref{as7}. So
applying $\pi^\stk_{\fObj_\A},\pi^\stk_{\fM(I,\tl,\ka)_\A}$ or
$\pi^\stk_{\fM(K,\tl,\mu)_\A}$ to these identities and using Theorem
\ref{as2thm3}(b), we see that the identities of
\S\ref{as5}--\S\ref{as7} hold with $\pi^\stk_{\fObj_\A}
(\bdss^\al(\tau))$ in place of $\dss^\al(\tau)$, and so on.

But from \eq{as7eq7} we see the $\dss(*,\tau)$ determine the
$\dsi,\dst,\dssb,\dsib,\dstb(*,\tau)$, so as
$\pi^\stk_{\fM(I,\tl,\ka)_\A}(\bdss(I,\tl,\ka,\tau))
=\dss(I,\tl,\ka,\tau)$ we see that $\pi^\stk_{\fM(I,\tl, \ka)_\A}$
takes $\bdsi,\ldots,\bdstb(I,\tl,\ka,\tau)$ to
$\dsi,\ldots,\dstb(I,\tl,\ka,\tau)$. The claim for
$\bs\si(I)_*\,\bdss,\ldots,\bdstb\ab(I,\ab\tl,\ab\ka,\tau)$ now
follows from Theorem \ref{as2thm3}(b), and for $\bdsi^\al,
\bdst^\al,\bep^\al (\tau)$ from the corresponding identities.
\end{proof}

Here are stack function analogues of material in \S\ref{as71}
and~\S\ref{as73}.

\begin{dfn} Let Assumption \ref{as3ass} hold, and $(\tau,T,\le)$ be
a permissible weak stability condition on $\A$. Define $\Q$-vector
subspaces $\bHp_\tau,\bHt_\tau$ in $\SF(\fObj_\A)$ by
\begin{align*}
\bHp_\tau&=\bigl\langle\bs\si(I)_*\,\bdss(I,\pr,\ka,\tau):
\text{$(I,\pr,\ka)$ is $\A$-data}\bigr\rangle{}_\Q,\\
\bHt_\tau&=\bigl\langle\bde_{[0]},
\bdss^{\al_1}(\tau)*\cdots*\bdss^{\al_n}(\tau):\al_1,\ldots,\al_n\in
C(\A)\bigr\rangle{}_\Q.
\end{align*}
Here $\langle\cdots\rangle_\Q$ is the set of all finite $\Q$-linear
combinations of the elements `$\,\cdots$'. From \eq{as8eq2} we see
that $\bs\si(I)_*\,\bdss(I,\pr,\ka,\tau)=P_\sIp\bigl(\bdss^{\ka(i)}
(\tau):i\in I\bigr)$, giving
\begin{equation*}
\bHp_\tau=\bigl\langle P_\sIp\bigl(\bdss^{\ka(i)}(\tau):i\in
I\bigr): \text{$(I,\pr,\ka)$ is $\A$-data}\bigr\rangle{}_\Q.
\end{equation*}
It follows from \cite[Th.~5.4]{Joyc4} that $\bHp_\tau$ is closed
under the operations~$P_\sIp$.

If we were to work instead in $\oSF(\fObj_\A,\Up,\La)$, for
instance, it might be better to define
$\bHp_{\tau,\Up,\La},\bHt_{\tau,\Up,\La}$ to be the $\La$-submodules
with the above generators, and then
$\bHp_{\tau,\Up,\La},\bHt_{\tau,\Up,\La}$ will be $\La$-algebras
rather than $\Q$-algebras.
\label{as8def2}
\end{dfn}

In \cite{Joyc5} we will show that if $(\tau,T,\le)$ and
$(\ti\tau,\ti T,\le)$ are permissible weak stability conditions on
$\A$, then (under some finiteness conditions) we have
$\bHp_\tau=\bHp_{\ti\tau}$ and $\bHt_\tau=\bHt_{\ti\tau}$, so that
$\bHp_\tau,\bHt_\tau$ are {\it independent of the choice of weak
stability condition} $(\tau,T,\le)$. We generalize \eq{as7eq9}
and~\eq{as7eq19}.

\begin{thm} $\bHp_\tau,\bHt_\tau$ are subalgebras of\/
$\SFa(\fObj_\A)$ with\/ $\bHt_\tau\!\subseteq\!\bHp_\tau$, and
\ea
\begin{split}
\bHp_\tau&=\bigl\langle\bs\si(I)_*\,\bdssb(I,\pr,
\ka,\tau):\text{$(I,\pr,\ka)$ is $\A$-data}\bigr\rangle{}_\Q\\
&=\bigl\langle\bs\si(I)_*\,\bdsi(I,\pr,
\ka,\tau):\text{$(I,\pr,\ka)$ is $\A$-data}\bigr\rangle{}_\Q\\
&=\bigl\langle\bs\si(I)_*\,\bdsib(I,\pr,
\ka,\tau):\text{$(I,\pr,\ka)$ is $\A$-data}\bigr\rangle{}_\Q\\
&=\bigl\langle\bs\si(I)_*\,\bdst(I,\pr,
\ka,\tau):\text{$(I,\pr,\ka)$ is $\A$-data}\bigr\rangle{}_\Q\\
&=\bigl\langle\bs\si(I)_*\,\bdstb(I,\pr,
\ka,\tau):\text{$(I,\pr,\ka)$ is $\A$-data}\bigr\rangle{}_\Q,
\end{split}
\label{as8eq19}\\
\bHt_\tau&=\bigl\langle\bde_{[0]},\bep^{\al_1}(\tau)*\cdots*
\bep^{\al_n}(\tau):\al_1,\ldots,\al_n\in C(\A)\bigr\rangle{}_\Q,
\label{as8eq20}
\ea
supposing $(\tau,T,\le)$ is a stability condition in the last two
lines of\/ \eq{as8eq19}. When $\K$ has characteristic zero
$\pi^\stk_{\fObj_\A}:\bHp_\tau\ra\Hp_\tau$ and\/
$\pi^\stk_{\fObj_\A}:\bHt_\tau\ra\Ht_\tau$ are surjective
$\Q$-algebra morphisms.
\label{as8thm2}
\end{thm}

\begin{proof} Clearly $\bHt_\tau$ is the subalgebra of
$\SF(\fObj_\A)$ generated by the $\bdss^\al(\tau)$ for all $\al\in
C(\A)$. The analogue of \eq{as7eq4} implies that
$\bHt_\tau\!\subseteq\!\bHp_\tau$. We have
$\bdss^\al(\tau)=\io_{\fObj_\A}(\dss^\al(\tau))$, so
$\bdss^\al(\tau)\in\SFa(\fObj_\A)$ as $\io_{\fObj_\A}$ maps
$\CF(\fObj_\A)\ra\SFa(\fObj_\A)$ by \cite[Def.~5.5]{Joyc4}. Also
$\bs\si(I)_*\,\bdss(I,\pr,\ka,\tau)=P_\sIp\bigl(\bdss^{\ka(i)}
(\tau):i\in I\bigr)$ by \eq{as8eq2}, so
$\bs\si(I)_*\,\bdss(I,\pr,\ka,\tau)\in\SFa(\fObj_\A)$ by
\cite[Prop.~5.6]{Joyc4}, and $\bHt_\tau\!\subseteq\!\bHp_\tau
\subseteq\SFa(\fObj_\A)$. Since $\bHt_\tau$ is closed under the
$P_\sIp$ and $*=P_{\sst(\{1,2\},\le)}$, it is closed under $*$, and
is a subalgebra of $\SFa(\fObj_\A)$. Equation \eq{as8eq19} follows
from applying $\bs\si(I)_*$ or $\bs\si(K)_*$ to \eq{as8eq3},
\eq{as8eq5} and \eq{as8eq8}--\eq{as8eq17}, and \eq{as8eq20} from
\eq{as8eq9} and \eq{as8eq18}, as for \eq{as7eq19}. Finally,
Corollary \ref{as8cor1} implies $\pi^\stk_{\fObj_\A}$ induces
surjective maps $\bHp_\tau\ra\Hp_\tau$ and $\bHt_\tau\ra\Ht_\tau$,
which are $\Q$-algebra morphisms as \eq{as3eq9} is.
\end{proof}

The multiplication relations in $\bHp_\tau$ for the six spanning
sets $\bs\si(I)_* \,\bdss,\ldots,\ab\bdstb(I,\pr,\ka,\tau)$ are
given by the analogues of \eq{as7eq10} and \eq{as7eq11}. That is,
for $(I,\pr,\ka)$, $(J,\ls,\la)$, $(K,\tl,\mu)$ as defined before
\eq{as7eq10}, using \cite[Th.~5.4]{Joyc4} in place of
\cite[Th.~4.22]{Joyc4} shows the analogue of \eq{as7eq10} holds:
\begin{equation*}
\bigl(\bs\si(I)_*\,\bdss(I,\pr,\ka,\tau)\bigr)*
\bigl(\bs\si(J)_*\,\bdss(J,\ls,\la,\tau)\bigr)
=\bs\si(K)_*\,\bdss(K,\tl,\mu,\tau).
\end{equation*}
From this and identities \eq{as8eq3}, \eq{as8eq5}, \eq{as8eq8} and
\eq{as8eq10}--\eq{as8eq17} we can deduce multiplication relations
for the $\bs\si(I)_*\,\bdsi,\ldots,\bdstb(*)$. But as
\eq{as8eq3}--\eq{as8eq17} are analogues of constructible functions
identities these relations are exactly the analogues of the
constructible function relations~\eq{as7eq10}--\eq{as7eq11}.

Next we extend the Lie algebra material of \S\ref{as7}. The
following will be a key tool in proving elements of
$\bHp_\tau,\bHt_\tau$ lie in the Lie algebra~$\SFai(\fObj_\A)$.

\begin{thm} In Definition \ref{as8def1} we have~$\bdsi^\al(\tau)=
\Pi^\vi_1(\bdss^\al(\tau))$.
\label{as8thm3}
\end{thm}

\begin{proof} We shall combine \eq{as8eq1} with the definition of
$\Pi^\vi_1$ in \cite[\S 5.2]{Joyc2}, and show that the resulting
formula for $\Pi^\vi_1(\bdss^\al(\tau))$ agrees term-by-term with
the definition of $\bdsi^\al(\tau)$ in \eq{as8eq4}. Apply
\cite[Prop.~5.7]{Joyc4} with the constructible set
$S\subseteq\fObj_\A(\K)$ equal to $\Oss^\al(\tau)$. This gives a
finite decomposition $\Oss^\al(\tau)=\coprod_{l\in L}\fF_l(\K)$ and
1-isomorphisms $\fF_l\cong[U_l/A_l^\t]$, for $U_l$ a quasiprojective
$\K$-variety and $A_l$ a finite-dimensional $\K$-algebra, such that
if $u\in U_l(\K)$ projects to $[X]\in\fObj_\A(\K)$ then there exists
a subalgebra $B_u$ of $A_l$ with $\Stab_{A_l^\t}(u)=B_u^\t$ and an
isomorphism $B_u\cong\End(X)$ compatible
with~$\Stab_{A_l^\t}(u)\cong\Iso_\K([X])\cong\Aut(X)$.

Write $\rho_l:[U_l/A_l^\t]\ra\fObj_\A$ for the composition of
$\fF_l\cong[U_l/A_l^\t]$ and the inclusion $\fF_l\ra\fObj_\A$. Then
the definition \cite[Def.~3.2]{Joyc2} of $\bde_C$ implies that
\e
\bdss^\al(\tau)=\ts\sum_{l\in L}\bigl[\bigl([U_l/A_l^\t],
\rho_l\bigr)\bigr].
\label{as8eq21}
\e
There exists a subalgebra $C_l$ of $A_l$ isomorphic as an algebra to
$\K^{r_l}$, where $r_l=\rk\,A_l^\t$, and $C_l^\t\cong(\K^\t)^{r_l}$
is a maximal torus of $A_l^\t$. If $u\in U_l(\K)$ then
$\Stab_{A_l^\t}(u)\cap C_l^\t= D_u^\t$, where $D_u=B_u\cap C_l$ is a
subalgebra of $C_l$, for $B_u$ as above. It is now easy to see, in
the notation of \cite[\S 5.2]{Joyc2}, that
\e
\cP(U_l,C_l^\t),\cQ(A_l^\t,C_l^\t),\cR(U_l,A_l^\t,C_l^\t)\subseteq
\bigl\{D^\t:\text{$D\subseteq C_l$ a subalgebra}\bigr\}.
\label{as8eq22}
\e

It is a consequence of the proof in \cite[\S 5]{Joyc2} that the
definition of $\Pi^\vi_1$ is independent of choices, that in
defining $\Pi^\vi_1$ we can replace $\cP,\cQ,\cR(\cdots)$ by larger
sets of $\K$-subgroups of $C_l^\t$ closed under intersection. So, we
can define $\Pi^\vi_1(\bdss^\al(\tau))$ using the representation
\eq{as8eq21} and replacing the l.h.s.\ of \eq{as8eq22} by the
r.h.s.\ of \eq{as8eq22}. This involves a sum over $l\in L$ and
$P,Q,R$ in the r.h.s.\ of \eq{as8eq22} with $R\subseteq P\cap Q$ and
$\dim R=1$ of a term with coefficient~$M_{A_l^\t}^{U_l}(P,Q,R)$.

We can simplify this sum in four ways. Firstly, the only $R$ in the
r.h.s.\ of \eq{as8eq22} with $\dim R=1$ is $\{\la\id_{C_l}:\la\in
\K^\t\}$, so we fix $R$ to be this. Secondly, by
\cite[Lem.~5.9]{Joyc2} if $M_{A_l^\t}^{U_l}(P,Q,R)\ne 0$ then $P,Q$
are the smallest elements of their sets containing $P\cap Q$, so as
$P,Q$ take values in the same set we can restrict to $P=Q=D^\t$.
Thirdly, if $D\subseteq C_l\cong\K^{r_l}$ is a subalgebra with $\dim
D=n$ then $D\cong\K^n$ and explicit calculation with the definitions
of \cite[\S 5.2]{Joyc2} shows that
\begin{equation*}
M_{A_l^\t}^{U_l}(D^\t,D^\t,R)=
\left\vert\frac{N_{A_l^\t}(C_l^\t)}{C_{A_l^\t}(D^\t)\!\cap\!
N_{A_l^\t}(C_l^\t)}\right\vert^{-1}\!\cdot(-1)^n(n-1)!,
\end{equation*}
computing $M_{A_l^\t}^{U_l}(\cdots)$ using the r.h.s.\ of
\eq{as8eq22} in place of $\cP,\cQ,\cR(\cdots)$. Fourthly, we choose
an algebra isomorphism $\mu:\K^n\ra D$. The number of such
isomorphisms is $n!$, so to compensate we divide by $n!$, which
together with the factor $(-1)^n(n-1)!$ above yields $(-1)^n/n$.
Combining these simplifications yields:
\ea
&\Pi^\vi_1(\bdss^\al(\tau))=\sum_{n\ge 1}\frac{(-1)^n}{n}\,\cdot
\label{as8eq23}
\\[-3pt]
&\raisebox{-11pt}{\begin{Large}$\displaystyle\biggl[$\end{Large}}\,
\sum_{l\in L}\,\,
\sum_{\substack{\text{injective algebra}\\ \text{morphisms}\\
\text{$\mu:\K^n\ra C_l$}}}
\begin{gathered}[t]
\big\vert N_{A_l^\t}(C_l^\t)/C_{A_l^\t}\bigl(\mu((\K^\t)^n\bigr)
\cap N_{A_l^\t}(C_l^\t)\big\vert^{-1}\cdot \\[-3pt]
\bigl[\bigl([U_l^{\mu((\K^\t)^n)}/C_{A_l^\t}(\mu(\K^\t)^n)],\rho_l\ci
\io^{\mu((\K^\t)^n))}\bigr)\bigr]
\end{gathered}
\,\raisebox{-11pt}{\begin{Large}$\displaystyle\biggr]$\end{Large}}.
\nonumber
\ea

Let $n,l,\mu$ be as in \eq{as8eq23} and $u\in U_l^{\mu((\K^\t)^n)}$
project to $[X]\in\fF_l(\K)\subseteq\fObj_\A(\K)$. Then the morphism
$\mu:(\K^\t)^n\ra\Stab_{A_l^\t}(u)\cong\Aut(X)$ induces a splitting
$X\cong X_1\op\cdots\op X_n$, with
$\mu(\ga_1,\ldots,\ga_n)\cong\ga_1\id_{X_1}+\cdots+\ga_n\id_{X_n}$,
and $X_i\not\cong 0$. Conversely, one can show that any
$[X]\in\fF_l(\K)$ and splitting $X\cong X_1\op\cdots\op X_n$ with
$X_i\not\cong 0$ come from such $\mu,u$, and the possible choices of
$\mu$ are all conjugate under the Weyl group $W_{A_l^\t}$ of
$A_l^\t$, and having chosen $\mu$ the possible choices of $u$ form a
$C_{A_l^\t}(\mu(\K^\t)^n)$-orbit in $U_l(\K)$. The orbit of $\mu$
under $W_{A_l^\t}$ is finite and isomorphic to~$N_{A_l^\t}(C_l^\t)/
C_{A_l^\t}\bigl(\mu((\K^\t)^n\bigr)\cap N_{A_l^\t}(C_l^\t)$.

Now a splitting $X\cong X_1\op\cdots\op X_n$ is equivalent to a
$(\{1,\ldots,n\},\bu)$-configuration $(\si,\io,\pi)$ with
$\si(\{1,\ldots,n\})=X$ and $\si(\{i\})=X_i$, up to canonical
isomorphism. Thus we see that the bottom line $[\cdots]$ of
\eq{as8eq23} is equal as a stack function to
$[(\fG^\al_n,\bs\si(\{1,\ldots,n\}))]$, where $\fG^\al_n$ is the
open $\K$-substack of points $[(\si,\io,\pi)]$ in
$\fM(\{1,\ldots,n\},\bu)_\A$ with
$[\si(\{1,\ldots,n\})]\in\Oss^\al(\tau)$ and $\si(\{i\})\not\cong 0$
for all $i=1,\ldots,n$. The factor $\vert N_{A_l^\t}(C_l^\t)
/C_{A_l^\t}(\mu((\K^\t)^n)\cap N_{A_l^\t}(C_l^\t)\vert^{-1}$ exactly
cancels the multiplicity of choices of $\mu$ to make this true.

Let $[(\si,\io,\pi)]\in\M(\{1,\ldots,n\},\bu)_\A$ with
$\si(\{i\})\not\cong 0$ for all $i=1,\ldots,n$. Define
$\ka:\{1,\ldots,n\}\ra C(\A)$ by $\ka(i)=[\si(\{i\})]$, so that
$(\si,\io,\pi)$ is an $(\{1,\ldots,n\},\bu,\ka)$-configuration. As
$\si(\{1,\ldots,n\})\cong\si(\{1\})\op\cdots\op\si(\{n\})$, it is
easy to show that $[\si(\{1,\ldots,n\})]\!\in\!\Oss^\al(\tau)$ if
and only if $\ka(\{1,\ldots,n\})\!=\!\al$, $\tau\!\ci\!\ka\!\equiv
\!\tau(\al)$ and $\si(\{i\})$ is $\tau$-semistable for all $i$, that
is, $[(\si,\io,\pi)]\!\in\!\Mss(\{1,\ldots,n\},\bu,\ka,\tau)_\A$.
Thus
\begin{equation*}
\fG^\al_n(\K)=\coprod_{\begin{subarray}{l}
\text{$\ka:\{1,\ldots,n\}\ra C(\A):$}\\
\text{$\ka(\{1,\ldots,n\})\!=\!\al$, $\tau\!\ci\!\ka\!\equiv
\!\tau(\al)$}\end{subarray}}\Mss(\{1,\ldots,n\},\bu,\ka,\tau)_\A.
\end{equation*}
Hence $[(\fG^\al_n,\bs\si(\{1,\ldots,n\}))]$ equals the second sum
in \eq{as8eq4}. But it also equals the bottom line of \eq{as8eq23},
so comparing \eq{as8eq4}, \eq{as8eq23} completes the proof.
\end{proof}

The theorem enables us to interpret the stack functions
$\bdsi^\al(\tau),\bdsi(I,\tl,\ka,\tau)$. Since $\bdss^\al(\tau)$ is
the `characteristic function' of $\Oss^\al(\tau)$ and $\Pi^\vi_1$ is
the projection to stack functions `supported on virtual
indecomposables', we should understand $\bdsi^\al(\tau)$ as the
`characteristic function of $\tau$-semistable virtual
indecomposables in class $\al$', and $\bdsi(I,\tl,\ka,\tau)$ as the
`characteristic function of $(I,\tl,\ka)$-configurations
$[(\si,\io,\pi)]$ with each $\si(\{i\})$ $\tau$-semistable and
virtual indecomposable'. Note that because `virtual indecomposable'
stack functions can have nonzero components over decomposable
objects, $\bdsi^\al(\tau), \bdsi(I,\tl,\ka,\tau)$ will generally not
be supported on~$\Osi^\al(\tau),\Msi(I,\tl,\ka,\tau)_\A$.

It remains to interpret $\bdst^\al(\tau)$ and $\bdst,\bdssb,\bdsib,
\bdstb(I,\tl,\ka,\tau)$. These are all defined by analogues of
constructible functions equations in \S\ref{as5}--\S\ref{as6} that
were proved using $\chi(\K^m)=1$. Since the spaces $\SF(\cdots)$ do
not set $[\K^m]=1$, the $\bdst^\al(\tau)$ and $\bdst,\ldots,\bdstb
(I,\tl,\ka,\tau)$ do {\it not\/} have a nice interpretation
in~$\SF(\cdots)$.

However, in the spaces $\oSF(\fF,\Th,\Om)$ the relations do set
$[\K^m]=1$, so here the identities have the same interpretations as
their constructible function analogues, but using ideas of `virtual
$\tau$-stable' and `virtual best configuration'. Thus, we interpret
$\bdst^\al(\tau)$ in $\oSFa(\fObj_\A,\Th,\Om)$ as the
`characteristic function of virtual $\tau$-stables in class $\al$',
and $\bdssb(I,\tl,\ka,\tau)$ in $\oSF(\fM(I,\tl,\ka)_\A,\Th,\Om)$ as
the `characteristic function of virtual best
$(I,\tl,\ka)$-configurations $[(\si,\io,\pi)]$ with each
$\si(\{i\})$ $\tau$-semistable', and so on.

This suggests that if we wish to define invariants `counting
$\tau$-stables in class $\al$' we should apply some linear map to
$\bdst^\al(\tau)$ in $\oSFa(\fObj_\A,\Th,\Om)$, but we should not
work in larger spaces such as $\SF(\fObj_\A)$, as the result might
not mean what we want it to mean. The same applies to `counting best
configurations'.

Combining Theorem \ref{as8thm3} with the ideas of \S\ref{as74} we
prove:

\begin{thm} In Definition \ref{as8def1}, for all\/ $k\ge 0$ we have
\begin{gather}
\begin{aligned}
\Pi^\vi_k&\bigl(\bs\si(I)_*\,\bdsib(I,\tl,\ka,\tau)\bigr)\\
=&\begin{cases} \bs\si(I)_*\,\bdsib(I,\tl,\ka,\tau), &
\text{$(I,\tl)$ has $k$ connected components,}\\ 0, &
\text{otherwise,} \end{cases}
\end{aligned}
\label{as8eq24}
\\
\begin{aligned}
\Pi^\vi_k&\bigl(\bs\si(I)_*\,\bdstb(I,\tl,\ka,\tau)\bigr)\\
=&\begin{cases} \bs\si(I)_*\,\bdstb(I,\tl,\ka,\tau), &
\text{$(I,\tl)$ has $k$ connected components,}\\
0, & \text{otherwise,} \end{cases}
\end{aligned}
\label{as8eq25}
\end{gather}
supposing $(\tau,T,\le)$ is a stability condition in \eq{as8eq25}.
Also~$\bep^\al(\tau)\!\in\!\SFai(\fObj_\A)$.
\label{as8thm4}
\end{thm}

\begin{proof} Make the convention that the constants $C_{\ldots},
D_{\ldots},E_{\ldots},F_{\ldots}$ below lie in $\Q$ and depend only
on their subscripts up to isomorphism. In \cite[Th.~5.16]{Joyc4}, if
$(I,\tl)$ is a finite poset and $f_i\in\SFai(\fObj_\A)$ for $i\in
I$, we write $\Pi^\vi_k(P_\sIt(f_i:i\in I))$ as a $\Q$-linear
combination of $P_\sIp(f_i:i\in I)$ over partial orders $\pr$ on $I$
dominated by $\tl$. Since $\bdsi(I,\tl,\ka,\tau)=P_\sIt
(\bdsi^{\ka(i)}(\tau):i\in I)$ and $\bdsi^{\ka(i)}(\tau)\in
\SFai(\fObj_\A)$ by Theorem \ref{as8thm3}, this implies a universal
formula
\ea
\sum_{\substack{\text{p.o.s $\pr$ on $I:$}\\ \text{$\tl$ dominates
$\pr$}}}C_{\tl,\pr,k}\cdot \bs\si(I)_*\,\bdsi(I,\pr,\ka,\tau)&=
\Pi^\vi_k\bigl(\bs\si(I)_*\,\bdsi(I,\tl,\ka,\tau)\bigr). \nonumber
\allowdisplaybreaks\\
\intertext{Combining this with \eq{as8eq5} and \eq{as8eq11} gives}
\sum_{\substack{\text{p.o.s $\pr$ on $I:$}\\ \text{$\tl$ dominates
$\pr$}}}D_{\tl,\pr,k}\cdot \bs\si(I)_*\,\bdsib(I,\pr,\ka,\tau)&=
\Pi^\vi_k\bigl(\bs\si(I)_*\,\bdsib(I,\tl,\ka,\tau)\bigr).
\label{as8eq26}
\ea

In \cite[Th.~5.17]{Joyc4} we show that if $f\in\SFa(\fObj_\A)$ with
$\Pi^\vi_1(f)=f$ and $\cha\K=0$ then $\pi^\stk_{\fObj_\A}(f)$ is
supported on points $[X]$ for $0\not\cong X$ indecomposable. A
generalization of the same proof shows that if $\Pi^\vi_k(f)=f$ then
$\pi^\stk_{\fObj_\A}(f)$ is supported on points $[X_1\op\cdots\op
X_k]$ for $0\not\cong X_a$ indecomposable. Since
$(\Pi^\vi_k)^2=\Pi^\vi_k$, we see that for any $f\in\SFa(\fObj_\A)$,
$\pi^\stk_{\fObj_\A}\bigl(\Pi^\vi_k(f)\bigr)$ is the component of
$\pi^\stk_{\fObj_\A}(f)$ supported on points $[X_1\op\cdots\op X_k]$
for $0\not\cong X_a$ indecomposable.

Applying this to $f=\bs\si(I)_*\,\bdsib(I,\tl,\ka,\tau)$, so that
$\pi^\stk_{\fObj_\A}(f)=\CF^\stk(\bs\si(I))\ab
\dsib(I,\tl,\ka,\tau)$ by Corollary \ref{as8cor1}, and using
Proposition \ref{as7prop2} and \eq{as8eq26} shows that
\e
\begin{split}
&\ts\sum_{\,\text{p.o.s $\pr$ on $I$: $\tl$ dominates $\pr$}\,}
D_{\tl,\pr,k}\cdot\CF^\stk(\bs\si(I))\dsib(I,\pr,\ka,\tau)\\
&\quad=\pi^\stk_{\fObj_\A}\bigl[\Pi^\vi_k\bigl(\bs\si(I)_*\,
\bdsib(I,\tl,\ka,\tau)\bigr)\bigr]\\
&\quad=\begin{cases} \CF^\stk(\bs\si(I))\dsib(I,\tl,\ka,\tau), &
\text{$(I,\tl)$ has $k$
connected components,} \\
0, & \text{otherwise.} \end{cases}
\end{split}
\label{as8eq27}
\e
Now the difference between the top and bottom lines of \eq{as8eq27}
is a {\it universal linear relation} on the
$\CF^\stk(\bs\si(I))\dsib(I,\pr,\ka,\tau)$. Theorem \ref{as7thm3}
shows that there exist no such universal linear relations with
nonzero coefficients. Therefore $D_{\tl,\pr,k}$ is 1 when $\tl=\pr$
and $(I,\tl)$ has $k$ connected components, and 0 otherwise.
Equation \eq{as8eq24} now follows from~\eq{as8eq26}.

Next we prove \eq{as8eq25}. Substituting \eq{as5eq8} into
\eq{as5eq14} into \eq{as6eq3} into \eq{as6eq14} and applying
$\CF^\stk(\bs\si(I))$ gives a universal formula
\e
\begin{aligned}
\sum_{\substack{\text{iso. classes of\/ $\A$-data $(J,\ls,\la)$
and}\\
\text{surjective $\phi\!:\!J\!\ra\!I$: $i\ls j$ implies
$\phi(i)\!\tl\!\phi(j)$,}\\
\text{$\la(\phi^{-1}(i))\!=\!\ka(k)$ for $i\!\in\!I$,
$\tau\!\ci\!\ka\!\ci\!\phi\!\equiv\!\tau\!\ci\!\la$}}
\!\!\!\!\!\!\!\!\!\!\!\!\!\!\!\!\!\!\!\!\!\! }
\!\!\!\!\!\!\!\!\!\!\!\!
\begin{aligned}[t]
E_{J,\ls,I,\tl,\phi}\CF^\stk(\bs\si(J))\dsib(J,\ls,\la,\tau)&=\\
\CF^\stk(\bs\si(I))\dstb(I,\tl,&\ka,\tau).
\end{aligned}
\end{aligned}
\label{as8eq28}
\e
Fix $(I,\tl,\ka)$ and $(J,\ls,\mu)$ in \eq{as8eq28}, and let
$(I,\tl),(J,\ls)$ have $k,l$ connected components. Then by
Proposition \ref{as7prop2}, the terms on the right and left hand
sides of \eq{as8eq28} are supported on points $[X_1\op\cdots\op
X_k]$ and $[Y_1\op\cdots\op Y_l]$ in $\fObj_\A(\K)$ respectively,
with all $X_a,Y_b$ indecomposable. So, for fixed $k\ne l$, consider
the sum of all terms on the l.h.s.\ of \eq{as8eq28} in which
$(J,\ls)$ has $l$ connected components. This is simply the component
of \eq{as8eq28} supported on $[Y_1\op\cdots\op Y_l]$ for $Y_b$
indecomposable, and as $k\ne l$ the r.h.s.\ of \eq{as8eq28} is zero
on such points. Thus restricting to $(J,\ls)$ with $l$ connected
components gives a universal identity of the form \eq{as7eq27}.
Theorem \ref{as7thm3} therefore shows that $E_{J,\ls,I,\tl,\phi}=0$
if~$k\ne l$.

Similarly, substituting \eq{as8eq11} into \eq{as8eq13} into
\eq{as8eq3} into \eq{as8eq17} and applying $\bs\si(I)_*$ gives the
stack function analogue of \eq{as8eq28}, with the same
$E_{J,\ls,I,\tl,\phi}$. This writes $\bs\si(I)_*\,\bdstb
(I,\tl,\ka,\tau)$ as a linear combination of $\bs\si(J)_*\,
\bdsib(J,\ls,\la,\tau)$, over $(J,\ls)$ with the same number of
connected components as $(I,\tl)$. But \eq{as8eq24} shows
$\Pi^\vi_k$ is the identity on these terms if this number of
connected components is $k$, and 0 otherwise. Equation \eq{as8eq25}
follows.

Finally, substituting \eq{as5eq8} into \eq{as5eq14} into \eq{as7eq4}
into \eq{as7eq17} gives an identity
\e
\begin{aligned}
\sum_{\substack{\text{iso. classes of\/ $\A$-data $(I,\pr,\ka)$:}\\
\text{$\ka(I)=\al$, $\tau\ci\ka\equiv\tau(\al)$}}}
\!\!\!\!\!\!\!\!\!\!
F_{I,\pr}\CF^\stk(\bs\si(I))\dsib(I,\pr,\ka,\tau)&=\ep^\al(\tau).
\end{aligned}
\label{as8eq29}
\e
Using Proposition \ref{as7prop2} and Theorem \ref{as7thm2}, the same
method shows $F_{I,\pr}=0$ unless $(I,\pr)$ is connected.
Substituting \eq{as8eq11} into \eq{as8eq13} into the analogue of
\eq{as7eq4} into \eq{as8eq9} gives the stack function version of
\eq{as8eq29}, writing $\bep^\al(\tau)$ as a linear combination of
$\bs\si(I)_*\,\bdsib(I,\pr,\ka,\tau)$ with $(I,\pr)$ connected. By
\eq{as8eq24}, $\Pi^\vi_1$ is the identity on each term, so
$\Pi^\vi_1\bigl(\bep^\al(\tau)\bigr)=\bep^\al(\tau)$,
and~$\bep^\al(\tau)\in\SFai(\fObj_\A)$.
\end{proof}

By \eq{as8eq19} and \eq{as8eq24} $\bHp_\tau$ is spanned by
eigenvectors of $\Pi^\vi_k$, proving:

\begin{cor} In Definition \ref{as8def2}, $\bHp_\tau$ is closed under
$\Pi^\vi_k$ for all\/~$k\ge 0$.
\label{as8cor2}
\end{cor}

In general $\bHt_\tau$ is not closed under $\Pi^\vi_k$ for $k>0$. We
can now define and study Lie algebras $\bLp_\tau,\bLt_\tau$, the
analogues of~$\Lp_\tau,\Lt_\tau$.

\begin{dfn} Let Assumption \ref{as3ass} hold, and $(\tau,T,\le)$ be a
permissible weak stability condition on $\A$. Define
$\bLp_\tau=\bHp_\tau\cap\SFai(\fObj_\A)$. Then $\bLp_\tau$ is a {\it
Lie subalgebra} of $\SFai(\fObj_\A)$, since Theorem \ref{as8thm2}
implies $\bHp_\tau$ is a Lie algebra. From \eq{as8eq19},
\eq{as8eq24} and \eq{as8eq25} we see that
\e
\begin{split}
\bLp_\tau\!&=\!\bigl\langle\bs\si(I)_*\,\bdsib(I,\pr,\ka,\tau):
\text{$(I,\pr,\ka)$ $\A$-data, $(I,\pr)$
connected\/}\bigr\rangle{}_\Q\\
&=\!\bigl\langle\bs\si(I)_*\,\bdstb(I,\pr,\ka,\tau):
\text{$(I,\pr,\ka)$ $\A$-data, $(I,\pr)$
connected\/}\bigr\rangle{}_\Q,
\end{split}
\label{as8eq30}
\e
supposing $(\tau,T,\le)$ is a stability condition in the second
line. Using \eq{as7eq16}, \eq{as8eq30}, Corollary \ref{as8cor1} and
\eq{as3eq10} a Lie algebra morphism, we find
$\pi^\stk_{\fObj_\A}:\bLp_\tau\ra\Lp_\tau$ is a surjective Lie
algebra morphism when $\cha\K=0$. Also $\bLp_\tau$ generates
$\bHp_\tau$ as in Proposition \ref{as7prop2}, so there is a natural,
surjective $\Q$-algebra morphism $\bar\Phi^{\rm
pa}_\tau:U(\bLp_\tau)\ra\bHp_\tau$. As we have no analogue of
Proposition \ref{as3prop2} we cannot show $\bar\Phi^{\rm pa}_\tau$
is an isomorphism, but the ideas of \S\ref{as74} imply there is no
nontrivial `universal' kernel of $\bar\Phi^{\rm pa}_\tau$ generated
by universal multiplicative relations on~$\bLp_\tau$.

Motivated by Corollary \ref{as7cor2}, and using Theorem
\ref{as8thm4}, define $\bLt_\tau$ to be the Lie subalgebra of
$\SFai(\fObj_\A)$ generated by the $\bep^\al(\tau)$ for all $\al\in
C(\A)$. Then $\bLt_\tau\subseteq\bLp_\tau$. Using Corollaries
\ref{as7cor2} and \ref{as8cor1} and \eq{as3eq10} a Lie algebra
morphism, we see that $\pi^\stk_{\fObj_\A}:\bLt_\tau\ra\Lt_\tau$ is
a surjective Lie algebra morphism. Equation \eq{as8eq20} implies
$\bLt_\tau$ generates $\bHt_\tau$, so there is a natural, surjective
$\Q$-algebra morphism $\bar\Phi^{\rm to}_\tau:U(\bLt_\tau)\ra
\bHt_\tau$, but as above we cannot prove $\bar\Phi^{\rm to}_\tau$ is
an isomorphism. As we have no stack function analogue of Proposition
\ref{as3prop2}, and $\bHt_\tau$ may not be closed under $\Pi^\vi_k$,
we also cannot prove that~$\bLt_\tau=\bHt_\tau\cap\SFai(\fObj_\A)$.
\label{as8def3}
\end{dfn}

In \cite{Joyc5} we will apply these ideas as follows. Under extra
assumptions on $\A$, in \cite[\S 6]{Joyc4} we defined (Lie) algebra
morphisms $\Phi^{\sst\La}\ci\Pi^{\Up,\La}_{\fObj_\A},\ldots$ from
$\SFa(\fObj_\A)$ or $\SFai(\fObj_\A)$ to some explicit algebras
$A(\A,\La,\chi),\ldots,C(\A,\Om,\chi)$. Restricting these yields
(Lie) algebra morphisms from $\bHp_\tau,\bHt_\tau$
or~$\bLp_\tau,\bLt_\tau$.

We shall regard these maps $\bHp_\tau\ra A(\A,\La,\chi),\ldots$ as
encoding {\it systems of invariants} that `count'
$\tau$-(semi)stable objects and configurations. The fact that the
maps are morphisms implies {\it multiplicative relations} upon these
invariants, and also that the map is determined by its values on a
generating set for the (Lie) algebra, such as the $\bep^\al(\tau)$
for $\bHt_\tau$ or $\bLt_\tau$. The identities of
\S\ref{as5}--\S\ref{as7} imply identities on the invariants, and the
results of \cite{Joyc5} yield {\it transformation laws} for the
invariants between different stability
conditions~$(\tau,T,\le),(\ti\tau,\ti T,\le)$.

In particular, if $P$ is a {\it Calabi--Yau $3$-fold\/} and
$\A=\coh(P)$, then \cite[\S 6.6]{Joyc4} defined a Lie algebra
morphism $\Psi^{\sst\Om}\ci\bar\Pi^{\Th,\Om}_{\fObj_\A}:\SFai
(\fObj_\A)\ra C(\A,\Om,\frac{1}{2}\bar\chi)$. Restricting this to
$\bLp_\tau$ and $\bLt_\tau$ yields interesting invariants `counting'
$\tau$-semistable sheaves on Calabi--Yau 3-folds, with attractive
transformation laws, which may be related to {\it Donaldson--Thomas
invariants}. This is one reward for the work we put in to construct
$\bLp_\tau,\bLt_\tau$ and show they lie in~$\SFai(\fObj_\A)$.

\medskip

\noindent{\small\sc The Mathematical Institute, 24-29 St. Giles,
Oxford, OX1 3LB, U.K.}

\noindent{\small\sc E-mail: \tt joyce@maths.ox.ac.uk}

\end{document}